\documentclass[mnsc,nonblindrev]{informs3}

\OneAndAHalfSpacedXI 

\usepackage{natbib}
 \bibpunct[, ]{(}{)}{,}{a}{}{,}%
 \def\bibfont{\small}%
\usepackage[ruled]{algorithm2e}
\usepackage{booktabs}
\usepackage{enumitem}
\usepackage{multirow}
\newenvironment{behavior-condition}[1]{\par\noindent\underline{Behavioral Conditions:}\space#1}{}
\usepackage{subcaption}
\usepackage[title]{appendix}
\usepackage{hyperref}

\TheoremsNumberedThrough     
\ECRepeatTheorems

\EquationsNumberedThrough    

\MANUSCRIPTNO{}
\usepackage{bm}
\usepackage{makecell}
\usepackage{mdframed}
\usepackage[dvipsnames]{xcolor}
\usepackage{threeparttable}

\usepackage[normalem]{ulem}

\newcommand{\calS}{\mathcal{S}}

\newcommand{\bmS}{{\bm S }}

\newcommand{\dstar}{$ ^{\ast\ast} $}
\renewcommand{\star}{$ ^{\ast} $}
\newenvironment{proofref}[1]{{\bf Proof of {#1}.\/}}{\hfill $\Box$\vskip 0.01in}
\let\originalPr\Pr
\renewcommand{\Pr}[1]{\originalPr\!\left(#1\right)}

\begin{document}



\RUNTITLE{A Behavioral Model for Exploration vs. Exploitation}

\TITLE{A Behavioral Model for Exploration vs. Exploitation: Theoretical Framework and Experimental Evidence}

\ARTICLEAUTHORS{%
\AUTHOR{Jingying Ding}
\AFF{School of Information Management and Engineering,
	Shanghai University of Finance and Economics,  \EMAIL{dingjingying@mail.shufe.edu.cn}} 
\AUTHOR{Yifan Feng}
\AFF{Department of Analytics and Operations,
   National University of Singapore, \EMAIL{yifan.feng@nus.edu.sg}}
\AUTHOR{Ying Rong}
\AFF{Antai College of Economics and Management,
  Shanghai Jiao Tong University,  \EMAIL{yrong@sjtu.edu.cn}}
} 

\ABSTRACT{
	How do people navigate the exploration-exploitation (EE) trade-off when making repeated choices with unknown rewards? We study this question through the lens of multi-armed bandit problems and introduce a novel behavioral model, \textit{Quantal Choice with Adaptive Reduction of Exploration} (QCARE). It generalizes Thompson Sampling, allowing for a principled way to quantify the EE trade-off and reflect human decision-making patterns. The model adaptively reduces exploration as information accumulates, with the reduction rate serving as a parameter to quantify the EE trade-off dynamics. We theoretically analyze how varying reduction rates influence decision quality, shedding light on the distinct effects of ``over-exploration'' and ``under-exploration.'' Empirically, we validate QCARE using behavioral data collected from human participants. QCARE not only captures critical behavioral patterns in the EE trade-off but also performs competitively with, and often outperforms, alternative models in out-of-sample prediction. Our analysis reveals a behavioral tendency toward over-exploration.
}


\HISTORY{This version: \today.}

\maketitle

%

\section{Introduction}

The exploration-exploitation (EE) trade-off is a fundamental challenge in sequential decision making under uncertainty. Consider an online platform allocating a limited stream of user traffic among products, sellers, advertisements, or content whose performance is initially unknown. Routing traffic to alternatives that have performed well so far exploits current evidence and generates immediate value, whereas routing traffic to less-tested alternatives produces information that may reveal better opportunities. Too much exploration sacrifices current performance; too little may leave promising alternatives undiscovered. Effective sequential decision making therefore requires balancing the value of learning about uncertain alternatives against the value of acting on what has already been learned.

The stochastic multi-armed bandit (MAB) problem provides a canonical formalization of this tension. Each alternative is represented as an arm with an unknown but fixed reward distribution. In each period, the decision maker selects one arm and observes only its realized reward. Because each choice generates both a payoff and information, maximizing cumulative reward requires learning about uncertain arms without sacrificing too much current value. Policies such as the Gittins index, Upper Confidence Bound, and Thompson Sampling operationalize this balance through dynamic indices, optimism bonuses, and posterior sampling, respectively. Together, these policies provide normative prescriptions for how a reward-seeking algorithm should balance exploration and exploitation under specified objectives and assumptions.

Yet not all important EE decisions are delegated to algorithms. Similar selective-feedback problems arise when individuals or teams make sequential choices themselves. For example, a clinical team allocating comparable patients among established therapies with uncertain response rates observes the outcome of the treatment actually administered, but not the counterfactual outcomes under alternative therapies. Similarly, a consumer repeatedly choosing among unfamiliar experience goods learns about the quality of the purchased alternative but does not observe the experience that unchosen alternatives would have provided. Related trade-offs also arise in organizational resource allocation in contexts such as experimental trials (\citealp{zhuo2025navigating}). These examples illustrate that the EE trade-off is not only an algorithm-design problem but also a behavioral and organizational decision problem.

Human decision makers, however, need not manage the EE trade-off in the same manner as algorithms designed to maximize cumulative reward. Whereas such algorithms follow prespecified updating and choice rules, human decisions may also reflect behavioral heuristics, heterogeneous preferences, and so on. Relatively simple behavioral rules can sometimes fit experimental bandit choices better than Gittins-index-based policies \citep{GansKnoxCroson2007}. Yet this evidence leaves open how to parsimoniously characterize changes in human decision makers' propensity to sample currently non-leading alternatives as arm-specific experience accumulates and how those changes affect decision quality. Motivated by this gap, we study the following questions: 

\begin{itemize}
\item \textit{Is there a parsimonious family of bandit policies that spans the EE spectrum with a single interpretable parameter, so that ``how much exploration'' becomes a measurable quantity, with a rigorous theoretical characterization of the learning performance?}  
\item \textit{Can the same family be used as an empirical model, so that given observed behavioral data, one can use the calibrated model to reveal new insights regarding how human beings balance the EE trade-off?}
\end{itemize}

\subsection{Roadmap and Summary of Contributions}\label{subsec:summary-contributions}

\noindent \textbf{Roadmap.} We answer both questions by introducing \textit{Quantal Choice with Adaptive Reduction of Exploration} (QCARE), a family of bandit policies indexed by a single parameter $\alpha$, which we call the \textit{reduction rate of exploration}. Our analysis proceeds in three steps, which also organize the contributions summarized below. First, guided by model-free evidence from incentivized experiments (Section~\ref{subsec:behavioral-data-descriptive-analysis}), we develop QCARE as a behavioral model of repeated choice under the EE trade-off (Section~\ref{sec:QCARE-Intro}). Second, we characterize how $\alpha$ shapes learning and regret, which identifies the reward-maximizing reduction rate and gives over- and under-exploration a precise meaning (Section~\ref{sec:QCARE-Theory}). Third, we estimate $\alpha$ for each participant, evaluate QCARE's predictive performance against alternative models, and investigate the potential behavioral patterns (Section~\ref{sec:structure-estimation-behavioral-data}). 

\vspace{0.2 cm}
\noindent \textbf{Modeling contributions.}  The model we introduce leverages concepts from the quantal choice framework, where the decision maker relies on a (randomized) score system and chooses the arm with the highest attraction score. At a high level, the score is arm- and history-dependent and takes the following form:
\begin{align*}
\textsc{Score} \ =\  \textsc{Historical Performance} \ +\ \textsc{Weight}  \ \times  \ \textsc{Random Shock}.
\end{align*}

We refer to Equation~\eqref{eq:QCARE-definition} for the precise formula. In short, the historical-performance term favors arms with higher average rewards and thus captures exploitation, whereas the random shock allows arms with weaker records to be chosen and thus captures exploration. Finally, the weight on this component determines how much exploration remains after an arm has accumulated experience. An important feature of our model is that the weight shrinks with the number of times the arm has been pulled. This feature is motivated by our experimental evidence and gives QCARE its name. The larger $\alpha$ is, the faster exploration is reduced.  The specific form of the weight is set up so that QCARE nests Gaussian Thompson Sampling as the benchmark case when \(\alpha=0.5\); see Proposition~\ref{prop:TS_QCARE_equa}. Correspondingly, $ \alpha $ values smaller (resp. larger) than this benchmark describe policies that explore more (resp. less) aggressively than Thompson Sampling.

By connecting Thompson Sampling with quantal choice models, QCARE combines three desirable features. First, it captures the learning effect parsimoniously, with a single interpretable parameter. The traditional approach to dynamic choice modeling derives choice probabilities from dynamic programming \citep[e.g.,][]{rust1987optimal,erdem1996decision} and therefore suffers from the curse of dimensionality. In contrast, QCARE's scores are \textit{myopic}, which keeps both the theoretical analysis and the empirical estimation relatively simple.

Second, it captures a spectrum of exploration-exploitation behavior, which allows us to capture choice behavior with potentially bounded rationality. It also does so in a \textit{principled} way, as our analysis characterizes how $\alpha$ shapes learning and regret across this spectrum. The notions of over- and under-exploration, which play central roles in the paper, thus emerge from this analysis.
	
Third, it is a randomized policy that captures the learning effect, which makes it empirically friendly to estimate from individual choice data. Our empirical validation shows that QCARE predicts human choices well out of sample.

\vspace{.2 cm}
\noindent \textbf{Technical contributions.} We first establish basic properties of QCARE. In Proposition \ref{prop:QCARE-basic-properties}, we show that QCARE displays comparative statics properties that are qualitatively consistent with our lab evidence. In Theorem \ref{thm:UB_sublinear}, we show that all $ \alpha  > 0$ enable the decision maker to converge to the optimal arm in the long run. That distinguishes QCARE from behavioral policies with persistent exploration, such as the fixed-\(\varepsilon\)-greedy benchmark considered
in our empirical analysis. Furthermore, it suggests that every $ \alpha > 0 $ is ``plausible'' -- at least when $ T $ is large -- since all such values correspond to long-run-average optimal policies.

Our main theoretical results characterize how $\alpha$ affects both instance-dependent and instance-independent (i.e., worst-case) regret order over the entire range $\alpha>0$, which yields a comprehensive regret landscape summarized in Table~\ref{tab:major-theoretical-results}. At $\alpha=0.5$, QCARE inherits the standard regret of Gaussian Thompson Sampling. When $\alpha<0.5$, the regret is of order $(\log T)^{1/(2\alpha)}$ on every fixed instance, and the instance-independent regret is of order $T^{1-\alpha}$ up to a logarithmic factor (Theorem~\ref{thm:smaller0.5}). When $\alpha>0.5$, its regret becomes nearly linear, $\Omega(T^{1-o(1)})$, on every interior Bernoulli instance (Theorem~\ref{thm:lb_greater0.5}). Together, these results identify $\alpha=0.5$ as the sweet spot of the EE trade-off asymptotically. Otherwise, regret deteriorates gradually as $\alpha$ falls below $0.5$ (due to over-exploration) and sharply as $\alpha$ rises above it (due to under-exploration).

	\begin{table}[htbp]
		\centering
		\caption{Summary of Regret Order Bounds for QCARE}
		\label{tab:major-theoretical-results}
			\begin{threeparttable}
				\small
				\setlength{\tabcolsep}{3pt}
				\renewcommand{\arraystretch}{1.5}
				\begin{tabular*}{\linewidth}{%
						@{}l l@{\extracolsep{\fill}}c c c@{}
					}
					\toprule
					& &
					\makebox[0.24\linewidth]{$0<\alpha<0.5\,{}^{\ast}$}
					& \makebox[0.13\linewidth]{%
						\footnotesize $\alpha=0.5\,{}^{\dagger}$}
					& \makebox[0.24\linewidth]{$\alpha>0.5$} \\
					\midrule
					
					\multirow{2}{*}{\shortstack[l]{%
							Instance-dependent\\$\mathcal{R}^{\Theta}(T)$}}
					& Upper
					& $O\big((\log T)^{1/(2\alpha)}\big)$
					& {\footnotesize $O(\log T)$}
					& $o(T)^{\ddagger}$ \\
					& Lower
					& $\Omega\big((\log T)^{1/(2\alpha)}\big)$
					& {\footnotesize $\Omega(\log T)$}
					& $\Omega\big(T^{1-o(1)}\big)^{\S}$ \\
					\midrule
					
					\multirow{2}{*}{\shortstack[l]{%
							Instance-independent\\$\mathcal{R}(T)$}}
					& Upper
					& $\widetilde{O}\big(T^{1-\alpha}\big)$
					& {\footnotesize $O(\sqrt{T})$}
					& $O(T)$ \\
					& Lower
					& $\Omega\big(T^{1-\alpha}\big)$
					& {\footnotesize $\Omega(\sqrt{T})$}
					& $\Omega\big(T^{1-o(1)}\big)^{\S}$ \\
					\bottomrule
				\end{tabular*}
				
				\begin{tablenotes}[flushleft]
					\item[] {\scriptsize\raggedright
						\textit{Note:}
						$^{\ast}$~Theorem~\ref{thm:smaller0.5}.
						$^{\dagger}$~Proposition~\ref{prop:TS_QCARE_equa} plus the regret bounds from \citet{lai1985asymptotically} and \citet{agrawal2017near}.
						$^{\ddagger}$~Theorem~\ref{thm:UB_sublinear}.
						$^{\S}$~Theorem~\ref{thm:lb_greater0.5}.
						\par}
				\end{tablenotes}
			\end{threeparttable}%
	\end{table}

Beyond the asymptotic theory, we also conduct extensive numerical studies to examine whether these asymptotic insights carry over to finite horizons; see Section~\ref{subsec:finite-analysis}. In short, the effects of over- and under- exploration hurt in the same way, and hence a reward-maximizing $\alpha$ still exists, but it is generally larger than $0.5$ and depends on the configuration, since there is less time to explore.

\vspace{.2 cm}
\noindent \textbf{Experimental and empirical contributions.} We conducted experiments to collect behavioral data regarding how real human beings make decisions in the MAB problem. We use the data to empirically validate our model. In particular, we compare QCARE with other models drawn from different streams of literature. QCARE displays a strong capability to capture human behavior, especially in terms of out-of-sample prediction power. We also illustrate how the QCARE-generated reward distributions align well with the real data after controlling for $ \alpha $.

Comparing the estimates with the reward-maximizing benchmark reveals a robust pattern: participants tend to \textit{over-explore}. That is, their propensity to choose a non-leading arm declines with experience more slowly than under the reward-maximizing level.
While such bias can be partly rationalized by the risk aversion of the participants, their behavior of not choosing the leading arm appears to reflect a mix of random behavioral errors (in the same spirit as \citealp{Su2008}) layered on top of the ``intrinsic'' exploration required for dynamic learning.

\section{Literature Review}

Despite the salient role of the EE trade-off in decision quality, investigation into decision makers' behavior under the MAB framework is limited. The underlying difficulty is the \textit{dual} challenge of (i) the intricate EE trade-off in the MAB problem itself (for the decision maker) and (ii) the need for machinery to ``learn the learning'' (for the observer). Accordingly, our work bridges MAB policy design and analysis, which addresses the decision maker's problem, with behavioral science and operations management, which provide the tools for learning how people learn.

\vspace{0.2 cm}
\noindent \textbf{MAB policy design and analysis.}
The MAB problem is a classical model that inherently captures the exploration-exploitation trade-off. Various algorithms for solving stochastic MAB problems have been proposed, including the Gittins index~\citep{gittins1979bandit,gittins2011multi}, Upper Confidence Bound (UCB)~\citep{lai1985asymptotically, auer2002finite}, Thompson Sampling~\citep{thompson1933likelihood}, among others. The main goal of those algorithms is to maximize the expected reward (or minimize the regret) in different senses. A particularly relevant algorithm to our paper is Thompson Sampling. It is asymptotically order optimal in terms of the instance-independent regret \citep{agrawal2013further, agrawal2017near} as $ T $ grows, as well as a few other regret metrics such as Bayesian regret \citep{russo2014learning} and instance-dependent regret (\citealp{kaufmann2012thompson}). There is also an emerging literature on a deeper understanding of the system dynamics and reward distribution of those policies, such as diffusion-scaling behavior (\citealp{kalvit2021closer,fan2021diffusion,kuang2024weak}) and tail risk (\citealp{simchi2023regret,fan2021fragility}), to name a few. 

However, (asymptotically) optimal policies may not represent people's behavior well. Motivated by the need to capture people's highly heterogeneous behavior, we develop a class of MAB policies (which may or may not be asymptotically optimal) that span the wide EE trade-off spectrum, and analyze their system dynamics and decision qualities. As such, we believe our analysis makes a theoretical contribution to the analysis of MAB policies in its own right.

\vspace{0.2 cm}
\noindent \textbf{(Dynamic) quantal choice modeling and behavioral science.} As mentioned previously, this paper involves developing MAB policy families to understand human behavior in the decision-making process. That essentially boils down to developing a dynamic discrete choice model with the learning effect. A natural choice model framework is the quantal choice theory. While the literature on static quantal choice models is very rich (e.g., \citealp{mcfadden1976quantal} for models with full rationality, and \citealp{mckelvey1995quantal,Su2008,ChenSuZhao2012, LiChenRong2020} for models with bounded rationality), there is much less work on dynamic models (e.g., \citealp{rust1987optimal} for a representative rational model and \citealp{mckelvey1998quantal} for a behavioral one). 

Dynamic choice models that capture learning effects are even more limited. Among those that do, one approach is to incorporate learning into Rust-type models, i.e., the agent is modeled as a forward-looking rational decision maker who solves a Bayesian dynamic program (\citealp{erdem1996decision,ching2013learning}).
This approach has good theoretical foundations, but it has a few major limitations that motivate us to consider other directions. First, it suffers from the curse of dimensionality and can be computationally prohibitive.\footnote{To put things into perspective, solving for the Bayes optimal policy for the MAB problem is itself challenging (see \citealp{min2019thompson}). Empirical estimation from the data would be even harder. (By analogy, the structural estimation problem in \citealp{rust1987optimal} is much more complicated than the bus engine replacement problem itself.)} Hence, not surprisingly, the \textit{actual} implementation needs various numerical heuristics such as simulation and value function interpolation (\citealp{erdem1996decision}). Second, those models assume that the observer/researcher only observes the decision maker's actions. In our experimental environment, the realized arm rewards are also observed, which is valuable information that should be taken advantage of. Perhaps a more fundamental limitation of this approach is that the decision maker is assumed to be rational (utility-maximizing). It is well documented in the literature that incorporating broader behavioral patterns can better capture real data  (e.g., \citealp{banks1997experimental,GansKnoxCroson2007}). Many researchers have also used experiments to test the qualitative implications of Bayes optimal policies (e.g., Gittins index policy) and found various deviations (\citealp{horowitz1975experimental,meyer1995sequential,anderson2001behavioral}).

That naturally leads us to consider the behavioral approach. There are well-celebrated behavioral learning models in the literature, such as Experience-Weighted Attraction (EWA) (\citealp{camerer1999experience}). However, EWA is based on \textit{full-information feedback} (i.e., the decision maker observes the realized rewards of all arms) while MAB requires \textit{bandit feedback} (i.e., the decision maker only observes the realized rewards of the pulled arm).
Stylized behavioral policies have been considered in various streams of literature, such as behavioral operations management (\citealp{GansKnoxCroson2007}), marketing (\citealp{FerecatuBruyn2022}), cognitive psychology (\citealp{bouneffouf2017bandit}), and computer science (\citealp{zhang2013forgetful}).
Usually, those policies are heuristic and lack formal theoretical understanding, such as how to measure the EE trade-off and how that connects to decision qualities. 
Within this stream, our main innovation lies in applying the toolbox and concepts in online learning to behavioral modeling. To our knowledge, this direction is quite new. The only exception is the work of \citet{mauersberger2022thompson}, who use Thompson Sampling as a behavioral model and observe promising empirical evidence.

\section{The Multi-Armed Bandit Problem and Experiments}

\textbf{Formulation.} We first describe the standard definition of the multi-armed bandit problem. In this problem, the decision maker (DM) faces $ N $ choices, which are referred to as \textit{arms}. Each arm $ i $ corresponds to an i.i.d. sequence of stochastic rewards, denoted by $ \{r_i(t)\}_t $, which is generated from a fixed reward distribution. We refer to the ground truth reward distributions of the arms as an \textit{instance} $ \Theta $ of the MAB problem, which is unknown to the DM. We assume that the rewards are supported on $ [0,1] $. We assume that the reward sequences are mutually independent across arms. We use  $\mu_i$ to denote the expected reward of arm $i$, and let $ i^\ast := \argmax_i \mu_i $ be the optimal arm with the highest mean. For technical simplicity, we focus on cases where the optimal arm is unique. We use $\mu^*:=\max_i \mu_i$ to denote the expected reward of the optimal arm, and let $\Delta_{i}:=\mu^*-\mu_i$  represent the reward gap between arm $ i $ and the optimal arm.

At each time $t \in  [T] = \{1, \ldots, T\}$, the DM chooses one arm to pull, which we refer to as the DM's \textit{action} and denote by $ a(t) $. After that, the DM receives the reward of the pulled arm, $r_{a(t)}(t)$. As such, an admissible \textit{policy} of the DM is governed by a mapping from a historical path $ h_{t-1} := \big(a(1), r_{a(1)}(1), \ldots, a(t-1), r_{a(t-1)}(t-1)\big) $ to a (possibly randomized) action $ a(t) $. The DM's goal is to select the arms to maximize the total rewards over the time horizon $T$. We follow the standard convention in the online learning literature and adopt \textit{regret} to measure decision quality in terms of collecting rewards. One could understand this concept as a way to properly normalize the expected rewards of a policy so that it is easier to measure its performance across instances. The regret of a given policy is defined as the difference in the total expected rewards between always pulling the optimal arm $i^*$ and following the policy. Specifically, for every MAB instance $ \Theta $, the (instance-specific) regret $\mathcal{R}^\Theta (T)$ is defined by 
\begin{align*}
\mathcal{R}^\Theta(T)=\mu^{*} T - \mathbb{E}\left[\sum_{t=1}^{T}r_{a(t)}(t)\right].
\end{align*}

We also define the instance-independent (worst-case) regret as $\mathcal{R}(T):= \max_{\Theta} \mathcal{R}^\Theta (T).$ An effective policy that achieves higher expected reward corresponds to a smaller regret.

\subsection{The MAB Experiments }\label{subsec:experiment-design}

To capture how people \textit{actually} balance the EE trade-off in the MAB problem, we collected real human behavioral data based on a series of experiments. We first explain the details of our experimental design and then present descriptive findings from the experimental data that further motivate our quantitative behavioral model. 

\vskip 0.2 cm
\noindent \textbf{Setup.} Our main experiment was conducted in a controlled laboratory setting and is also referred to as the \textit{baseline} experiment throughout the paper. Following common practice in both behavioral economics and behavioral operations literature \citep{bolton2012managers}, we recruited students from a university behavioral economics laboratory as proxies for human behavior.
We investigated two-arm MAB problems with Bernoulli rewards and examined two MAB instances, namely $ (\mu_1, \mu_2) = (0.4, 0.5) $ and $ (\mu_1, \mu_2) = (0.6, 0.4) $. We also considered two time horizons, $T=100$ and $T=200$. Thus, the baseline experiment included four distinct configurations.

\vspace{0.2 cm}
\noindent \textbf{Procedure, information, and incentives.} We conducted our experiments using the oTree system \citep{chen2016otree}. Once recruited, each subject was first presented with an overview of the game. In this overview, the subject learned basic information about the game, such as the rules, the number of arms, and the fact that rewards followed fixed binary distributions. However, no other information about the rewards, such as the mean reward of either arm, was provided. In this way, the design preserved the information structure of a standard MAB problem. The subject was also informed of the goal of the MAB game, namely, to accumulate as many points as possible, and how much money they could earn for each additional point. Each subject was assigned to only one of the four configurations and played the MAB game at most once to avoid spillover or cross-learning effects. We also employed attention checks to ensure response quality.

Once the MAB game started, the participant was presented with a complete history of actions and realized rewards of each period, so that the participant could implement any admissible MAB policy at least in principle. Summary statistics of the history were also provided separately, including (i) the action and realized reward of the previous period, (ii) the cumulative reward, and (iii) the number of plays and average points of each arm. Finally, the subject was reminded of the period they were in and of the total number of periods. Based on that information, the participant was asked to choose an arm, and then the MAB game moved on to the next period. Once the MAB game finished, the subjects' payments were settled. They earned one point each time the selected arm resulted in a realized reward of one. The total earnings consisted of two components: a show-up fee and a bonus proportional to the points earned. In some experimental conditions, the point-to-cash conversion rate differed across reward environments. This design choice was used to equalize the expected earnings as if the optimal arm had been known, thereby making participant incentives more comparable across conditions. 
We refer the reader to Appendix \ref{section:appendix:baseline experiment} for more details about the procedure, such as the interface and payment calculations, among others.

\vspace{0.2 cm}
\noindent \textbf{Replications and auxiliary experiments.} Beyond the baseline experiment, we also conducted two replications of the same experiment to confirm its core findings.  First, to account for potential limitations in the subject sampling, we conducted an \textit{online replication} with a more diverse participant population, including nationality, gender, and occupation, among others. The other configurations were largely the same except for a slightly different point-to-cash conversion rate. The online replication was facilitated by the Prolific platform, which is recognized for acquiring high-quality data in online human-subject research (\citealp{douglas2023data,albert2023comparing}). Second, since the previous results were surprisingly consistent, we moved on to the \textit{preregistered replication}, which replicated the online replication at a larger scale so as to credibly confirm the robustness of our findings.\footnote{Preregistration link: \url{https://osf.io/ukew8/?view_only=37cca7563fa44da7bbeaac7124459684}. The preregistration was hosted on the Open Science Framework (OSF). We chose it because it is reportedly more robust to activities such as ``preregistration hacking'' than other platforms (\citealp{haroz2022comparison}).}

Furthermore, we conducted three sets of auxiliary experiments for robustness checks in different directions. In the \textit{three-arm extension}, we conducted the experiment in a three-arm environment. In the \textit{no-summary-table extension}, we partially hid the reward summary statistics from the participants. Finally, we conduct a \textit{targeted robustness check} to directly examine the participants' behavioral differences between a control group and three treatment conditions involving practice rounds, higher stakes, and a remaining-period reminder, which provide exogenous shocks to the participants' familiarity with the task, monetary incentives, and the awareness of the remaining horizon, respectively.

Put together, all the experiments comprise 19 configurations, with 801 subjects and 109,800 observations of actions in total. Table~\ref{tb:all-experiment-configurations} summarizes the number of subjects in each experimental configuration. The sample size is in line with (if not larger than) common academic practice; see the representative examples in the \textit{Management Science replication project} \citep{DavisFlickerHyndmanKatokKepplerLeiderLongTong2023}. Because the qualitative findings are similar across samples, the main text focuses on the baseline experiment. Other results for the online replication, preregistered replication, three-arm extension, and no-summary-table extension, together with the targeted robustness checks, are reported in Appendix~\ref{section:appendix:auxiliary experiments}.
\begin{table}[htbp]
	\centering
		\caption{Summary of Experiments}
		\label{tb:all-experiment-configurations}
		\small
		\begin{tabular}{lcccc}
			\toprule
			Experiment 
			& Reward mean vector
			& \(T=100\)
			& \(T=200\)
			& Total \\
			
			\midrule
			
			\multicolumn{5}{c}{\textit{Main Experiments and Replications}} \\[0.3em]
			
			\multirow{2}{*}{Laboratory (Baseline)}
			& \((0.4,0.5)\) & 35 & 33 & 68 \\
			& \((0.6,0.4)\) & 29 & 26 & 55 \\
			\addlinespace
			
			\multirow{2}{*}{Online replication}
			& \((0.4,0.5)\) & 32 & 30 & 62 \\
			& \((0.6,0.4)\) & 33 & 28 & 61 \\
			\addlinespace
			
			\multirow{2}{*}{Preregistered  replication}
			& \((0.4,0.5)\) & 76 & 75 & 151 \\
			& \((0.6,0.4)\) & 78 & 74 & 152 \\
			\midrule
			
			\multicolumn{5}{c}{\textit{Auxiliary Experiments}} \\[0.3em]
			
			Three-arm extension
			& \((0.6,0.4,0.2)\) & 30 & 31 & 61 \\

			No-summary-table extension
			& \((0.6,0.4)\) & 29 & -- & 29 \\
			
			Targeted robustness checks
			& \((0.6,0.4)\) &  &  &  \\[0.2em]
			
			{\footnotesize \quad Control}
			& & 42  & -- & 42 \\
			
			{\footnotesize \quad Practice rounds} 
			& & 40 & -- & 40 \\
			
			{\footnotesize \quad Higher stakes}
			& & 40 & -- & 40 \\
			
			{\footnotesize \quad Remaining-period reminder}
			& & 40 & -- & 40 \\
			\midrule
			
			Overall
			&  & 504 & 297 & 801 \\
			\bottomrule
		\end{tabular}
		
		\begin{tablenotes}[flushleft]
			\scriptsize
			\item[] \textit{Note:} The main text focuses on the baseline experiment. Analyses of the replications and auxiliary experiments are reported in Appendix~\ref{section:appendix:auxiliary experiments}.
		\end{tablenotes}
\end{table}

\subsection{Lab Evidence of the EE Trade-off}\label{subsec:behavioral-data-descriptive-analysis}

To motivate our model, we start with some high-level evidence regarding whether and how real people navigate the EE trade-off. We illustrate this through subjects' tendency to choose the \textit{leading} arm, i.e., the arm with the higher empirical mean of realized point rewards. (Note that the leading arm does not have to be the optimal arm, but just the one with better performance at the time of choosing.) Naturally, one would expect that if decisions are purely driven by exploitation, the actions should be purely ``greedy,'' i.e., the leading arm is always chosen. Conversely, if subjects are solely driven by exploration, the leading arm will be as likely as, if not less likely than, the other arms to be chosen. Next, we describe what we observe from the experiments.

\vskip 0.2 cm
\noindent\textbf{Overall frequency of choosing the leading arm.}  Our experimental data reveal that only 2 out of 123 subjects chose the leading arm throughout the entire time horizon and 100 out of 123 subjects chose the leading arm in at least 50\% of the rounds. On average, subjects chose the leading arm 64.70\% of the time.\footnote{This number is calculated by first recording each individual's fraction of rounds of choosing the leading arm and then taking an average across individuals. One would get a similar figure if taking the average at the individual-period level, which is $ 64.17\% $.} That ratio is strictly less than one, yet significantly larger than $ 0.5 $ (Wilcoxon signed-rank test right tail, $ p  = 8.68 \times 10^{-15}$). These findings indicate that, at the aggregate level, the participants are choosing the arms in a way that is neither purely exploration nor pure exploitation, as expected.

\vspace{0.2 cm}
\noindent\textbf{What influences choosing the leading arm?}  The choice probability of the leading arm may not be a constant. Intuitively, the EE trade-off would move towards \textit{exploitation} (i.e., exploration is reduced) if the subject is faced with \textit{stronger evidence} that one arm is the superior one, and vice versa. To this end, we also provide some qualitative evidence regarding factors that impact the probability of choosing the leading arm.

\begin{figure}[htbp]
    \centering
    \begin{subfigure}[b]{0.33\textwidth}
        \centering
        \includegraphics[width=\textwidth]{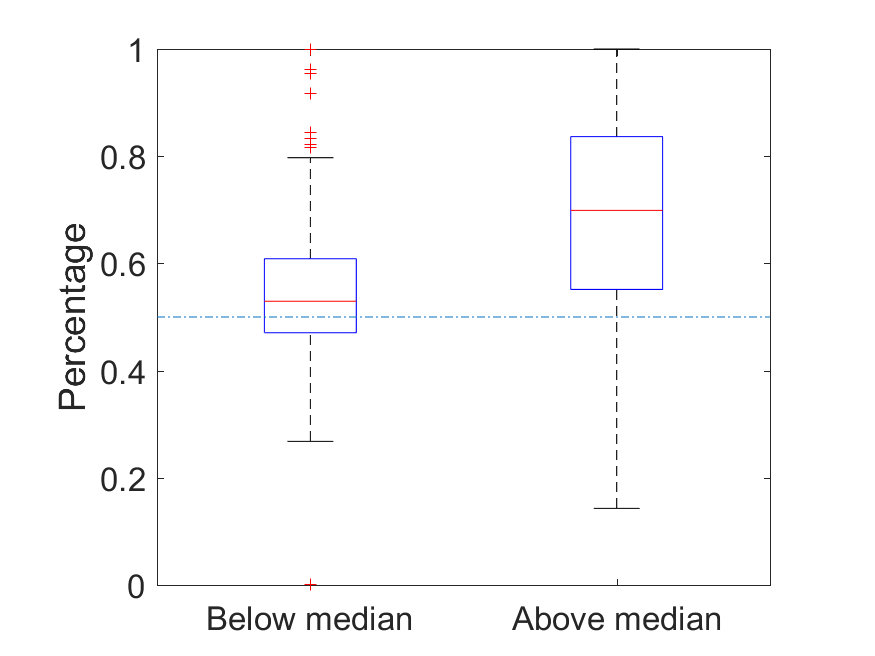}
        \caption{By empirical mean difference}
    \end{subfigure}
    \hspace{-2em}
    \begin{subfigure}[b]{0.33\textwidth}
        \centering
        \includegraphics[width=\textwidth]{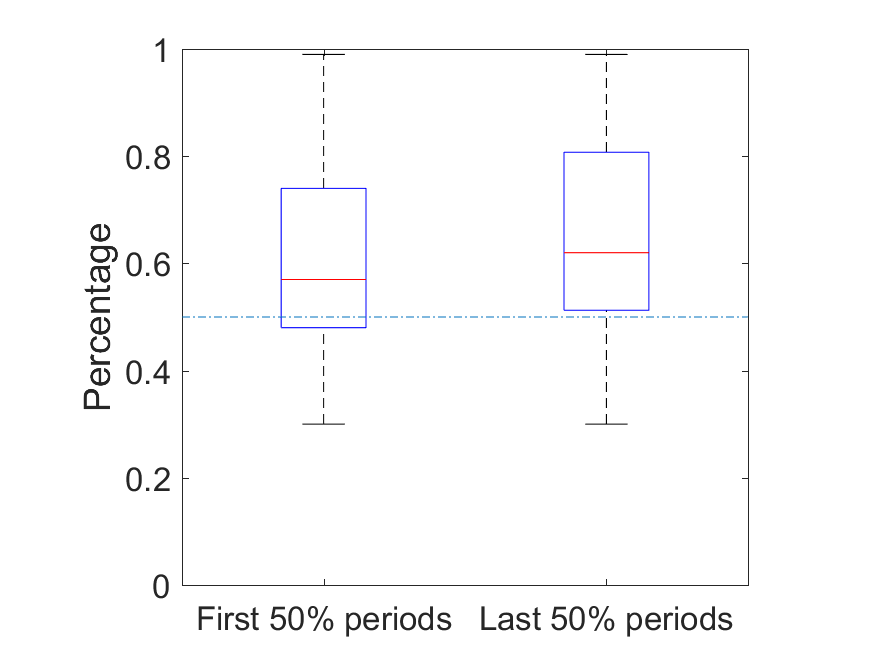}
        \caption{By period}
    \end{subfigure}
    \hspace{-2em}
    \begin{subfigure}[b]{0.33\textwidth}
        \centering
        \includegraphics[width=\textwidth]{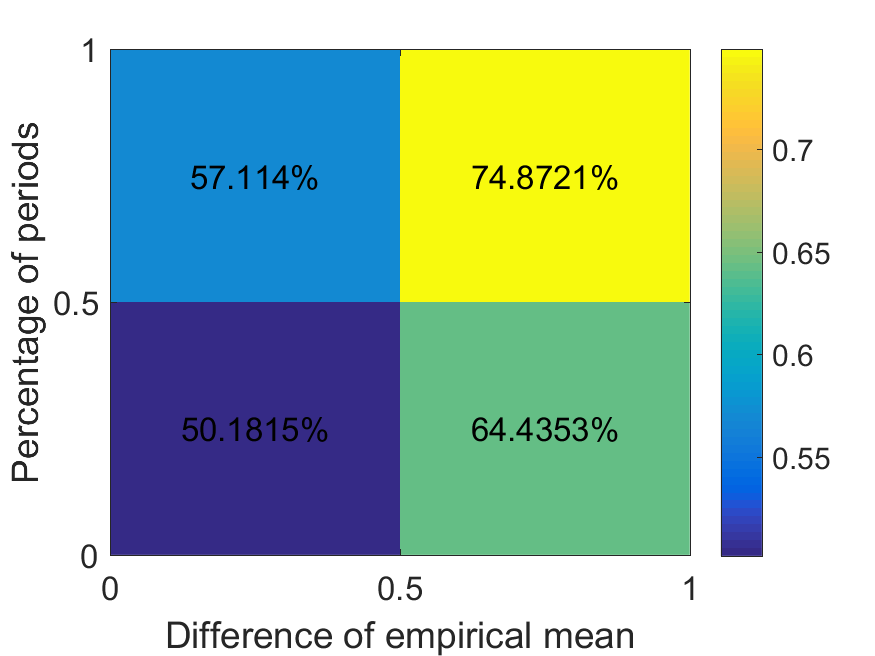}
        \caption{By both}
    \end{subfigure}
	
    \caption{Fraction of Periods Choosing the Leading Arm}
    \label{fig:box-num-larger}
 
\end{figure}

We summarize our findings in Figure~\ref{fig:box-num-larger}. For every fixed subject, we partition the time horizon into two groups: the first group contains periods corresponding to the highest $ 50\% $ of empirical mean gaps at that time, and the second group corresponds to the lowest $ 50\% $. We then calculate the fractions of choosing the leading arm for both groups; see Appendix \ref{section:appendix:example} for a concrete example of how we calculate the fraction. The distributions of both fractions across subjects are reported in Panel (a). We find a higher average fraction in the first group (two-sided Wilcoxon signed-rank test, $ p = 3.51 \times 10^{-7}$). In Panel (b), we go through a similar process but divide the time horizon into the first 50\% and last 50\% of periods. We observe a higher average fraction in the second group (two-sided Wilcoxon signed-rank test, $ p = 3.62 \times 10^{-5}$). Finally, in Panel (c), we partition the time horizon by both empirical mean gap and time period and report the average fraction in each quadrant. Quite intuitively, we find that the leading arm is more likely to be chosen when the gap is wider or based on more observations.  That is consistent with our expectation that the subjects are balancing the EE trade-off throughout the process.

\section{QCARE: A Quantitative Behavioral Model }\label{sec:QCARE-Intro}

We now proceed to explain our quantitative behavioral model to help provide deeper knowledge about people's EE trade-offs beyond descriptive statistics.

\subsection{Model Definition}
  
When faced with the task of choosing among multiple options, most popular policies (regardless of algorithmic or behavioral) are based on scoring systems. That is, a policy assigns a score to each choice and selects the one with the highest score. In the context of the MAB problem, there are many ways such scores can be generated. For instance, the Gittins index is a scoring system that corresponds to a Bayes optimal policy in the infinite-horizon discounted setting. The Thompson Sampling policy assigns a score to each arm by randomly drawing from the posterior distribution of its reward parameters. Similarly, as suggested by \citet{GansKnoxCroson2007}, humans may also behave according to a (randomized) score rule.

In this paper, we follow the convention of basing the choices on scoring rules. 
At any period $ t $ and given any arm $ i $, we consider two quantities below.
\begin{enumerate}
	\item We use $k_i(t):= \sum_{\ell=1}^{t-1} \mathbb{I}\{a(\ell) = i\}$ to denote the number of times this arm has been pulled. It quantifies the length of experience the learner has with the arm.
	\item We represent this arm's historical reward by $\hat{\mu}_i(t):= \frac{\sum_{\ell=1}^{t-1}r_i(\ell)\mathbb{I}\{a(\ell)=i\}}{k_i(t)+1}$, where $r_i(\ell)$ is this arm's (observed) reward at period $\ell$. This term quantifies arm $ i $'s historical performance.
\end{enumerate}
The notation for both terms is consistent with the earlier convention (e.g., \citealp{agrawal2017near}). The Quantal Choice with Adaptive Reduction of Exploration (QCARE) model we propose depicts a scoring system defined below:
\begin{equation}\label{eq:QCARE-definition}
\theta_i(t)=\hat{\mu}_i(t)+ \frac{1}{(k_i(t) + 1)^{\alpha}} \epsilon_{it},
\end{equation}
where $\epsilon_{it}$ is an i.i.d.\ standard normal random variable, independent across arms and periods and independent of the reward processes. The decision maker then chooses the arm with the highest score in every period, i.e., $ a(t)=\argmax_i \{\theta_i(t)\} $.

Let us walk through every component of the score $ \theta_i(t) $. The first term, $ \hat{\mu}_i(t) $, gives higher scores to the arms with better historical performance. In this way, it signifies the exploitation aspect of behavior. The second term corresponds to a random shock that enables exploration behavior. The random shock shrinks with $ k_i(t) $, i.e., the decision maker's number of interactions with this particular arm. When the experience is limited, pulling this arm can still be encouraged, even if the historical performance is suboptimal, as reflected in the large random shock term. When there have been sufficient interactions, i.e., when $ k_i(t) $ is large, the attractiveness of this arm is almost purely determined by its historical performance.

\subsection{Model Interpretation and Discussion}
\label{subsec: model_interpretation}

Under our model, the balance between exploration and exploitation is captured through the relative magnitude between the exploitation and exploration terms, or more specifically, how fast the random shock shrinks. If the random shock shrinks too fast, the decision maker might be misled by the historical performance of an arm based on too few interactions. Conversely, if the random shock shrinks too slowly (or does not shrink at all), the decision maker may end up spending too much time on suboptimal arms despite overwhelming evidence that a certain arm is the best.

\vspace{0.2 cm}
\noindent \textbf{The reduction rate of exploration.} We use a parametric form to present how the random shock shrinks over time, namely, $\big(k_i(t) + 1\big)^{-\alpha}$. This form is parameterized by a single parameter, $ \alpha $, which we call the \textit{reduction rate of exploration}. As its name suggests, a larger $ \alpha $ corresponds to a faster shrinkage of the random shock term, and hence a higher weight on exploitation given any system state. In the extreme case, if $ \alpha $ approaches infinity, the exploration term becomes negligible and hence QCARE becomes purely \textit{greedy}. That is, one simply chooses the leading arm exclusively. Conversely, if $ \alpha \leq 0 $, the random shocks remain significant regardless of the learner's historical experience. In this case, QCARE never settles on any arm.

The introduction of $ \alpha $ allows us to quantify the EE trade-off using a single easy-to-interpret parameter. In later sections, we formally discuss how different values of $ \alpha $ create different system dynamics and decision qualities. In addition, we will estimate the implied-$ \alpha $ values from the behavioral data collected from the experiments. That helps shed light on the behavioral patterns when \textit{human beings} manage the balance between exploration and exploitation.

\vskip 0.2 cm
\begin{remark} {\sf 
		\textbf{Time-position effects.} Note that QCARE captures the reduction of exploration in an ``any-time'' manner by keeping $ \alpha $ a single constant. That is, the reduction of exploration does not depend on the explicit time position \(t\), but only on the arm pull counts. That makes QCARE similar to the Thompson Sampling policy, which we will discuss more later. In principle, the time-position effect can matter, too. That can be captured by letting $ \alpha $ be a function of $ t $ rather than a constant. To assess whether it is sensible to do so, we perform a diagnostic test, where we split the decision sequences into two halves and estimate \(\alpha\) separately. We do not find robust evidence that the $ \alpha $ values are different in the first and second halves; see Appendix~\ref{section:appendix:subject pools} for more details. We therefore keep the constant-$ \alpha $ setting as our main model, and leave time-position dependence as an extension for future research.		
	}
\end{remark}

\vspace{0.2 cm}
\noindent \textbf{Relating to Thompson Sampling.} QCARE has an intimate connection to Thompson Sampling (TS), a well-celebrated policy from the online learning literature. The basic idea of TS is to maintain the posterior distributions of the arm reward parameters and randomly pull an arm according to its posterior probability of being optimal; see Appendix \ref{appendix:sec:algorihims} for a more detailed algorithmic description. TS is asymptotically optimal in terms of various notions of regret as $ T $ grows, such as instance-independent regret \citep{agrawal2013further, agrawal2017near}, instance-dependent regret (\citealp{kaufmann2012thompson}), and Bayesian regret \citep{russo2014learning}. QCARE generalizes TS; see the following result for a formal statement.

\begin{proposition}\label{prop:TS_QCARE_equa}
	When $\alpha=0.5$, QCARE is equivalent to Gaussian Thompson Sampling.
\end{proposition}

This observation has two implications. First, it enables us to view QCARE as a parameterized generalization of (Gaussian) Thompson Sampling:  when  $ \alpha < 0.5 $ (resp. $ \alpha > 0.5 $), our model captures a policy that explores more (resp. less) aggressively than TS. In this way, we establish a novel connection between Thompson Sampling and quantal response models, thus bridging the seemingly distinct online learning and behavioral literature. Second, given the order-optimality of Gaussian Thompson Sampling, we may interpret $ \alpha = 0.5 $ as the optimal benchmark for the EE trade-off when $ T $ is infinite. 

We also note that Thompson Sampling can be specified with different priors and posterior sampling rules. To that end, we include another Thompson Sampling specification, i.e., Thompson Sampling with a Beta prior (TS-B), in the empirical comparison in Section~\ref{subsec:empirical analysis}. In addition, a few other generalizations of TS have been considered in the online learning/psychology literature. For example, \citet{chapelle2011empirical} reshaped the posterior of TS and showed that shrinking the posterior to reduce exploration can benefit the non-asymptotic performance of TS. \citet{bouneffouf2017bandit} extended TS to a parametric version that allows for the incorporation of biases associated with mental disorders. %
Although, in principle, they could also be used to parameterize the EE trade-off and be fit to behavioral data, it is important to note that those methods are numerical and heuristic in nature. In particular, none of them provides a rigorous treatment of quantifying over- and under-exploration and how that affects the decision quality.

\vskip 0.2 cm
\noindent\textbf{Relating to UCB-type policies.} QCARE is also similar to Upper Confidence Bound (UCB) policies (\citealp{lai1985asymptotically,auer2002finite,jamieson2014lil}), since both rely on an index that decomposes into an exploitation and an exploration component. The two differ in how exploration is generated, as can be seen from the generalized score
	\begin{equation}\label{eq:generalized-score}
	\theta_i(t)=\hat{\mu}_i(t)+\frac{m_t}{(k_i(t)+1)^{\alpha}}+\frac{\sigma\,\epsilon_{it}}{(k_i(t)+1)^{\alpha}},
	\end{equation}
where the shock follows \(N(m_t,\sigma^2)\) for some \(m_t\ge0\) and \(\sigma\ge0\). The second term is a deterministic bias toward arms with fewer pulls, which induces exploration through \textit{optimism}. The third term induces exploration through \textit{randomization}. With \(\alpha=0.5\), \(m_t=c\sqrt{\log t}\), and \(\sigma=0\), the score reduces to the UCB index \(\hat{\mu}_i(t)+c\sqrt{\log t/(k_i(t)+1)}\), and with \(m_t=0\) and \(\sigma=1\), it reduces to QCARE. These different exploration mechanisms can nevertheless produce similar long-run learning dynamics. For example, UCB and Gaussian Thompson Sampling (QCARE with \(\alpha=0.5\)) both attain logarithmic instance-dependent regret, and their instance-independent regrets differ by at most a logarithmic factor \citep{auer2002finite,agrawal2017near,han2024ucb}. \citet{fan2022typical} also establish the same strong laws of large numbers and central limit theorems for their suboptimal-arm pull counts.

It is, however, worth noting that QCARE appears to be a more natural empirical model than UCB. UCB is a deterministic policy and assigns zero probability to every arm other than the one with the highest index, so its likelihood vanishes whenever a single observed choice deviates from its prediction. While it is possible to modify UCB by adding noise, it is unclear how to do so in a principled way.\footnote{The usual way is to add a softmax or probit layer where the noise does not shrink with the arm pulls. But that makes the modified policy unable to achieve complete learning. Alternatively, a new model may need to be developed, such as the general model above, which involves more parameters and can be hard to analyze.} 
In comparison, QCARE naturally admits structural noise that shrinks with arm pull counts. It is governed by a single parameter \(\alpha\) that is interpretable, easy to estimate, and directly linked to decision quality, as we will demonstrate later.

\vskip 0.2 cm
\noindent\textbf{Relating to Gittins index.} The Gittins index is different from QCARE in that it is normative in nature, as it leads to the Bayes optimal policy in the discounted Bayesian MAB problem. The Gittins index is closely related to UCB, as it adds to an arm's posterior mean a deterministic exploration bonus that shrinks as the arm is sampled. In fact, \citet{russo2021note} further shows that UCB can be viewed as a limiting approximation of the Gittins index as the discount factor approaches one. Like UCB, it is deterministic and thus faces the same estimation difficulty noted above. Empirically, Gittins-index-based models fit human bandit choices poorly, even in an experiment designed so that the Gittins index rule is optimal \citep{GansKnoxCroson2007}.

\vskip 0.2 cm

 \noindent \textbf{Relating to quantal choice models.} QCARE is related to quantal choice models in that both frameworks use random shocks to describe stochastic choice. Classic foundations include the attraction model \citep{Luce1959} and the random utility model (RUM) framework \citep{mcfadden1976quantal}, under which an individual assigns each alternative \(i\) a random score \(\theta_i=u_i+\epsilon_i\) and chooses the alternative with the highest realized score. Different assumptions on the shock distribution lead to different choice models, such as multinomial logit, probit, and exponomial choice models \citep{alptekinouglu2016exponomial}. In this literature, random shocks are often interpreted either as unobserved utility components, including in dynamic discrete-choice models \citep{rust1987optimal}, or as behavioral errors/cognitive limitations, as in behavioral newsvendor models \citep{Su2008,SchweitzerCachon2000} and quantal response equilibrium models in games \citep{mckelvey1995quantal,mckelvey1998quantal}.

The key distinction is that sequential learning introduces an additional role for the random shock. In traditional quantal choice settings, utilities are usually structurally specified. In an MAB problem, however, the rewards of the arms are initially unknown, and some exploration is intrinsically needed for learning. Thus, the stochastic term in QCARE should not be interpreted only as an unobserved utility shock or a behavioral error. It also captures the agent's tendency to explore when empirical evidence is still limited. This is why QCARE allows the shock scale to decline with arm-level pull counts. The parameter \(\alpha\) governs how quickly this exploration component fades as the decision maker learns. In this sense, QCARE keeps the quantal-choice idea of stochastic choice, but adapts it to the exploration-exploitation trade-off in sequential learning.

\vskip 0.2 cm
\noindent \textbf{Relating to heuristic behavioral models.}  In Section~\ref{subsec:empirical analysis}, we also consider a collection of other behavioral models. While they are also based on (randomized) scores, they are heuristic in nature and lack formal theoretical understanding. A qualitative difference between QCARE and those policies is that, under those policies, the softmax layer implies that every arm is chosen with a non-diminishing probability in every period. These models therefore never stop exploring. They are hence unable to represent policies that are reasonably optimal, such as those that achieve complete learning. Our empirical validation also shows that QCARE displays favorable out-of-sample prediction power compared to those empirical heuristics.

\section{Theoretical Characterizations of QCARE}\label{sec:QCARE-Theory}

This section provides the theoretical characterization of QCARE. We first establish its qualitative choice properties, complete-learning result, and reduced-form interpretation of \(\alpha\). We then characterize regret across different values of \(\alpha\), provide key steps and intuition of a general Markovian-policy framework to explain the underlying mechanisms, and examine how the asymptotic insights extend to finite horizons.

\subsection{Qualitative Properties}\label{subsec:QCARE-proprieties}

In this subsection, we derive several qualitative properties of QCARE. Recall from Equation~\eqref{eq:QCARE-definition} that the arm-pulling probabilities only depend on the history through the profile of every arm's pull count and empirical mean reward. Let that be encoded by a state \(S=(\kappa,u)\), where \(\kappa=(\kappa_1,\ldots,\kappa_N)\in \mathbb{Z}_+^N \) records the arm-level pull counts and \(u=(u_1,\ldots,u_N) \in [0,1]^N\) records the corresponding empirical mean rewards. In that light, let us use \(\hat{i}^\ast:=\argmax_i u_i\) to represent the leading arm, and let \(\hat{Q}(S):=\Pr{a(t)=\hat{i}^\ast\mid S(t)=S} \) denote the probability of selecting it. The following result summarizes comparative-statics properties of QCARE in terms of \(\hat{Q}\).

\begin{proposition}\label{prop:QCARE-basic-properties} (Comparative statics)
	Pick any $\alpha>0$,  $ N >1 $, and state $ S  = (\kappa_1, \ldots, \kappa_N, u_{1}, \ldots, u_{N})$.
	The following facts about QCARE hold. 
	\begin{enumerate}
		\item  For every \(0<\delta\leq 1-u_{\hat i^\ast}\), we have  $\hat{Q}(S) \leq \hat{Q}(\kappa_1, \ldots, \kappa_N, u_{1}, \ldots, u_{\hat{i}^\ast} + \delta,\ldots, u_{N})$.
		\item For every integer \(\delta \geq1\) and every \(j\neq\hat i^\ast\), we have $\hat{Q}(S) \leq \hat{Q}(\kappa_1, \ldots, \kappa_j+\delta, \ldots, \kappa_N, u_{1}, \ldots,  u_{N})$.	
	\end{enumerate}
\end{proposition}

Proposition~\ref{prop:QCARE-basic-properties} reveals some structural properties regarding how QCARE balances exploration with exploitation. The emphasis on exploitation (i.e., probability of choosing the leading arm) becomes more salient when the leading arm has established a larger empirical mean gap compared to the other arms, or when the gap is based on more observations of the non-leading arms. These properties are qualitatively consistent with the descriptive findings from the experimental evidence summarized in Figure \ref{fig:box-num-larger}.

\noindent\textbf{A reduced-form representation of QCARE.} In the long run, as the arm pulls accumulate, $ \hat{Q}(S) $ also admits a simple reduced-form representation that cleanly reveals the effect of $ \alpha $. For simplicity, consider the two-arm case, and let the probability of pulling the non-leading arm be denoted as \(Q_{\mathrm{non\text{-}lead}}(S):=1-\hat{Q}(S)\). It then admits the following form:
	\begin{align}\label{eq:reduced form}
	\frac{1}{2}\log\log\left(\tfrac{1}{Q_{\mathrm{non\text{-}lead}}(S)}\right)
	=
	\alpha\log\left(\min\{\kappa_1+1,\kappa_2+1\}\right)+O(1).
	\end{align}
It can be seen that, on the transformed scale above, $ \alpha $ corresponds to the slope at which the non-leading arm choice probability shrinks with arm pulls. All other quantities, such as the empirical mean gap, enter only through the intercept term. This representation therefore provides a direct interpretation of $\alpha$ as the parameter governing the reduction rate of exploration. The formal derivations, as well as generalizations to the multi-arm case, are contained in Appendix~\ref{appendix:reduced-form-qcare}.

Next, we explore another important asymptotic property of QCARE. To put it briefly, QCARE can ``learn'' the optimal arm and converge to it, at least in the long run. We formally state this property in the result below.
\begin{theorem}\label{thm:UB_sublinear} (Complete learning)
	Pick any $\alpha>0$,  $ N >1 $, and the bandit problem instance $ \Theta $. The following facts about QCARE hold. 
	\begin{itemize}
		\item The arm-pulling process converges to the optimal arm. That is, $ \Pr{a(t) = i^\ast} \to 1 $ as $ t \to \infty $;
		\item The regret satisfies $\mathcal{R}^\Theta(T)= o(T)$. As a result, QCARE with parameter $ \alpha $ represents a long-run-average optimal policy, i.e., one that maximizes $ \limsup_{T \rightarrow\infty}\  \frac{1}{T} \ \mathbb{E}\left[\sum_{t=1}^{T}r_{a(t)}(t)\right] $.
	\end{itemize} 
\end{theorem}

The result above suggests that all values of $\alpha > 0$ correspond to ``plausible'' behavioral MAB policies as they all achieve complete learning and are long-run-average optimal. This sense of optimality dates back to (\citealp{robbins1952some}).  While this property is not particularly stringent, it should not be taken for granted. For example, the fixed-\(\varepsilon\)-greedy policy \citep{auer2002finite} used as a benchmark in our empirical analysis is not long-run-average optimal because it keeps randomizing with a non-vanishing probability and therefore incurs a linear regret. In contrast,
QCARE possesses a minimally reasonable level of exploration to gradually learn the optimal arm. The fact that people do tend to learn over time is reflected both in the aforementioned descriptive evidence and in the favorable prediction power of QCARE relative to the fixed-\(\varepsilon\)-greedy benchmark, as will be discussed later.

Beyond intuitive qualities, QCARE offers a principled (and also parsimonious) way to quantify the EE trade-off by measuring how much one keeps exploring in the face of accumulating evidence, which we develop next.

\subsection{Quantifying the EE Trade-off: An Asymptotic Theory}\label{sec:Asymptotic Theory}

In this subsection, we provide a formal treatment of how $ \alpha $ relates to the decision qualities, especially the effects of ``over'' and ``under'' exploration.

\begin{figure}[htbp]
    \centering
    \begin{subfigure}[b]{0.4\textwidth}
        \centering
        \includegraphics[width=\textwidth]{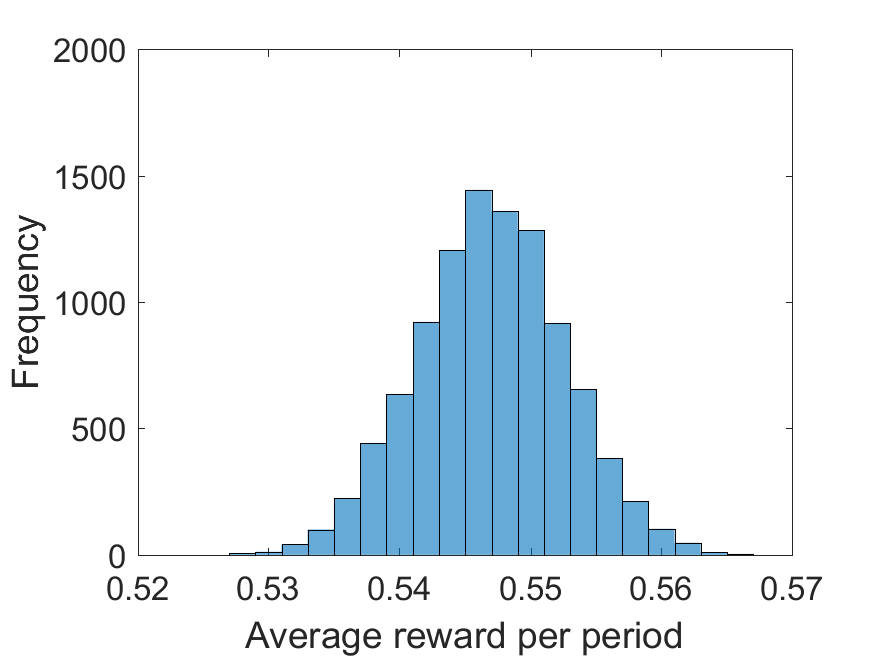}
        \caption{$\alpha = 0.2$, $T = 10^4$}
    \end{subfigure}
    \hspace{1cm}
    \begin{subfigure}[b]{0.4\textwidth}
        \centering
        \includegraphics[width=\textwidth]{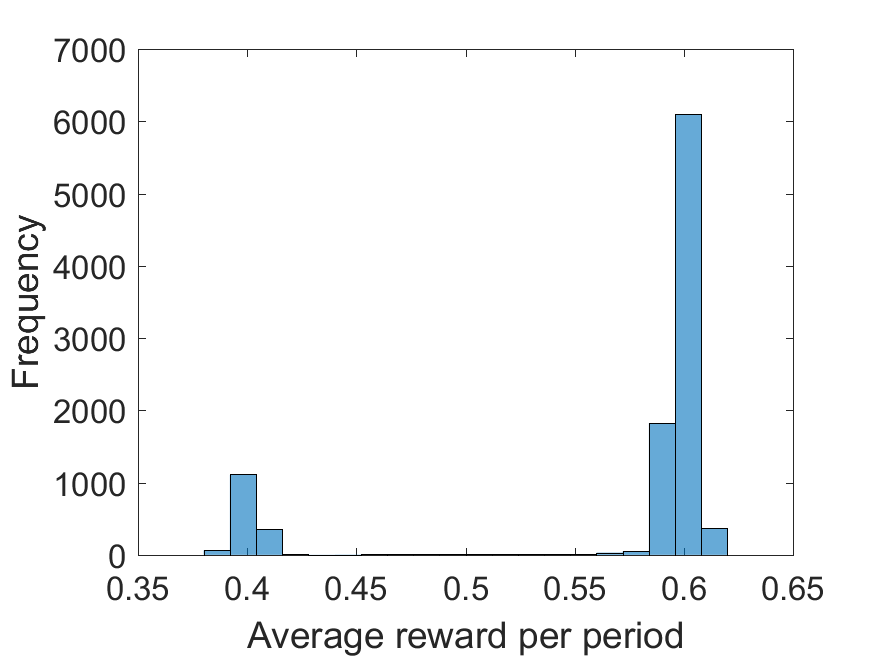}
        \caption{$\alpha = 2$, $T = 10^4$}
    \end{subfigure}
    
    \vspace{0.5cm}

    \begin{subfigure}[b]{0.4\textwidth}
        \centering
        \includegraphics[width=\textwidth]{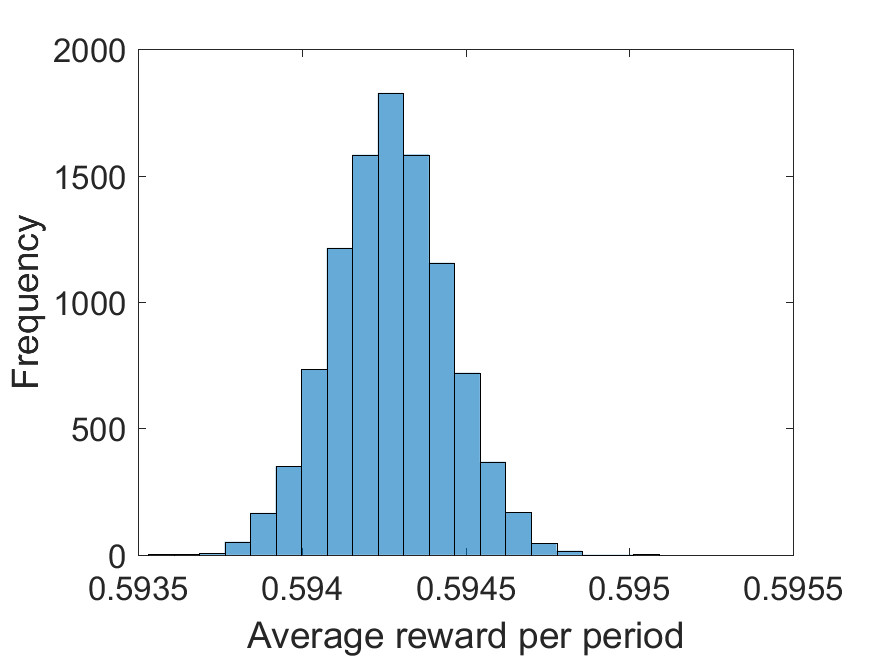}
        \caption{$\alpha = 0.2$, $T = 10^7$}
    \end{subfigure}
    \hspace{1cm}
    \begin{subfigure}[b]{0.4\textwidth}
        \centering
        \includegraphics[width=\textwidth]{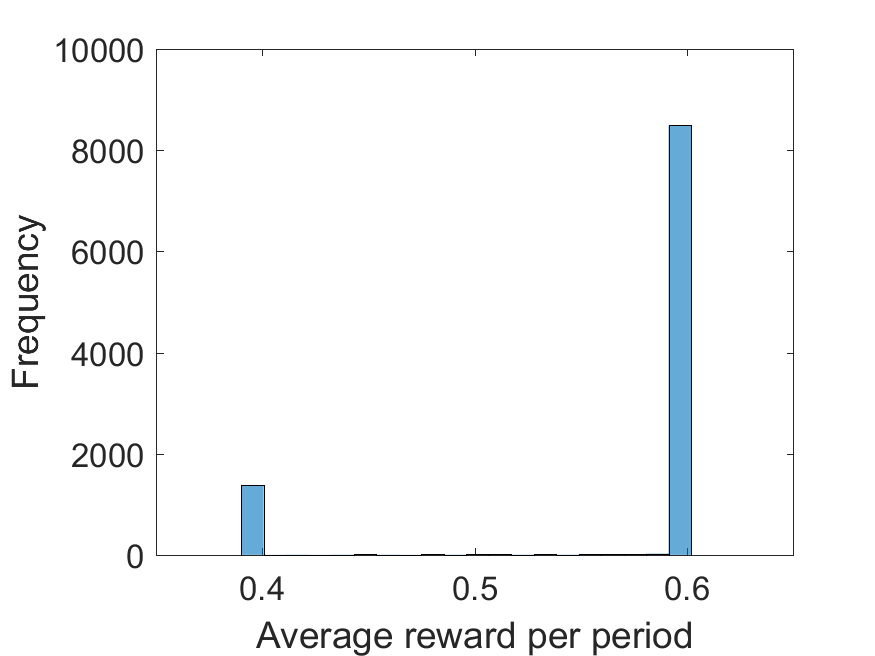}
        \caption{$\alpha = 2$, $T = 10^7$}
    \end{subfigure}
 
 \caption{Reward Distribution Under QCARE With Different $\alpha$}
	\label{fig:reward-bimode}
\end{figure}

\vspace{0.2 cm}
\noindent \textbf{Illustrative plots.}  To build intuition, let us first use simulation to illustrate how $ \alpha $ influences the reward distribution under QCARE. Figure~\ref{fig:reward-bimode} depicts the histograms of the average rewards per period in a two-arm bandit problem with different values of $ \alpha $ and $ T $ based on a simulation study of 10,000 sample paths. The arms' reward distributions are fixed to be Bernoulli with means $ (0.6, 0.4) $. In the figure, Panels (a) and (c) correspond to $\alpha=0.2$, representing high exploration. We can observe that the realized reward is distributed around a mean value less than $ 0.6 $, and the distribution tightens as $ T $ becomes large. That reflects how heavy exploration, although helping to learn the optimal arm eventually, may spend too much time on the inferior arm. Conversely, Panels (b) and (d) correspond to $ \alpha = 2 $,  representing low exploration. Here we observe a distinctly bimodal pattern in the reward distribution: around $15\%$ of the sample paths display an average reward of around $ 0.4 $, the expected reward of the inferior arm. This indicates that while most samples quickly concentrate around the optimal value of $ 0.6 $, a significant portion of the samples are stuck in the inferior arm due to insufficient exploration. We refer to this as finite-horizon premature lock-in. These two patterns preview the theoretical analysis below. An overly slow reduction rate of exploration leads to persistent sampling of inferior arms, while overly fast exploration reduction can create prolonged finite-horizon lock-in and substantial lower-tail risk.

\vspace{0.2 cm}
\noindent \textbf{The regret landscape of QCARE.}
We now present our main theoretical results, which characterize QCARE's performance in the asymptotic setting as \(T\) approaches infinity. To this end, we find it helpful to formally define the standard Big O notation. For every pair of functions $ f(\cdot), g(\cdot) : \mathbb{Z}_+ \to \mathbb{R}_+ $, we say that: 
\begin{align*}
f(T) = 
\begin{cases}
O(g(T)), \quad \text{ if }  &\limsup_{T \rightarrow\infty} f(T)/g(T) < \infty;\\
o(g(T)), \quad \text{ if }  &\limsup_{T \rightarrow\infty} f(T)/g(T) =0;\\
\Omega(g(T)), \quad \text{ if }  &\liminf_{T \rightarrow\infty} f(T)/g(T) > 0;\\
\omega(g(T)), \quad \text{ if }  &\liminf_{T \rightarrow\infty} f(T)/g(T) = \infty.\\
\end{cases}
\end{align*}
Recall that when $\alpha=0.5$, QCARE reduces to Gaussian Thompson Sampling and hence achieves the optimal regret order of $ O(\sqrt{T}) $. When $\alpha$ deviates from the optimal benchmark of $ 0.5 $, distinctly different patterns emerge.
We first consider the case where $ \alpha < 0.5 $.

\begin{theorem}\label{thm:smaller0.5}
	Let \(\alpha\in(0,0.5)\) and \(N>1\) be fixed. There exist constants \(0<c\le C<\infty\), depending only on $ \alpha $ and $ N $, such that the following statements hold.
	\begin{enumerate}[label=(\roman*)]
		\item \emph{(Instance-dependent regret)} For every instance \(\Theta\), there exists \(T_0=T_0(\alpha,N,\Theta)<\infty\) such that, for every \(T\ge T_0\),
		\[
		c\,(\log T)^{1/(2\alpha)}\sum_{i\neq i^\ast}\Delta_i^{1-1/\alpha}
		\;\le\; \mathcal{R}^{\Theta}(T) \;\le\;
		C\,(\log T)^{1/(2\alpha)}\sum_{i\neq i^\ast}\Delta_i^{1-1/\alpha}.
		\]

		\item \emph{(Instance-independent regret)} There exists \(T_0'=T_0'(\alpha,N)<\infty\) such that, for every \(T\ge T_0'\),
		\[
		c\,T^{1-\alpha} \;\le\; \mathcal{R}(T) \;\le\; C\,T^{1-\alpha}\sqrt{\log T}.
		\]
	\end{enumerate}
\end{theorem}

Theorem~\ref{thm:smaller0.5} characterizes the regret of QCARE, in both the instance-dependent and worst-case senses, with (nearly) matching upper and lower bounds. In the former case, the regret is of order \((\log T)^{1/(2\alpha)}\), and in the latter case, the regret is of order \(T^{1-\alpha}\) up to a logarithmic factor. In comparison, the $ \alpha = 0.5 $ benchmark achieves \(\mathcal{R}^\Theta (T) = O(\log T)\) and \(\mathcal{R} (T) = O(\sqrt{T})\); see \citet{agrawal2017near}. As a result, it can be seen that regardless of the regret measure, the regret deteriorates slowly as $ \alpha $ decreases below $ 0.5 $. That reflects the cost of over-exploration, as the arm-pulling probabilities concentrate on the leading arm too slowly. In fact, QCARE may continue to allocate non-negligible probability to multiple arms even when all arms have been pulled \(\Omega(T)\) times and their empirical reward differences are of order \(\Omega(T^{-\alpha})\). We defer the more in-depth analysis of the proof to Appendix~\ref{appendix:upper-bound-alpha-lt-half}.

Next, we turn to the case where $\alpha > 0.5$.  

\begin{theorem}\label{thm:lb_greater0.5}
Suppose $\alpha>0.5$ and $N>1$ are fixed. For every interior Bernoulli instance $\Theta$ satisfying $0 < \mu_i < \mu_{i^\ast} < 1$ for every $ i \neq i^\ast $, QCARE pulls suboptimal arms $T-o(T)$ times with probability $\Omega(T^{-o(1)})$. Consequently, for every such instance $\Theta$,
\[
\mathcal{R}^{\Theta}(T)=\Omega(T^{1-o(1)}),
\]
and the instance-independent regret satisfies
$
\mathcal{R}(T)=\Omega(T^{1-o(1)}).
$
\end{theorem}

\begin{figure}[ht]
	\centering
	\begin{subfigure}[b]{0.45\textwidth}
		\centering
		\includegraphics[width=\textwidth]{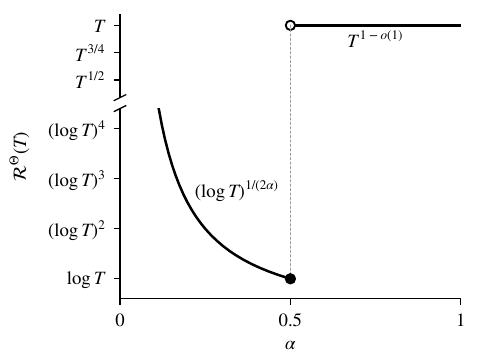}
		\caption{Instance-dependent regret}
		\label{fig:theory_f_instance}
	\end{subfigure}
	\begin{subfigure}[b]{0.45\textwidth}
		\centering
		\includegraphics[width=\textwidth]{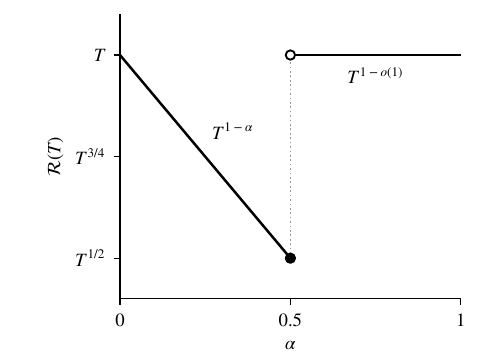}
		\caption{Instance-independent regret}
		\label{fig:theory_f_problem}
	\end{subfigure}
	\caption{Asymptotic Regret Order of QCARE as a Function of $\alpha$}
	\label{fig:theory_f}
\end{figure}

Together with the upper bounds $\mathcal{R}^{\Theta}(T)=o(T)$ and the trivial $\mathcal{R}(T)=O(T)$, Theorem~\ref{thm:lb_greater0.5} determines QCARE's regret for $ \alpha > 0.5 $. That is, although QCARE eventually learns the optimal arm, its regret jumps to a near-linear order once $\alpha$ exceeds $0.5$, however slightly. Intuitively, when $\alpha>0.5$, the arm-pulling probabilities concentrate on the leading arm too quickly. It is possible that an optimal arm and a suboptimal arm have each been pulled only $(\log T)^{\frac{1}{2\alpha-\delta}}$ times, while the pulling probability of the empirically non-leading arm has already diminished to order $O\left((\log T)^{\frac{1}{2\alpha-\delta}}/T\right)$. As a result, with probability $\Omega(T^{-o(1)})$, QCARE can be stuck on a suboptimal arm and thereby incur $\Omega(T)$ regret. We defer the more in-depth analysis of the proof to Appendix~\ref{appendix:formal-general-framework}.
Combining Theorems \ref{thm:smaller0.5} and \ref{thm:lb_greater0.5}, our asymptotic analysis covers the whole spectrum of $\alpha > 0$, and we illustrate the relationship between $ \alpha $ and asymptotic regret in Figure~\ref{fig:theory_f}.

\vspace{0.2 cm}
\noindent \textbf{Insights from the asymptotic analysis.} In summary, our analysis quantifies how decision quality is related to the EE trade-off asymptotically. We show that there is a \textit{sweet spot} of the EE trade-off, which supports the concepts of over- and under- exploration. It is also interesting to note that over- and under-exploration lead to different decision dynamics and harm the decision quality in different ways. The former results in insufficient focus on the leading arm, even with overwhelming supporting evidence, while the latter causes a non-negligible probability that the decision maker is trapped in the inferior arm.

\subsection{Key Proof Steps, Generalizations, and Intuition}\label{sec:Q_function}\label{sec:key_steps_generalizations}

We now summarize the general structural properties underlying the complete-learning and regret lower-bound results that extend beyond QCARE. The formal definitions, full proofs, and QCARE verification are deferred to Appendix~\ref{appendix:formal-general-framework}. The main idea is to treat QCARE as a special case of a broader class of Markovian MAB policies whose choice probabilities depend on the state
\( S(t)=\big(k_1(t),\ldots,k_N(t),\hat{\mu}_1(t),\ldots,\hat{\mu}_N(t)\big).
\)
Such a policy can be represented by a probability function \(Q=(Q_1,\ldots,Q_N)\), where \(Q_i(S(t))\) is the probability of pulling arm \(i\) in period \(t\).  We then identify several structural properties of \(Q\) that are sufficient to establish complete learning and regret lower bounds for these general Markovian policies. 

The first two properties, termed \textit{recurrence} and \textit{leading-arm convergence}, capture the complete learning phenomenon established in Theorem \ref{thm:UB_sublinear}. In short, the former property means that after an arm is pulled, there is always a non-diminishing probability that it is re-pulled. As a result, every arm is pulled infinitely many times in the long run. The latter property complements the former one by requiring the policy to concentrate on the leading arm given sufficiently strong empirical advantage of that arm. Together, these two properties are sufficient for complete learning.

We then introduce a property that describes overly heavy exploration, which we term \textit{\(\alpha\)-exploratoriness}. Roughly speaking, policies that satisfy this property allow situations where all arms have been pulled $ \Omega(T) $ times and their empirical reward differences are on the order of $ \Omega(T^{-\alpha}) $, yet the arm-pulling probabilities do not concentrate on the leading arm. In this sense, the policy does not concentrate on any particular arm despite overwhelming evidence that a certain arm is the best. For \(0<\alpha<0.5\), we show that this property yields the instance-independent lower bound \(\mathcal R(T)=\Omega(T^{1-\alpha})\). 

On the flip side, we introduce a property named \textit{\(\alpha\)-irreversibility} that describes under-exploration. Policies with this property allow situations where two arms have both been pulled only $ (\log T)^{\frac{1}{2\alpha - \delta}} $ times, yet the worse-performing arm's pulling probability has already diminished on the order of $ O \left((\log T)^{\frac{1}{2\alpha - \delta}}/T\right) $, or loosely translated, the probability of pulling the non-leading arms reduces super-exponentially fast with the pull count. That causes a non-negligible chance that the policy is mistakenly stuck in the inferior arm just because that arm performed better during the initial periods. For \(\alpha>0.5\), we show that, under the interior-Bernoulli instance conditions of Theorem~\ref{thm:under explore} and with a mild regularity condition that we call \textit{unradical}, this property yields a regret lower bound of order \( \Omega(T^{1-o(1)}) \). 

\subsection{From Asymptotics ($ T = \infty $) to Finite Horizons ($ T <  \infty $) }
\label{subsec:finite-analysis}

\begin{figure}
    \centering
    \hspace{-0.3 cm}
    \begin{subfigure}{0.23\textwidth}
        \centering
        \includegraphics[width=\textwidth]{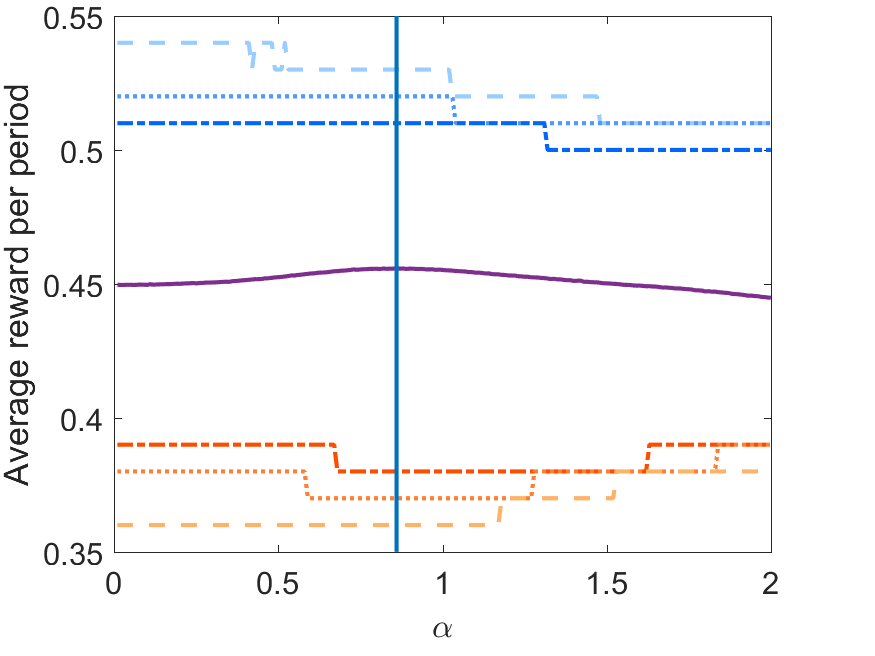}
        \caption{\footnotesize $\mu=(0.4, 0.5)$, $T=100$}
    \end{subfigure}
    \begin{subfigure}{0.23\textwidth}
        \centering
        \includegraphics[width=\textwidth]{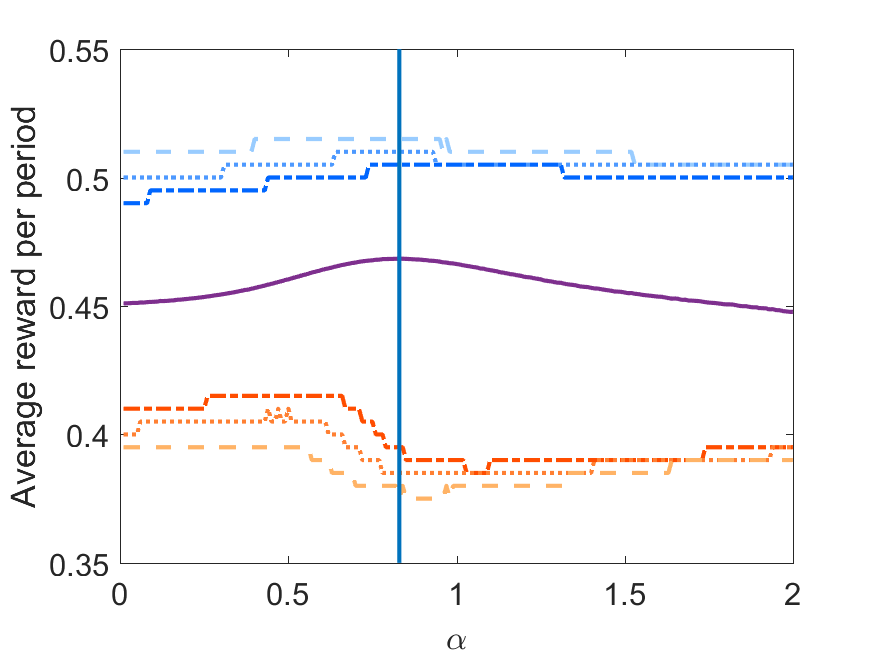}
        \caption{\footnotesize$\mu=(0.4, 0.5)$, $T=200$}
    \end{subfigure}
    \begin{subfigure}{0.23\textwidth}
        \centering
        \includegraphics[width=\textwidth]{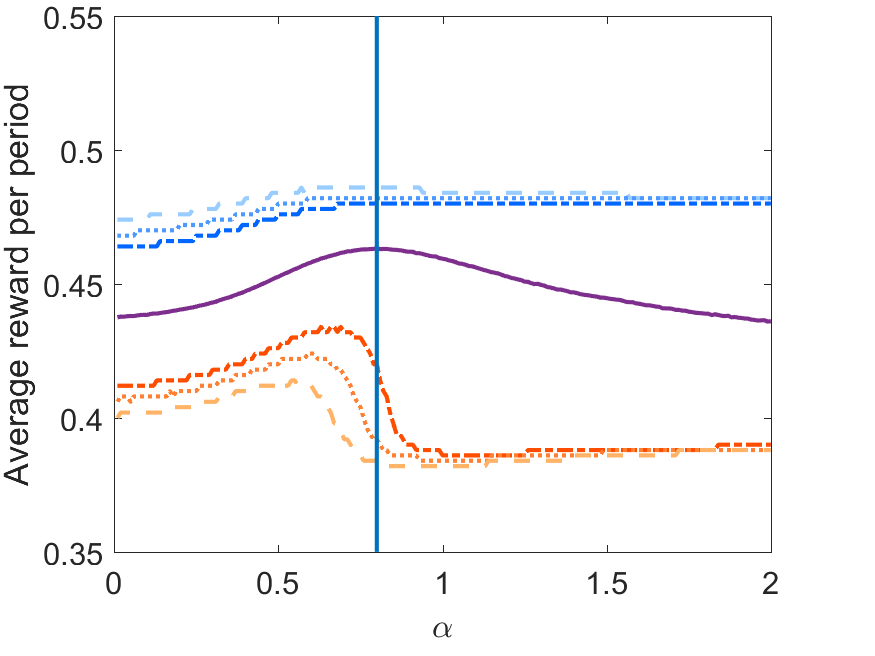}
        \caption{\footnotesize$\mu=(0.4, 0.5)$, $T=500$}
    \end{subfigure}
    \begin{subfigure}{0.24\textwidth}
        \centering
        \includegraphics[width=\textwidth]{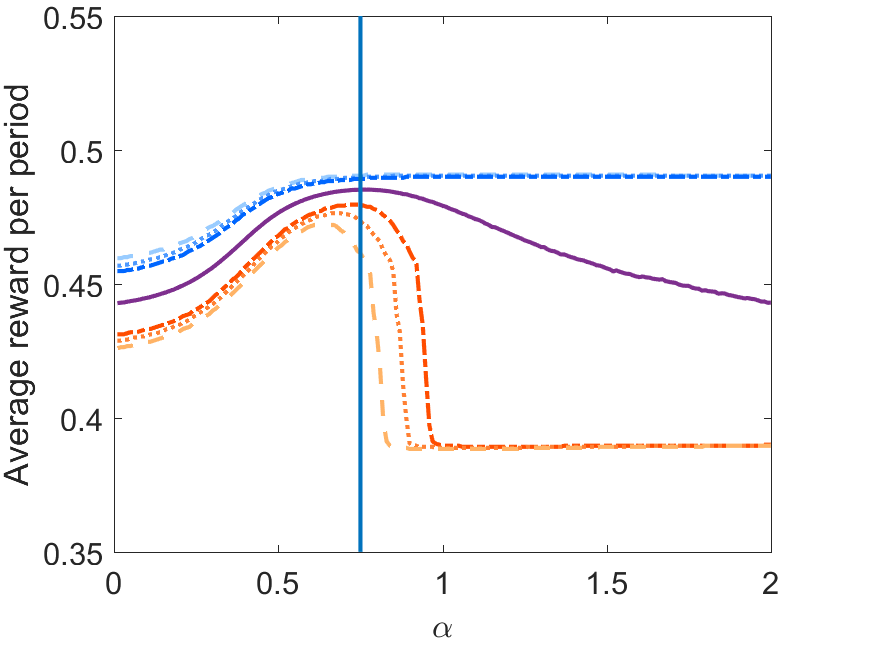}
        \caption{\footnotesize$\mu=(0.4, 0.5)$, $T=2500$}
    \end{subfigure}
	\hspace{-0.35 cm}
    \begin{subfigure}{0.06\textwidth}
        \centering
         \raisebox{0.8cm}{%
        \includegraphics[width=\linewidth]{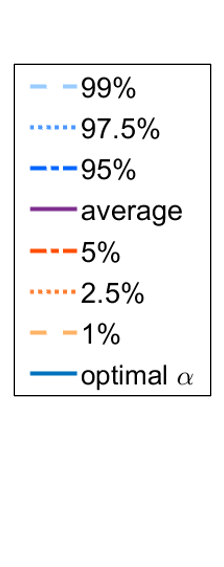}%
    }
    \end{subfigure}

	\hspace{-0.3 cm}
    \begin{subfigure}{0.23\textwidth}
        \centering
        \includegraphics[width=\textwidth]{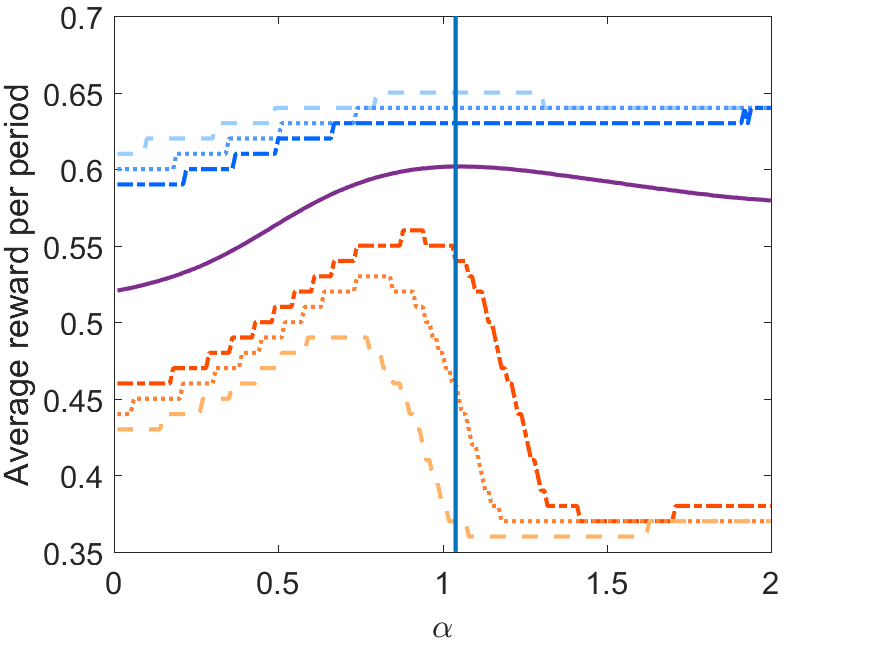}
        \caption{\footnotesize$\mu=(0.6,0.4)$, $T=100$}
    \end{subfigure}
    \begin{subfigure}{0.23\textwidth}
        \centering
        \includegraphics[width=\textwidth]{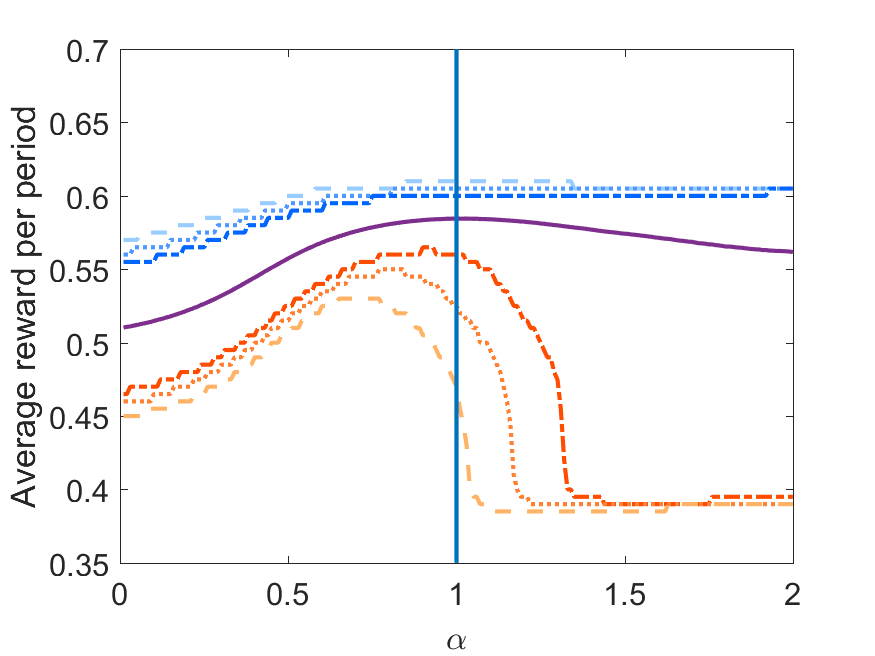}
        \caption{\footnotesize$\mu=(0.6, 0.4)$, $T=200$}
    \end{subfigure}
    \begin{subfigure}{0.23\textwidth}
        \centering
        \includegraphics[width=\textwidth]{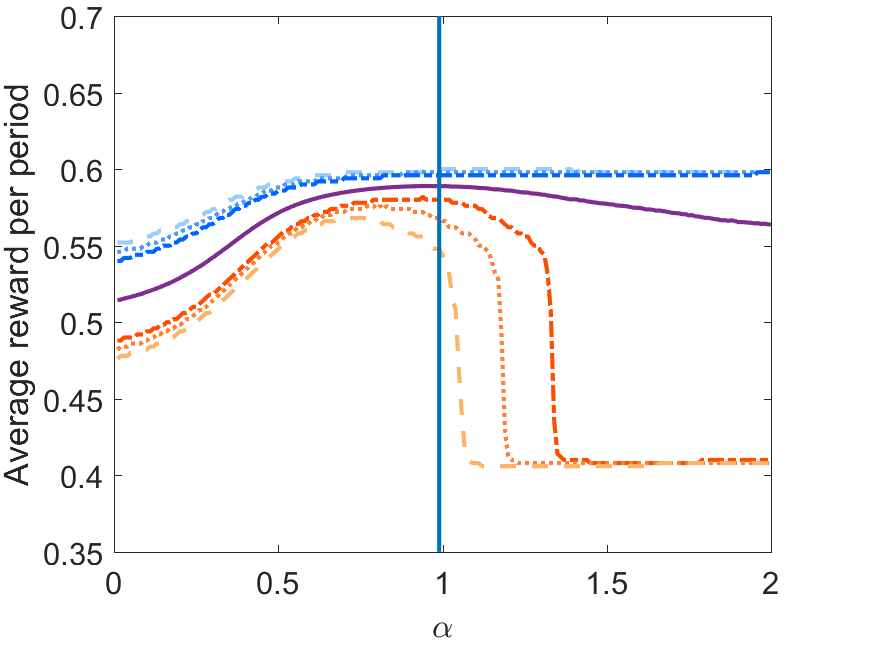}
        \caption{\footnotesize$\mu=(0.6,0.4)$, $T=500$}
    \end{subfigure}
    \begin{subfigure}{0.24\textwidth}
        \centering
        \includegraphics[width=\textwidth]{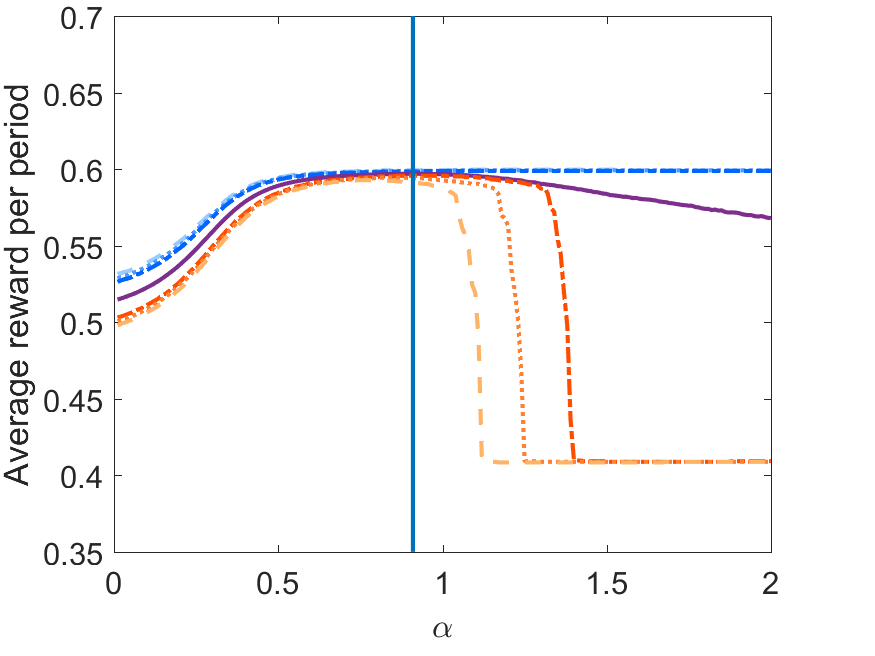}
        \caption{\footnotesize$\mu=(0.6,0.4)$, $T=2500$}
    \end{subfigure}
	\hspace{-0.35 cm}
    \begin{subfigure}{0.06\textwidth}
        \centering
         \raisebox{0.8cm}{%
        \includegraphics[width=\linewidth]{figs/legend2.png}%
    }
    \end{subfigure}

    \caption{Reward Distribution of QCARE For Finite $ T $}
    \label{fig:finite-reward}
\end{figure}

So far, our asymptotic analysis of QCARE explains how $ \alpha $ affects decision qualities when $ T = \infty $. In this subsection, we test whether our insights from the asymptotic theory can carry over to the finite-horizon setting, i.e., $ T < \infty $. To do so, we simulate the reward distributions of QCARE with a range of reward parameters ($ \mu \in  \{(0.4, 0.5), \, (0.6,0.4)\}$),  $ \alpha $ values ($ \alpha \in  \{0.01,0.02,\ldots,1.99,2\} $), and time horizons ($ T \in \{100,200, 500,2500\} $). Each simulation is based on 100,000 sample paths. We present our results in Figure~\ref{fig:finite-reward}. In this figure, the purple curve represents the average rewards generated by QCARE with the corresponding $\alpha$ values. The blue dashed curves depict the 95th, 97.5th, and 99th percentiles of the reward distribution, and the orange dashed curves represent the 1st, 2.5th, and 5th percentiles. The vertical line represents the ``optimal'' $ \alpha $ for every given configuration, which we obtain through a grid search.

Let us summarize our observations as we move from the asymptotic setting $ (T = \infty) $ to finite horizons $ (T < \infty) $. First, the mean reward is unimodal: it first increases and then decreases in $ \alpha $. As a result, there is still an ``optimal'' level of $ \alpha $ to maximize the expected reward. Once again, it represents the \textit{sweet spot} for the EE trade-off and justifies the notions of over- and under-exploration measured by $ \alpha $. The main difference here is that the optimal $ \alpha $ is generally larger than $ 0.5 $. That is intuitive since as the time horizon gets shorter, there is less time to explore and the decision maker should exploit more aggressively. This observation is consistent with the literature, which finds that although TS (i.e., $ \alpha = 0.5 $) is an order optimal policy, it tends to over-explore in the finite horizon; see \cite{min2019thompson,jin2023thompson} for related work. In general, while there is a simple answer regarding how much to explore in the infinite horizon, the optimal $ \alpha $ in the finite time horizon can be configuration-specific. That is quite expected since while the infinite-horizon setting tends to focus on the overall landscape of the EE trade-off and makes many factors negligible, those factors become material in the finite-horizon setting.

Another interesting finding concerns how the reward \textit{distribution} (e.g., its tail) is influenced by $ \alpha $. When $ \alpha $ is larger than the optimal level (and hence in the under-exploration region), both the dispersion and tail risk of the reward are significantly increased. That is because of the non-negligible probability of the extremely detrimental event where the decision maker is stuck in the inferior arm. In comparison, when $ \alpha $ is in the over-exploration region, there can be a noticeable drop in the mean reward for every fixed $ T $, especially when $ \alpha $ and $ T $ are small. However, both the dispersion and tail of the reward distributions are milder, and the drop in the mean reward can be quickly mitigated as $ T $ grows. Our observations are consistent with our insights from the previous analysis, and we revisit these observations when analyzing the behavioral patterns in the experimental data.

\section{Empirical Validation and Implications}\label{sec:structure-estimation-behavioral-data}

Empirically, the value of $\alpha$ is a behavioral parameter that can be estimated from data. In this section, we validate QCARE from various perspectives. First, we use synthetic data to demonstrate its ability to recover the ground truth parameter values. Then, we use real data to highlight QCARE's favorable capability to capture human behavior, especially in terms of out-of-sample prediction power.
We also illustrate how the QCARE-generated reward distributions align well with the behavioral data after controlling for $ \alpha $. Finally, we discuss a behavioral pattern we discover: people tend to ``over-explore'' compared to the optimal level that maximizes the expected reward.
\subsection{Empirical Analysis}\label{subsec:empirical analysis}

\vspace{0.2 cm}
\noindent \textbf{Estimation.} We first explain how we estimate $\alpha$ through maximum likelihood estimation (MLE) using the behavioral data collected in the aforementioned experiments. Since we collected the human behavioral data in two-arm MAB settings, we present the MLE procedure for $N=2$; the case with $N>2$ generalizes accordingly. To capture heterogeneity across subjects, we estimate the model parameter $\alpha$ separately for each individual.

Now, suppose we have gathered data for a specific subject denoted as $j$. It can be expressed in the form of $ \{\hat{\mu}^j_i(t), k^j_i(t), y_{ti}^{j}\}_{ijt} $, which we also refer to as a \textit{sample path}. Here,  $ \hat{\mu}^j_i(t)$ and $ k^j_i(t) $ are the state variables that subject $ j $ faces for arm $ i $ at period $ t $, and $ y_{ti}^{j} $ is an indicator variable representing whether arm $ i $ is pulled by subject $ j $ at period $ t $. Recall that $\epsilon_{it}^{j}$  are i.i.d.\ draws from the standard normal distribution. We denote the cumulative distribution function of a standard normal random variable by $\Phi(\cdot)$. Thus, the log-likelihood function for the parameter $ \alpha $ under QCARE is given by 
\begin{align*} 
	LL^j(\alpha) =\  &\sum_t \sum_{i\in \{1,2\}} y_{ti}^{j} \log \left(\Pr { \hat{\mu}^j_i(t)+ \frac{\epsilon_{it}^{j}}{(k^j_i(t)+1)^{\alpha}} \geq \hat{\mu}^j_{3-i}(t)+ \frac{\epsilon_{(3-i)t}^{j}}{(k^j_{3-i}(t)+1)^{\alpha}} } \right) \\
	=\  & \sum_t \sum_{i\in \{1,2\}} y_{ti}^{j}  \log \left( \Phi \left( \tfrac{\hat{\mu}^j_i(t) - \hat{\mu}^j_{3-i}(t)}{\sqrt{1/\left(k^j_i(t)+1 \right)^{2\alpha} + 1/\left(k^j_{3-i}(t)+1 \right)^{2\alpha}}} \right) \right).
\end{align*}
The estimated value of $ \alpha $ is thus calculated by solving $ \max_{\alpha} LL^j(\alpha) $.

To check whether the estimated $\alpha$ recovers the ground-truth value, we simulate the behavioral data according to QCARE for
various $ (T, \alpha) $ combinations, each with 10,000 sample paths. We estimate a separate $ \hat{\alpha} $  for every single sample path. In Table~\ref{tb:comp-alpha-estimation}, we report the distributions of $ \hat{\alpha} $ for every combination of $ (T, \alpha) $ by showing their average values and their 5th and 95th percentiles. As this table demonstrates, the ground truth $\alpha$ is always within the 5th and 95th percentiles of $\hat{\alpha}$, and the average $\hat{\alpha}$ is very close to the ground truth $\alpha$.

\begin{table}[htbp]
	\caption{Comparison of Ground-Truth and Estimated $\alpha$}
	\label{tb:comp-alpha-estimation}
	\centering{
		\begin{tabular}{cc|ccc|ccc}
			\toprule
		 &	 & \multicolumn{3}{c}{ $ \mu = (0.4, 0.5) $}& \multicolumn{3}{c}{$ \mu = (0.6,0.4) $}  \\
		 \midrule
   $ T $ & $\alpha$& average $\hat{\alpha}$ & \makecell{5th percentile \\of $\hat{\alpha}$} & \makecell{95th percentile \\of $\hat{\alpha}$} & average $\hat{\alpha}$ & \makecell{5th percentile \\of $\hat{\alpha}$} & \makecell{95th percentile \\of $\hat{\alpha}$} \\
			\midrule
   100   & 0.2 & 0.2015 & 0.0139 & 0.4881 & 0.1951 & 0.0139 & 0.4181 \\
200   & 0.2 & 0.1937 & 0.0139 & 0.3948 & 0.1924 & 0.0139 & 0.3171 \\
1000  & 0.2 & 0.1967 & 0.1227 & 0.2549 & 0.1999 & 0.1694 & 0.2315 \\
10000 & 0.2 & 0.2004 & 0.1927 & 0.2082 & 0.2005 & 0.1927 & 0.2082 \\
100   & 0.5 & 0.4827 & 0.2082 & 0.6591 & 0.4996 & 0.3637 & 0.6280 \\
200   & 0.5 & 0.4949 & 0.3948 & 0.5813 & 0.5014 & 0.4414 & 0.5658 \\
1000  & 0.5 & 0.5006 & 0.4803 & 0.5192 & 0.5011 & 0.4803 & 0.5192 \\
10000 & 0.5 & 0.5008 & 0.4958 & 0.5036 & 0.5005 & 0.4958 & 0.5036 \\
100   & 1   & 1.0388 & 0.8923 & 1.2421 & 1.0531 & 0.8845 & 1.3276 \\
200   & 1   & 1.0258 & 0.9234 & 1.1644 & 1.0371 & 0.9078 & 1.2265 \\
1000  & 1   & 1.0115 & 0.9467 & 1.0866 & 1.0216 & 0.9312 & 1.1371 \\
10000 & 1   & 1.0067 & 0.9622 & 1.0555 & 1.0138 & 0.9545 & 1.0983 \\
			\bottomrule
	\end{tabular}}
\end{table}

\noindent \textbf{Robustness to model misspecification.} Although QCARE can approximate a large family of policies in the long run, it is subject to model misspecification for finite $ T $, as is the case for any parametric model. It is thus valuable to assess whether it remains informative about the exploration-exploitation trade-off when the data-generating policy deviates from baseline QCARE.

We evaluate this using another simulation study. We consider a two-armed Bernoulli bandit with horizon \(T=200\) and reward means \((0.6,0.4)\), and we simulate 500 subjects for each policy and each parameter value. We consider three data-generating policies, in order of increasing departure from QCARE.

\begin{enumerate}
	\item \textit{QCARE with a different noise scale.} The first policy is QCARE itself,
	except that the Gaussian shocks have variance \(\sigma^2=4\). We write \(\alpha_0\) for its ground-truth reduction rate of exploration. As in QCARE, a smaller \(\alpha_0\) represents a stronger tendency to explore.
	
	\item \textit{Reduced-form exploration-reduction policy.} Inspired by the reduced-form representation~\eqref{eq:reduced form}, the second policy preserves the reduced-form representation above but replaces the intercept \(O(1)\) term with a fixed constant. In each period, it selects the non-leading arm with probability
		$
		p_t^{\mathrm{non\text{-}lead}}
		=\min\Bigl\{0.5,\ \exp\bigl(-c\,\min\{k_1(t)+1,k_2(t)+1\}^{2\alpha_0}\bigr)\Bigr\},
		$
		so that for sufficiently large $ k_1 $ and $ k_2 $,
		$$
		\tfrac{1}{2}\log\log\!\left(\tfrac{1}{p_t^{\mathrm{non\text{-}lead}}}\right)
		=\alpha_0\log \min\{k_1(t)+1,k_2(t)+1\}+\tfrac{1}{2}\log c.
		$$
		In the simulation, we take  \(c=0.05\). As before, a smaller \(\alpha_0\) corresponds to a stronger
		tendency to explore.

	\item \textit{Decaying-\(\varepsilon\)-greedy policy.} The third policy is a decaying-\(\varepsilon\)-greedy rule with \(\varepsilon_t=t^{-\beta}\). Here, exploration reduction is arm-agnostic, and a smaller \(\beta\) means slower reduction of exploration. Note that this policy does not necessarily satisfy the reduced-form representation \eqref{eq:reduced form}.
	
\end{enumerate}

\begin{figure}[htbp]
	\centering
	\includegraphics[width=\linewidth]{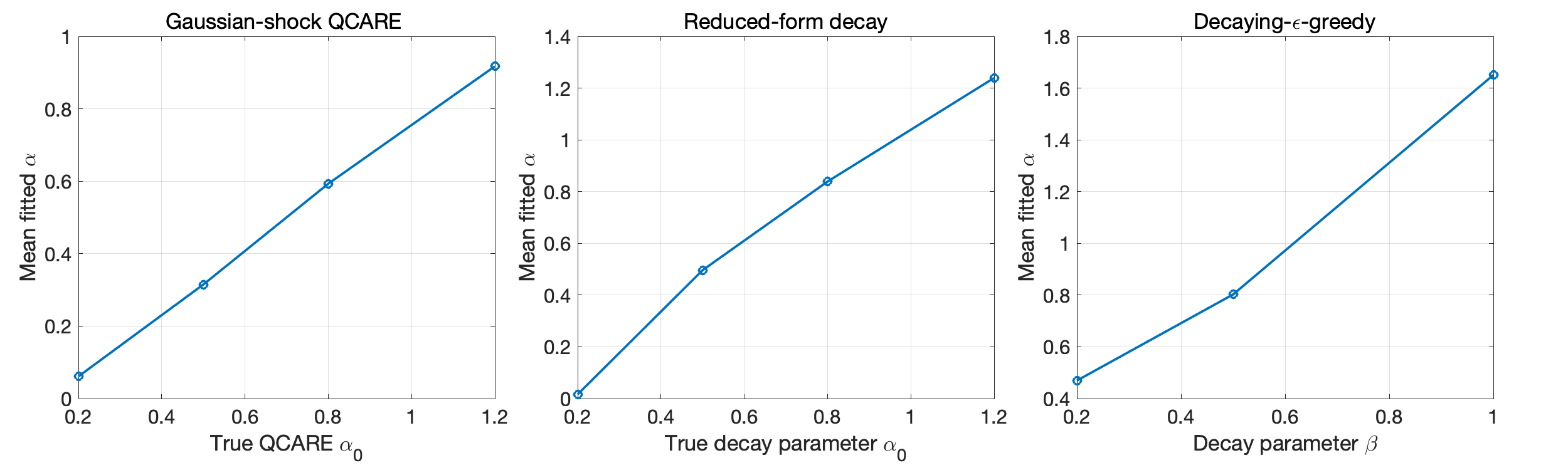}
	\caption{Estimated $ \alpha $ under Model Misspecification}
	\noindent\parbox{\linewidth}{\sf \scriptsize
		\textit{Note:} The left panel reports QCARE with noise scale $\sigma=2$; the middle panel reports the reduced-form exploration-reduction policy; and the right panel reports the decaying-$\varepsilon$-greedy policy.}
	\label{fig:misspecification-order-preservation}
\end{figure}

For every simulated path, we fit QCARE and compare the mean fitted $\hat{\alpha}$ with the ground truth exploration parameter. Figure~\ref{fig:misspecification-order-preservation} shows that the mean fitted $\hat{\alpha}$ increases monotonically with that parameter under all three data-generating policies. Thus, the estimated $\alpha$ can remain informative about relative exploration tendencies even when the data-generating policy departs from baseline QCARE. In this sense, the qualitative patterns of the EE trade-off are all correctly diagnosed.

\vspace{0.2 cm}
\noindent \textbf{Validation on real behavioral data.} Apart from simulations, we conducted a comprehensive study on how well QCARE performs in terms of explanation/prediction power in behavioral data from the experiments. 

We draw from different streams of literature and consider seven MAB policies as benchmarks. In the behavioral domain, a representative work is that by \citet{GansKnoxCroson2007}, who also investigated human behavior in a two-arm MAB problem in laboratory experiments.
Accordingly, we include the behavioral models from their study: \textit{Last-$ n $}, \textit{hot hand}, and \textit{exponential smoothing}. Under those models, the probability of selecting an arm for subject $ j $ is determined through a softmax formula:
\begin{equation}\label{eq:Gans-Prob}
P_j\{\text{choose arm $i$ at period $ t $}\} = \frac{e^{\beta_j I_{jt}^i}}{\sum_{l \in [N] } e^{\beta_j I_{jt}^l}},
\end{equation}
where $\beta_j$ is a parameter that needs to be estimated for subject $j$, and $I_{jt}^i$ is an index for subject $j$ and arm $i$ at time $t$. Each model corresponds to a different formula for $I_{jt}^i$, which we 
explain below.
\begin{itemize}
\item \textbf{Last-$ n $ model (L-$n_j$):}
This model operates on the assumption that the subject $j$ possesses limited memory of their past trials. The index $I_{jt}^i$  is determined by averaging the rewards obtained from the previous $n_j$ pulls of arm $ i $.
Here $n_j$ is a behavioral parameter to be estimated for every subject $ j $.

\item  \textbf{Hot hand model (HH-$n_j$):} This model operates on the ``hot hand'' concept, that is, once an arm has won recently, it is more likely to win again. Here the index $I_{jt}^i$ is constructed in the following way.
For arm $ i $ sampled at period $ t-1 $, $I_{jt}^i$ is set to be $ 0 $ if subject $j$ encounters consecutive losses on arm $ i $ within the last $n_j$ pulls, and $ 1 $ otherwise. For other arms, the index equals $ 1- I_{jt}^i$.
Here $n_j$ is a parameter to be estimated for every subject $ j $.

\item \textbf{Exponential smoothing model (ES-$\gamma_j$):} In this model, the index of arm $i$ is based on a time-discounted average of historical rewards. Specifically, let $t^{-}_i$ denote the last time arm $ i $ was pulled before period $ t $. Then we set $I^i_{jt} = \gamma_j r_i(t^-_i) + (1-\gamma_j)I^i_{j t^{-}_i}$.
Here the parameter $ \gamma_j  $ takes values in $ (0,1) $ and is to be estimated for every subject $ j $.
\end{itemize}
To investigate whether we can achieve the same descriptive/predictive power by a ``simpler'' model than QCARE, we consider a variation of QCARE.
\begin{itemize}
\item \textbf{QCARE-0:}  We adopt QCARE but restrict $ \alpha = 0 $. In this way, the exploration term does not shrink over time.
\end{itemize}
Finally, we include the following classical bandit algorithms from the online learning literature.
\begin{itemize}
\item \textbf{Fixed-$\varepsilon_j$-greedy:} In this model, at time period $t$, subject $j$ chooses the leading arm with probability $1-\varepsilon_j$ and each other arm with probability $\frac{\varepsilon_j}{N-1}$.
The value of $\varepsilon_j$ is to be estimated for each subject.

\item \textbf{Thompson Sampling (TS-B, TS-G):} We consider two types of Thompson Sampling algorithms: one with a Beta prior (TS-B) and the other with a Gaussian prior (TS-G). In these models, the subjects follow the Thompson Sampling approach with respective prior distributions; see Appendix \ref{appendix:sec:algorihims} for a more detailed algorithmic description.
\end{itemize}
Whenever applicable, the model parameters are estimated through MLE.

We compare both the in-sample and out-of-sample performances of all models. To accommodate the heterogeneity among individuals, we split the data into training and testing sets along the temporal dimension. That is, we use the first three quarters of periods to train/calibrate models for each subject, and we use the last quarter of periods for out-of-sample evaluation. We use two evaluation metrics: (i) the average log-likelihood values and (ii) head-to-head comparisons in terms of log-likelihood between QCARE and benchmark models, where we report the fraction of subjects for whom QCARE outperforms the benchmark. We also investigate whether the performance differences are statistically meaningful.\footnote{To do so, we employ a one-tailed paired $t$-test for the average log-likelihood comparison ($ H_0 $: the average log-likelihood of QCARE equals that of the benchmark; $ H_1 $: the average log-likelihood of QCARE is higher than that of the benchmark) and a one-sample proportion test for the head-to-head comparison ($ H_0 $: the proportion of subjects for whom QCARE outperforms the benchmark equals 50\%; $ H_1 $: this proportion exceeds 50\%).} To summarize, we have considered multiple factors such as various MAB policies, the performance metric, and the data to evaluate, and included $ 8 \times 2 \times 2 = 32 $ combinations of configurations in the comparison.

Table~\ref{tb:comp-predict-power} reports the results. QCARE dominates QCARE-0, TS-B, and TS-G in terms of both in-sample explanatory power and out-of-sample predictive power. Relative to $\varepsilon_j$-greedy, QCARE predicts better out of sample for a significant majority of subjects. The three heuristic models, L-$n_j$, HH-$n_j$, and ES-$\gamma_j$, fit the training data better than QCARE for a majority of subjects, which is unsurprising because each has two parameters to estimate, whereas QCARE has only one. Out of sample, however, QCARE outperforms L-$n_j$ and ES-$\gamma_j$ for a significant majority of subjects and achieves a significantly higher average log-likelihood than both. Against HH-$n_j$, the strongest heuristic benchmark, QCARE performs comparably out of sample. Taken together, these results suggest that QCARE provides favorable prediction power for human choices out of sample.\footnote{As a robustness check, we repeat the out-of-sample evaluation under alternative temporal train--test splits. The main comparative conclusions remain unchanged.}

\begin{table}[htbp]
	\caption{ Comparison Between QCARE and Other Models}
	\label{tb:comp-predict-power}
	\centering
			\begin{tabular}{cccccccccc}
			\toprule
			\footnotesize{Data} & \footnotesize{Metric}  & QCARE   & L-$n_j$ & HH-$n_j$ & ES-$\gamma_j$ & QCARE-0 & $\varepsilon_j$-greedy & TS-B  & TS-G \\
			\midrule
			 \multirow{2}{*}{IS}  & AvgLL  & -64.6      & -65.0     & -61.1    & -63.1     & -72.1\dstar    & -64.9  & -93.9\dstar & -71.3\dstar  \\
			 & Pct & --   & 42.3\% & 32.5\% & 33.3\% & 100.0\%\dstar & 52.8\% & 97.6\%\dstar & 98.4\%\dstar \\
			 \cmidrule{1-10}
			\multirow{2}{*}{OoS} & AvgLL   & -20.6    & -25.7\dstar    & -20.4     & -25.3\dstar    & -23.7\dstar      & -20.5  & -28.5\dstar & -22.5\dstar    \\
			& Pct & -- & 61.0\%\dstar & 52.8\% & 63.4\%\dstar & 84.6\%\dstar & 62.6\%\dstar & 78.0\%\dstar & 60.2\%\star \\
			\bottomrule
	\end{tabular}
	
	\vskip 0.1 cm
	\noindent\parbox{\linewidth}{\sf \scriptsize
		\textit{Note:} ``IS'' means in-sample (i.e., evaluated on the training data), and ``OoS'' means out-of-sample (i.e., evaluated on the test data). ``AvgLL'' refers to the average log-likelihood at the subject level by each approach. ``Pct'' indicates the corresponding percentage of subjects for whom QCARE outperforms the benchmark. The remarks ``**'' and ``*'' indicate $ p $-values less than $ 0.01 $ and $ 0.05 $, respectively.}
\end{table}

\begin{figure}
	\centering
	\begin{subfigure}{0.41\textwidth}
		\includegraphics[width=\textwidth]{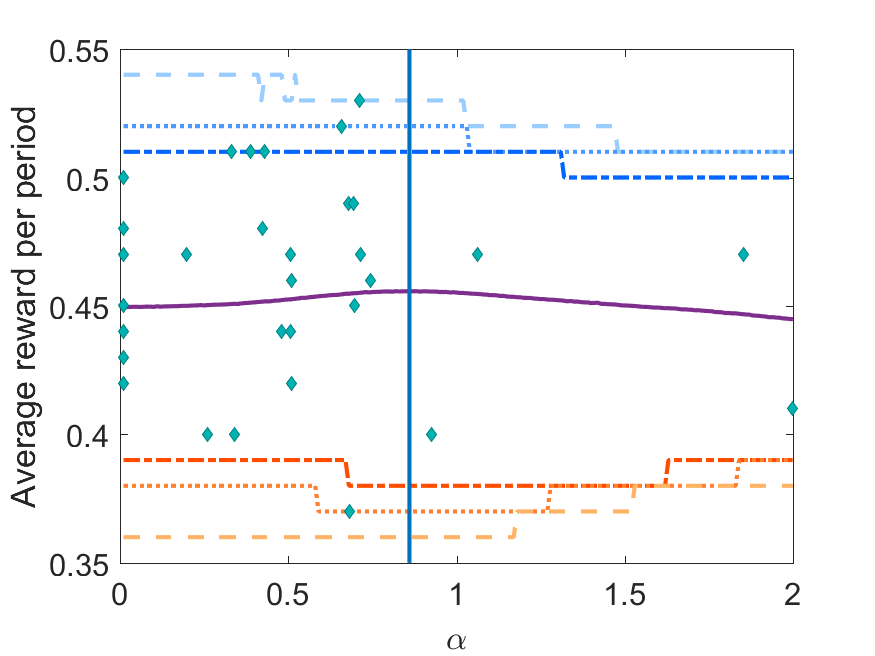}
		\subcaption{$\mu = (0.4, 0.5)$, $T=100$}
	\end{subfigure}
	\hspace{0.4 cm}
	\begin{subfigure}{0.41\textwidth}
		\includegraphics[width=\textwidth]{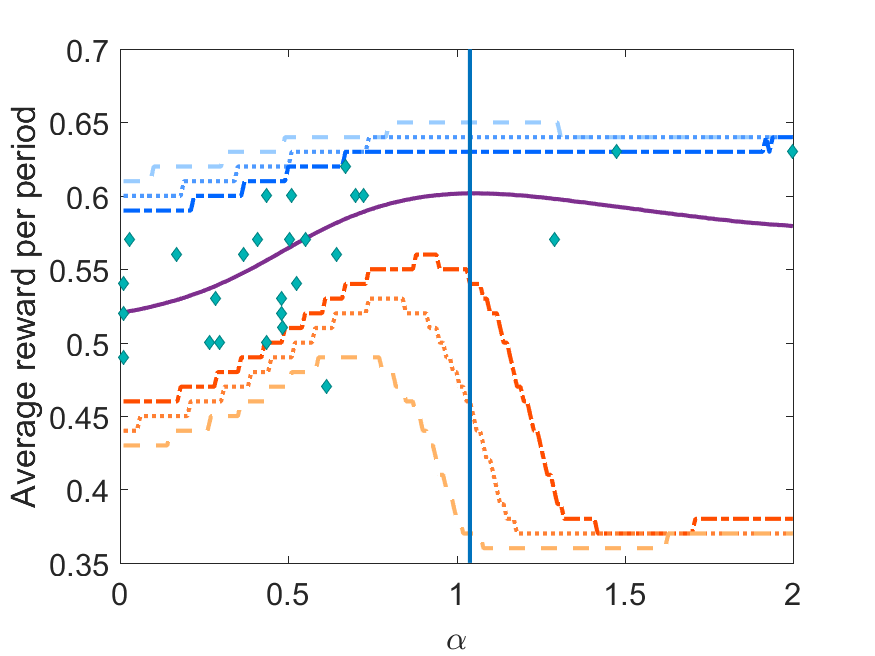}
		\subcaption{$\mu = (0.6,0.4)$, $T=100$}
	\end{subfigure}
	\begin{subfigure}{0.09\textwidth}
		\includegraphics[width=\textwidth]{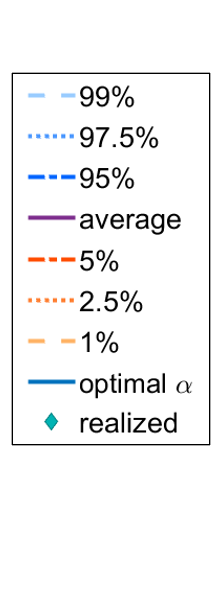}
	\end{subfigure}
	
	\begin{subfigure}{0.41\textwidth}
		\includegraphics[width=\textwidth]{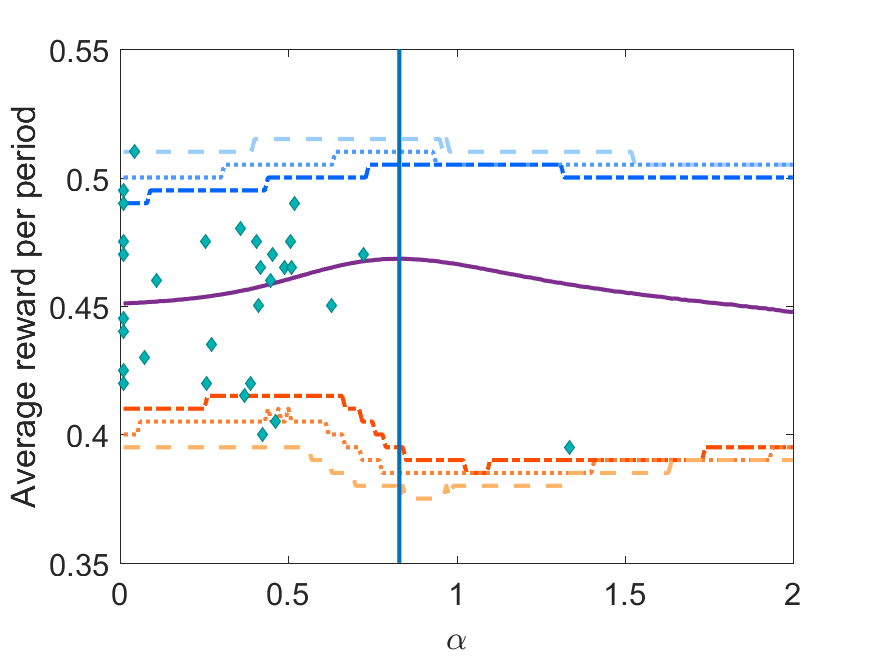}
		\subcaption{$\mu = (0.4, 0.5)$, $T=200$}
	\end{subfigure}
	\hspace{0.4 cm}
	\begin{subfigure}{0.41\textwidth}
		\includegraphics[width=\textwidth]{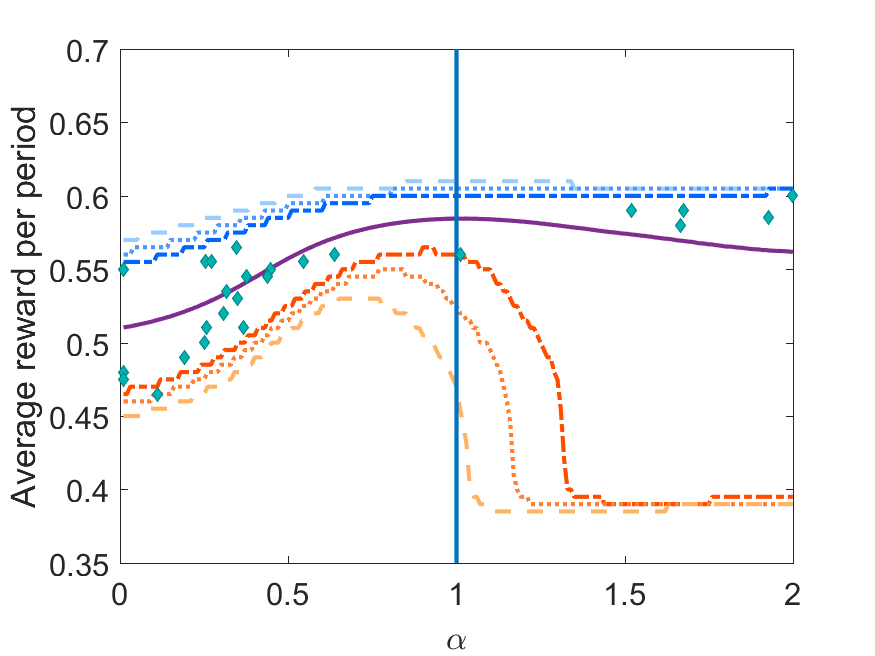}
		\subcaption{$\mu = (0.6,0.4)$, $T=200$}
	\end{subfigure}
	\begin{subfigure}{0.09\textwidth}
		\includegraphics[width=\textwidth]{figs/legend-beha.png}
	\end{subfigure}
	\caption{Comparison of Estimated \(\alpha\) by Individual Subjects}
	\noindent\parbox{\linewidth}{\sf \scriptsize
		\textit{Note:} The green dots correspond to subjects from the baseline experiment, with 123 subjects in total.}
	\label{fig:beha-individual}
\end{figure}

\vspace{0.2 cm}
\noindent \textbf{Implied alpha and reward distributions.} Finally, we use our QCARE model to address our original question: how do human beings \textit{actually} balance the EE trade-off? Figure~\ref{fig:beha-individual} summarizes our findings. In this figure, each green dot represents an individual subject and contains information about the subject's implied $ \hat{\alpha} $ value and realized reward. The green dots are then organized by experiment configurations (i.e., mean reward vector $ \mu $ and time horizon $ T $). For each configuration, we also follow the same method as in Figure \ref{fig:finite-reward} to display the reward distributions generated by a range of $ \alpha $ values. That is, we report the mean rewards and the 1st, 2.5th, 5th, 95th, 97.5th, and 99th percentiles of the reward distribution based on 10,000 simulated paths. Finally, the instance-specific optimal $ \alpha $ value for each configuration is represented by the blue vertical line.
Comparing the realized and simulated data, we find that the majority of the green points fall within the percentile bands under QCARE, further validating it as a reasonable model to capture human behavior in the MAB problem.

A perhaps even more interesting (and surprising) finding from Figure \ref{fig:beha-individual} is that a large majority of the green points lie to the left of the vertical line that represents the expected-reward-maximizing $ \alpha $. In fact, a significant portion of subjects have implied $ \alpha $ values less than $ 0.5 $, the optimal level in the asymptotic regime when $ T = \infty $.
It suggests that most subjects tend to exhibit an \textit{over-exploration} behavior. We find this pattern surprising since if every subject has their own idiosyncratic behavioral preferences/biases, one would expect the distribution of $ \alpha $ to be more balanced around the optimal level. 

\vspace{0.2 cm}

The behavioral pattern towards over-exploration in Figure \ref{fig:beha-individual} is robust across the configurations, such as time horizon and reward distributions. The results suggest that while those factors could potentially influence the EE trade-off, which is reflected in the different values of the optimal $ \alpha $, the fact that people tend to explore more aggressively than the reward-maximizing level seems consistent. The patterns are also robust across the various replications and auxiliary experiments as reported in Appendix~\ref{section:appendix:auxiliary experiments}. In particular, we conduct targeted robustness checks in which subjects receive five practice rounds before the task, face higher monetary stakes, or are reminded of the remaining periods near the end of the task. These conditions do not change the tendency toward over-exploration. Details are provided in Appendix~\ref{appendix:additional-experimental-checks}.

\subsection{Further Discussion}\label{subsec:further discussion}

\noindent \textbf{Prescriptive implications.}   The exploration-exploitation trade-off is a fundamental challenge in sequential decision making under uncertainty because each action affects both current performance and the information available for future decisions. Our findings show that human choices reflect this tension, in that the propensity to sample currently non-leading alternatives declines as arm-specific experience accumulates. The balance struck by human decision makers, however, can depart from that prescribed by reward-seeking algorithms. In our experiments, the predominant deviation is persistent exploration: the propensity to choose non-leading alternatives declines too slowly relative to reward-maximizing benchmarks. At first glance, this finding might appear to favor delegating such decisions entirely to bandit algorithms. Yet full delegation may be undesirable in many operational settings because frontline decision makers often possess tacit knowledge, qualitative signals, or operational constraints that are not encoded in the system's observed state. The appropriate prescription may therefore be to augment human judgment with behavioral diagnosis and timely decision support rather than to replace it.

QCARE suggests a few complementary forms of such support. First, organizations could use incentivized MAB-like calibration exercises, followed by individualized QCARE-based feedback, to educate frontline decision makers about the EE trade-off and their own behavioral tendencies. The resulting profile could show how quickly an individual's exploration declines with arm-specific experience and whether this pattern indicates persistent exploration or premature concentration relative to a task-appropriate benchmark. By making an otherwise latent tendency visible, such personalized feedback could provide a foundation for targeted reflection and deliberate behavioral adjustment. Second, QCARE could serve as a behavioral diagnostic tool within an operational decision-support system. When a decision history indicates persistent exploration, the system could ask: ``Do you possess contextual information not available to the system that makes further sampling of this alternative worthwhile, and what outcome would change your subsequent allocation?'' When the history indicates premature concentration, the system could instead flag alternatives supported by relatively limited evidence and invite the decision maker to reconsider whether the available evidence justifies their abandonment. In both cases, QCARE identifies and challenges a potentially problematic decision tendency, while the decision maker retains authority over the ultimate action.

Third, just as static choice models such as the Multinomial Logit (MNL) serve as the demand model for assortment and pricing decisions, QCARE provides the foundation for a dynamic choice model with learning, estimable from purchase histories and outcome data such as ratings. For example, customers learn about brands of uncertain quality through their own experiences. Over-exploration then implies that customers may settle on their best brand more slowly than one may expect. This affects the value of adding new products to an assortment, the need for trial-inducing promotions, and how strongly a recommender should steer each customer toward their most-preferred option.

\vspace{0.2 cm}
\noindent \textbf{Possible reasons for over-exploration.}
From the behavioral side, it would also be interesting to investigate the underlying mechanisms behind such behavioral patterns. Here, let us discuss a few plausible explanations. First, the pattern could be partially rationalized by the subjects' own goals and preferences, such as risk aversion. In fact, there is an emerging literature discussing the risk-return trade-off in MAB policies, and how more aggressive exploration could reduce the risk; see \citet{keller2019note}, \citet{zhu2020thompson}, \citet{simchi2023regret}, \citet{fan2021fragility} for related work. Our model displays the same pattern: with an $ \alpha $ slightly smaller than the optimal level, the risk can be reduced (based on various measures of risk, such as the dispersion of the reward distribution and the 1\%/2.5\%/5\% tail risks) at the expense of lower expected reward; see Figure \ref{fig:beha-individual}. In this sense, at least some subjects' behavior could be explained by their using over-exploration to trade expected rewards for less risk. 

However, risk aversion can only explain part of the observations.  When exploration is too heavy (i.e., when $ \alpha $ is too small), both the expected reward and the risk level deteriorate. In the real data, a noticeable portion of subjects could have simultaneously reduced their risk and increased their expected reward by using a larger $ \alpha $. It appears that the over-exploration pattern should also be attributed to the subjects' bounded rationality and cognitive limitations.\footnote{This is reasonable given the complexity of MAB problems. Even with advanced computational power, finding the Bayes optimal solution to the MAB problem can be challenging as $ T $ increases, due to the curse of dimensionality.} In other words, subjects' behavior of not choosing the leading arm appears to reflect a mix of random behavioral errors (in the same spirit as \citealp{Su2008}) layered on top of the ``intrinsic'' exploration required for dynamic learning.

We believe our finding opens the door to future research opportunities. For example, one could extend the QCARE model by exploring alternative parameterizations of exploration to handle different decision-making environments, such as those with more complex reward structures. Furthermore, applying QCARE to decision-making tasks in real-world scenarios with richer operational contexts could help validate the findings from the lab environments and deepen these insights. Finally, it would be interesting to see if one can use ``nudging'' to mitigate over-exploration.

\section{Conclusion}

This paper introduces a dynamic quantal choice model, QCARE, that parametrically extends Thompson Sampling to represent people's decision behavior in multi-armed bandit settings. The model quantifies the exploration-exploitation trade-off by incorporating a reduction rate of exploration (i.e., $ \alpha $) as more information becomes available. The relationship between different values of $ \alpha$ and decision quality is mathematically established through an asymptotic theory. A general Markovian-policy framework further establishes complete learning and the lower-bound mechanisms for over- and under-exploration beyond QCARE. The analysis thus helps us better understand the effects of over- and under- exploration. Based on behavioral data we collected from experiments, QCARE demonstrates a compelling ability to capture people's decision-making process: it replicates key qualitative patterns observed in subjects' behavior and also displays favorable explanatory/predictive power compared to other models. Finally, empirical analysis through QCARE offers a novel perspective on how over-exploration tends to be a behavioral pattern.

\bibliographystyle{ormsv080}
\bibliography{references}
\ECSwitch

\begin{appendices}

\section{Algorithmic Descriptions of QCARE and TS}\label{appendix:sec:algorihims}

We provide an algorithmic description of QCARE below.

\begin{mdframed}
	\captionsetup{labelformat=empty}
	\textbf{Quantal Choice with Adaptive Reduction of Exploration (QCARE)}
	\label{alg:algorithm1}
	
	\textbf{Input:} $ \alpha > 0$. 
	
	\textbf{Procedure:} Initialize with $k_i(1)=0, \hat{\mu}_i(1)=0$ for each arm $i\in [N]$. 
	
	At each time $t \in [T]$:
	\begin{enumerate}
		\item For each arm $i \in [N]$, sample $\epsilon_{it}$ independently from the $\mathcal{N} (0, 1)$ distribution. Set  $\theta_i(t) :=\hat{\mu}_i(t)+\epsilon_{it}/(k_i(t)+1)^{\alpha}$.
		\item Play arm $a(t):= \argmax_i \theta_i(t)$ and observe reward $r_{a(t)}(t)$.
		\item Update states:
		\begin{itemize}
			\item For arm $i=a(t)$, set $k_{i}(t+1)=k_{i}(t)+1$ and $\hat{\mu}_{i}(t+1)=\frac{\hat{\mu}_{i}(t)(k_{i}(t)+1)+r_{i}(t)}{k_{i}(t+1)+1}$;
			\item For other arms $i\neq a(t)$, set $k_{i}(t+1)=k_{i}(t)$ and
			$\hat{\mu}_{i}(t+1)=\hat{\mu}_{i}(t)$.
		\end{itemize}
	\end{enumerate}
\end{mdframed}

We also describe TS with Gaussian priors and TS with Beta priors, following the version in \citet{agrawal2013further}.

\begin{mdframed}
	\captionsetup{labelformat=empty}
	\textbf{Thompson Sampling using Gaussian Priors}
	\label{alg:algorithm2}

        For each arm $i \in [N]$ set $k_i=0$, $\hat{\mu}_i = 0$.
	
	At each time $t \in [T]$:
	\begin{enumerate}
		\item For each arm $i \in [N]$, sample $\theta_i(t)$ independently from the $\mathcal{N} (\hat{\mu}_i, \frac{1}{k_i+1})$ distribution. 
		\item Play arm $a(t):= \arg\max_i \theta_i(t)$ and observe reward $r_t$.
		\item Set $\hat{\mu}_{a(t)}:=\frac{\hat{\mu}_{a(t)}(k_{a(t)}+1)+r_{t}}{k_{a(t)}+2}$, $k_{a(t)}:=k_{a(t)}+1$.
	\end{enumerate}
	
\end{mdframed}

\begin{mdframed}
    \captionsetup{labelformat=empty}
		\textbf{Thompson Sampling using Beta priors}
		\label{alg:algorithm3}
		
		For each arm $i \in [N]$ set $S_i=0$, $F_i = 0$.
		
		At each time $t \in [T]$:
		\begin{enumerate}
			\item For each arm $i \in [N]$, sample $\theta_i(t)$ independently from the $\text{Beta}(S_i+1, F_i+1)$ distribution. 
			\item Play arm $a(t):= \arg\max_i \theta_i(t)$ and observe reward $r_t$.
			\item If $r_t = 1$, then $S_{a(t)}:=S_{a(t)}+1$, else $F_{a(t)}:=F_{a(t)}+1$. 
		\end{enumerate}
		
\end{mdframed}

\section{An Example for Aggregating the Choice Frequencies of the Leading Arm}
\label{section:appendix:example}

In the descriptive analysis of our experiment data, we consider a collection of metrics to evaluate subjects' tendency to choose the leading arm. Here we present a concrete example of how to construct the metrics. Suppose one subject plays a 10-round two-arm bandit game and the historical record is summarized in Table~\ref{tb: example_history} below. 
\begin{table}[htbp]
	\caption{Example: Historical Table of One Subject}
	\label{tb: example_history}
	\centering{
		\begin{tabular}{cccccc}
			\toprule
			$t$  & action & reward & $\hat{\mu}_1(t)$ & $\hat{\mu}_2(t)$ &\makecell{Leading arm\\ chosen? }  \\
			\midrule 
			1  & 2      & 0      & 0.000 & 0.000 & N \\
			2  & 2      & 1      & 0.000 & 0.000 & N \\
			3  & 1      & 1      & 0.000 & 0.333 & N \\
			4  & 1      & 1      & 0.500 & 0.333 & Y \\
			5  & 2      & 1      & 0.667 & 0.333 & N \\
			6  & 2      & 0      & 0.667 & 0.500 & N \\
			7  & 1      & 0      & 0.667 & 0.400 & Y \\
			8  & 1      & 0      & 0.500 & 0.400 & Y \\
			9  & 1      & 1      & 0.400 & 0.400 & N \\
			10 & 1      & 1      & 0.500 & 0.400 & Y \\
			\bottomrule
	\end{tabular}}
\end{table}
Let us walk through the calculations and values of our metrics.
\begin{enumerate}
	\item In periods 4, 7, 8, and 10, the subject chose the leading arm. The choice fraction of the leading arm is $40\%$.
	\item In periods 1, 2, 8, 9, 10, the empirical mean gaps are in the low group, while in periods 3-7, gaps are in the high group. The choice fraction of the leading arm is $40\%$ in the low group and also $40\%$ in the high group.
	\item In the first half (periods 1-5), the choice fraction of the leading arm is $20\%$. In the second half, the fraction is $60\%$.
\end{enumerate}
We calculate these metrics for every subject in the same way and then average them across subjects.

\section{Further Details about the Experiments}
\label{section:appendix:auxiliary experiments}

\subsection{The Laboratory (Baseline) Experiment}\label{section:appendix:baseline experiment}

\begin{figure}[ht]
 \centering
    \begin{subfigure}[b]{0.45\textwidth}
        \centering
        \includegraphics[width=\textwidth]{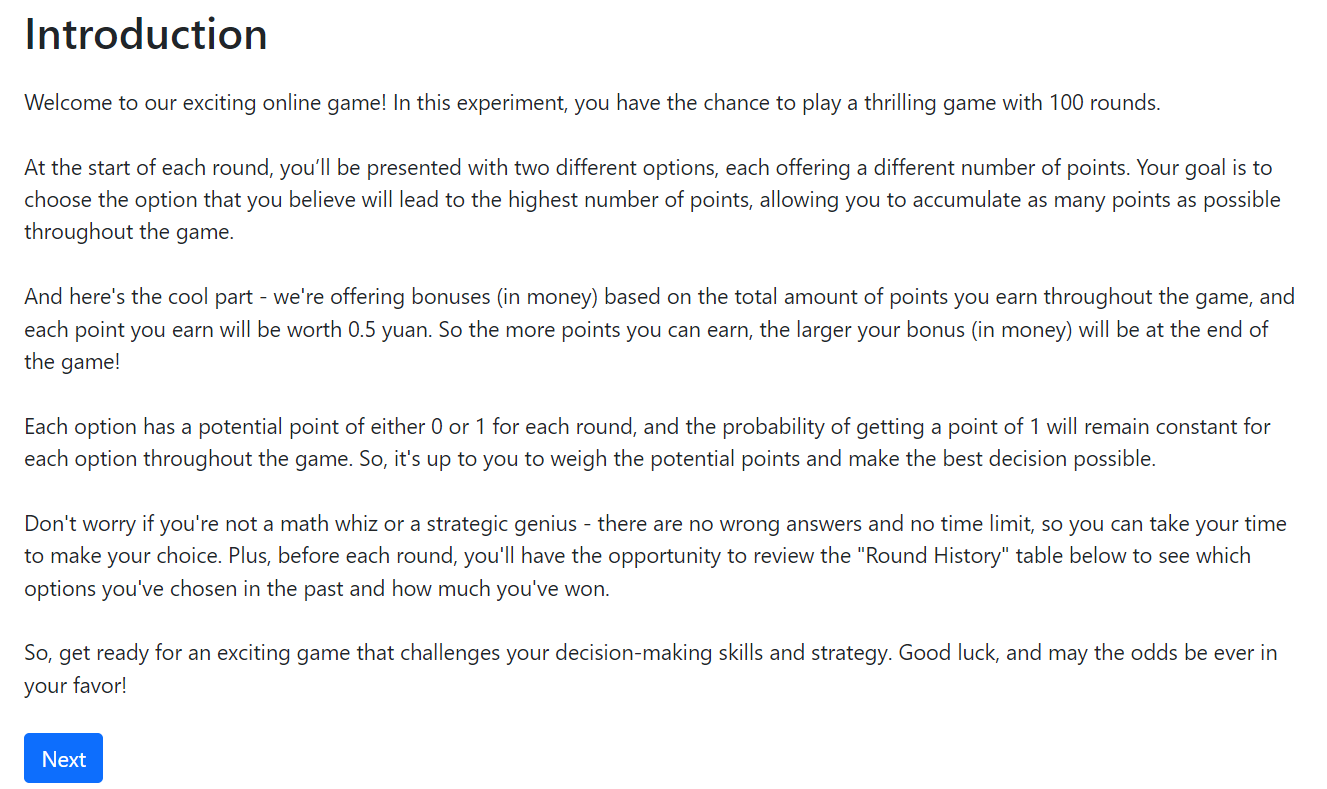}
        \caption{Introduction Page ($T = 100$)}
    \end{subfigure}
    \hspace{1cm}
    \begin{subfigure}[b]{0.45\textwidth}
        \centering
        \includegraphics[width=\textwidth]{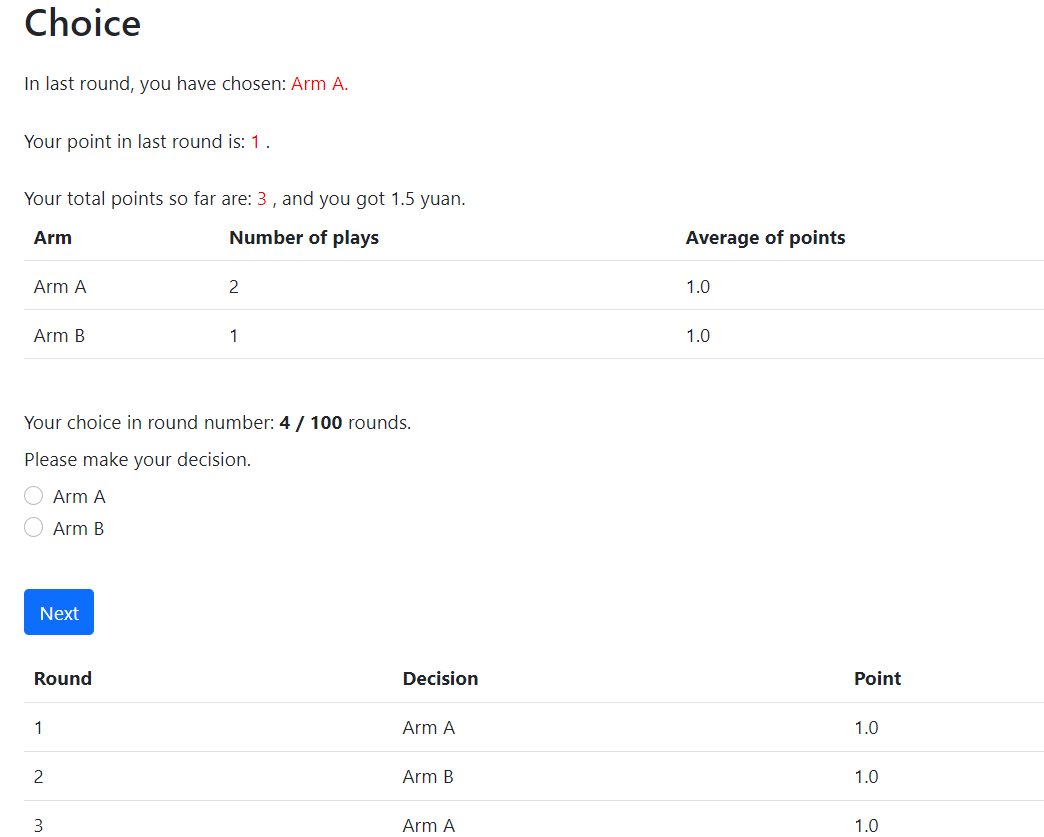}
        \caption{A Typical Round Page}
    \end{subfigure}

	\caption{Interface of the Baseline Experiment}
	\label{fig:exp-page-baseline}
\end{figure}

We present
in Figure~\ref{fig:exp-page-baseline} the experiment interfaces.

Let us also explain the payoff details. The show-up fee was $5$ Chinese Yuan (CNY). The point-to-cash conversion rule was set at \(0.5\) CNY per point for the reward-mean configuration \((0.4,0.5)\) and \(0.4\) CNY per point for the reward-mean configuration \((0.6,0.4)\). In this way, the earnings of playing the \textit{clairvoyant} policy (i.e., always selecting the arm with the larger mean) were approximately equal across parameter configurations. Based on those settings, the average total earnings per subject were approximately 5.23 US dollars (USD). On average, subjects took 13 minutes to complete the experiment, with a maximum duration of 24 minutes.

\subsection{Online Replication}
\label{section:appendix:online experiment}

We conducted an online replication facilitated by the Prolific platform to replicate our empirical analysis of the baseline experiment.

\begin{figure}
     \centering
    \begin{subfigure}[b]{0.45\textwidth}
        \centering
        \includegraphics[width=\textwidth]{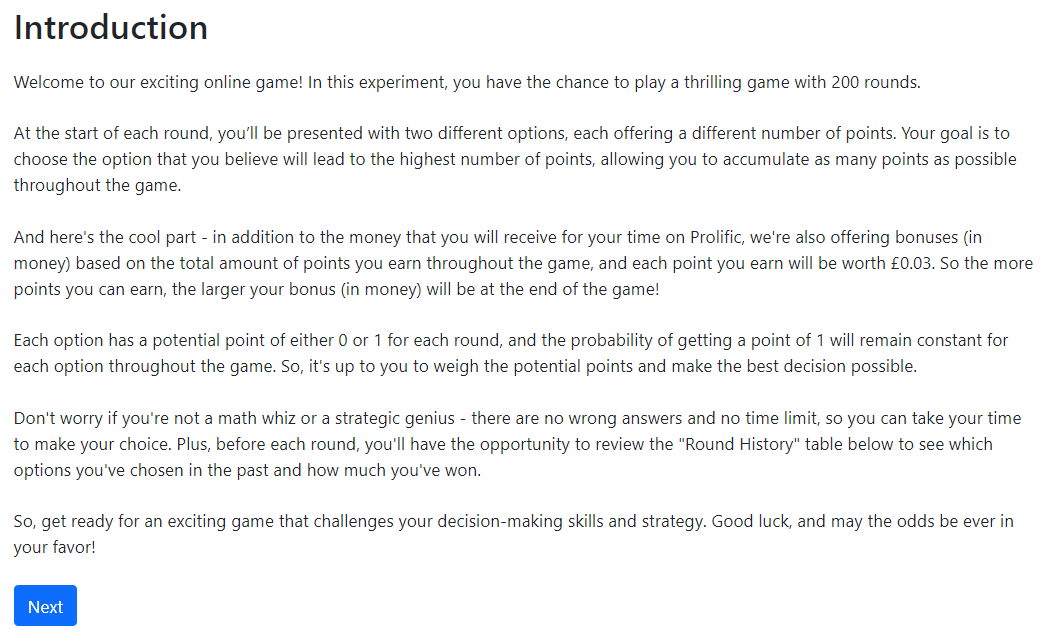}
        \caption{Introduction Page ($T = 200$)}
    \end{subfigure}
    \hspace{1cm}
    \begin{subfigure}[b]{0.45\textwidth}
        \centering
        \includegraphics[width=\textwidth]{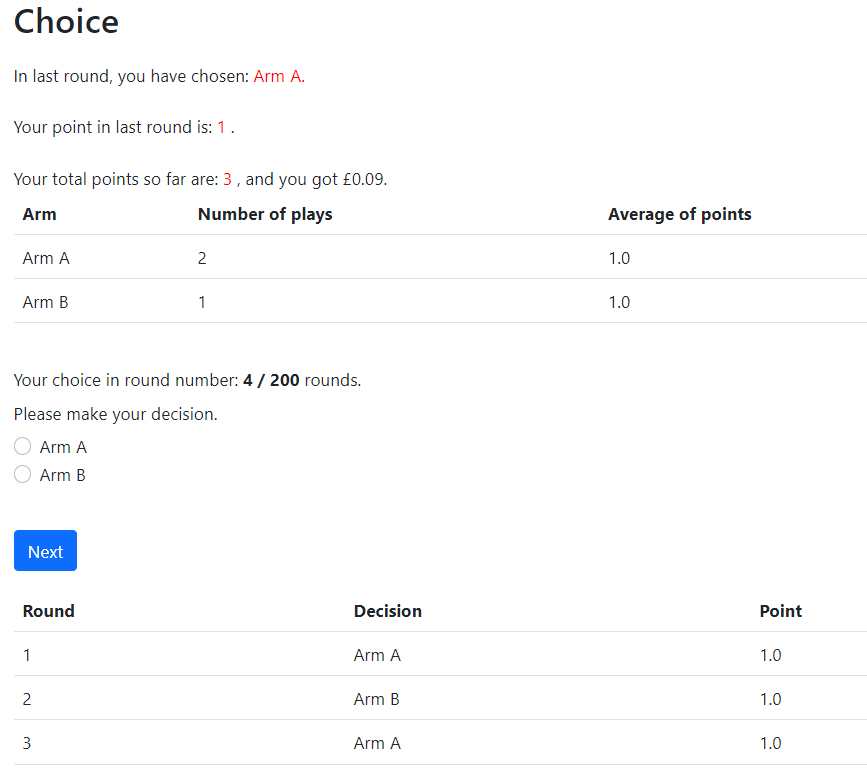}
        \caption{A Typical Round Page}
    \end{subfigure}
	\caption{Interface of the Online Replication}
	\label{fig:exp-page}
\end{figure}

The experimental procedure followed the same methodology as described in Section~\ref{subsec:experiment-design}. We show the interface of the introduction page and the round page during the experiment in Figure~\ref{fig:exp-page}. At the beginning of the experiment, the subjects were told the number of rounds and the basic information about the reward mechanism. Each participant's total earnings comprised two components: a show-up fee of 1.5 pounds and a bonus based on realized points. These points were subsequently converted into cash at the end of the session, with a conversion rate of 1 point equal to 0.03 British pounds. Based on this setting, the final average earnings per subject were 3.68 pounds (approximately 3.94 USD). On average, subjects took 14.85 minutes to complete the experiment, with a maximum duration of 36 minutes.

\noindent \textbf{Descriptive analysis.} 
Similar to the baseline experiment, the results from the online replication show that only 1 of 123 subjects chose the leading arm throughout the entire time horizon, and 87 out of 123 subjects chose the leading arm in at least 50\% of the rounds. On average, subjects chose the leading arm $59.76\%$ of the time.\footnote{If taking the average at the individual-period level, the figure is $ 59.10\% $.} That ratio is strictly less than one, yet significantly larger than $ 0.5 $ (Wilcoxon signed-rank test right tail, $ p $-value$ = $8.4440e-10). These findings suggest that human decision-making involves a balance between exploration and exploitation, rather than favoring either exploration or exploitation exclusively.

\begin{figure}[htbp]
    \centering
    \begin{subfigure}[b]{0.33\textwidth}
        \centering
\includegraphics[width=\textwidth]{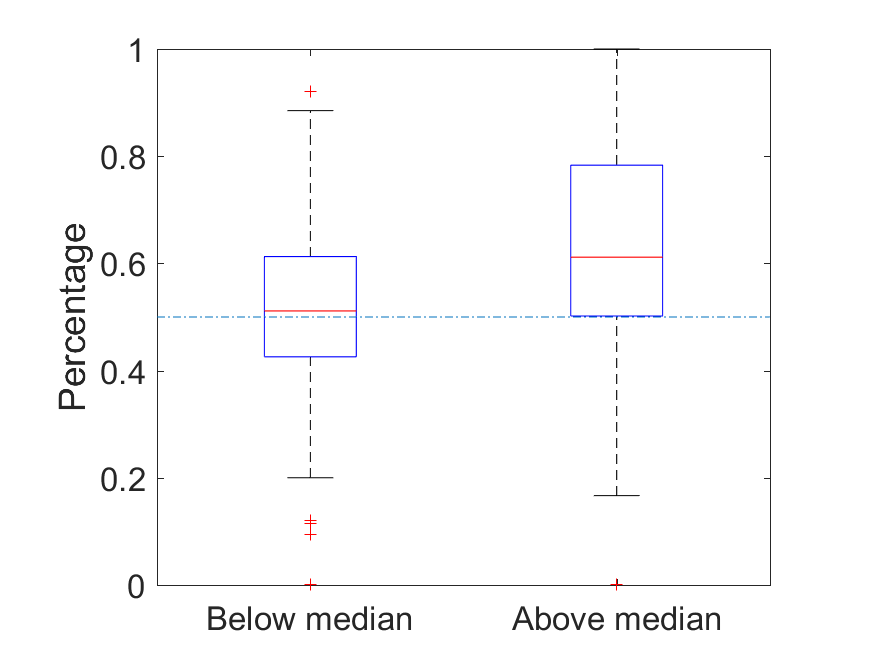}
        \caption{By empirical mean difference}
    \end{subfigure}
    \hspace{-2em}
    \begin{subfigure}[b]{0.33\textwidth}
        \centering
\includegraphics[width=\textwidth]{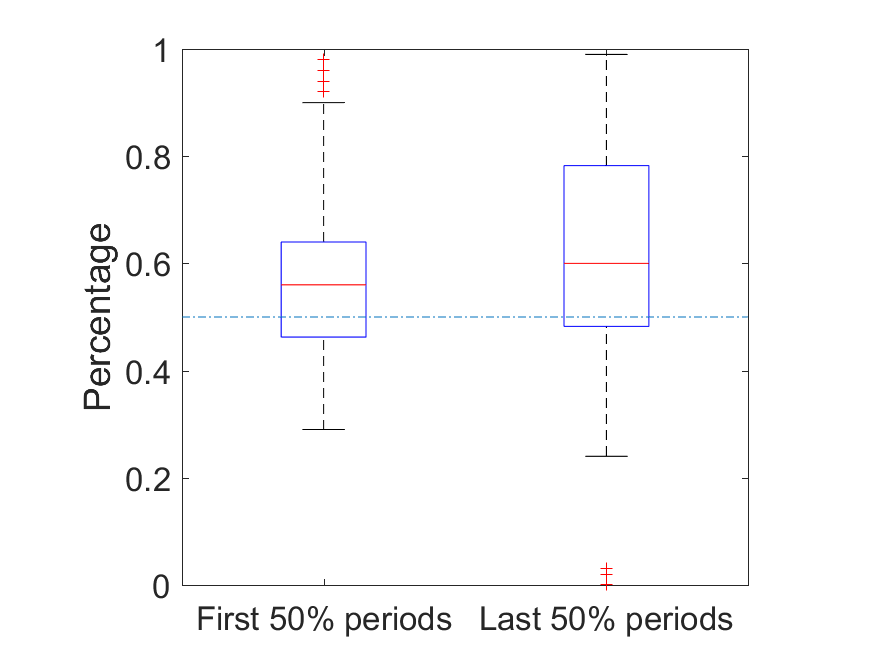}
        \caption{By period}
    \end{subfigure}
    \hspace{-2em}
    \begin{subfigure}[b]{0.33\textwidth}
        \centering
\includegraphics[width=\textwidth]{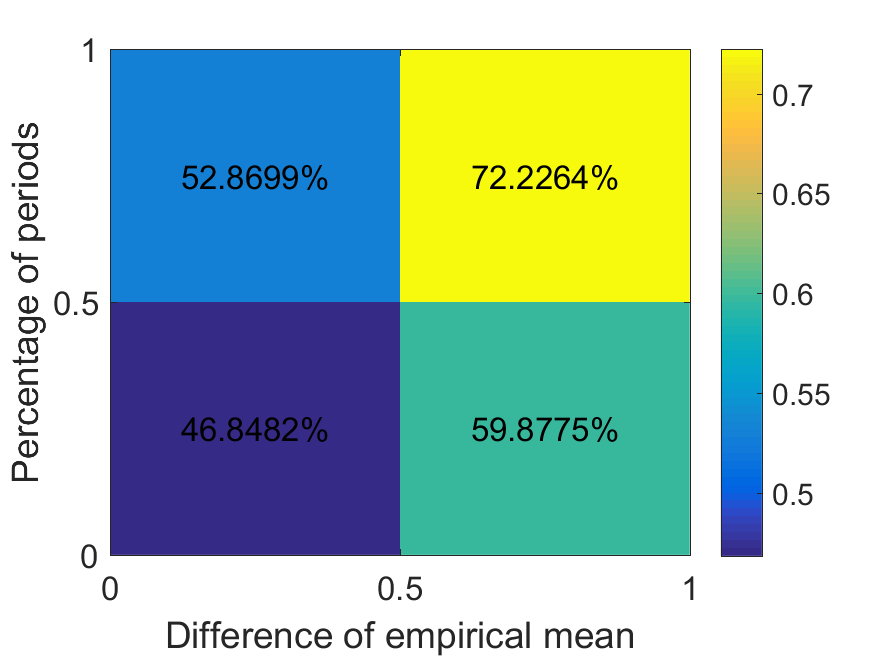}
        \caption{By both}
    \end{subfigure}
	\caption{Fraction of Periods Choosing the Leading Arm (Online Replication)}
	\label{fig:box-num-larger-2}
\end{figure}

Next, we use Figure~\ref{fig:box-num-larger-2} based on the online-replication data as a replication for Figure~\ref{fig:box-num-larger} in the main body.
In Panel (a), we observe a significantly higher average ratio when the difference is larger than the median ($ p $-value being 6.8480e-08 using a two-sided Wilcoxon signed-rank test). From Panel (b), we find a significantly higher average ratio in the second half of the time horizon ($ p $-value being 5.0270e-05 using a two-sided Wilcoxon signed-rank test). Panel (c) shows that the leading arm is more likely to be chosen when the empirical mean gap is wider or is based on more observations. Overall, all the plots in Figure~\ref{fig:box-num-larger-2} indicate the findings in Section~\ref{subsec:behavioral-data-descriptive-analysis} are robust.

\noindent \textbf{Empirical validation on real behavioral data.}  Similar to Section~\ref{sec:structure-estimation-behavioral-data}, we compare QCARE with the same benchmark models, but using the behavioral data gathered from the online replication. The detailed comparison results are presented in Table~\ref{tb:comp-predict-power-2}. The online-replication results are broadly consistent with the analysis of the baseline experiment. QCARE has significantly higher out-of-sample average log-likelihood than L-\(n_j\), ES-\(\gamma_j\), QCARE-0, TS-B, and TS-G. In the head-to-head comparison, QCARE significantly outperforms L-\(n_j\), QCARE-0, \(\varepsilon_j\)-greedy, TS-B, and TS-G. Thus, QCARE remains competitive with the strongest heuristic benchmarks and performs strongly against the remaining models.

\begin{table}[htbp]
	\caption{ Comparison Between QCARE and Other Models (Online Replication)}
	\label{tb:comp-predict-power-2}
	\centering{
		\begin{tabular}{cccccccccc}
			\toprule
			\footnotesize{Data} & \footnotesize{Metric}  & QCARE   & L-$n_j$ & HH-$n_j$ & ES-$\gamma_j$ & QCARE-0 & $\varepsilon_j$-greedy & TS-B  & TS-G \\
			\midrule
			 \multirow{2}{*}{IS}  & AvgLL  & -69.0      & -66.6     & -63.3    & -65.3     & -74.0\dstar    & -69.6  & -116.8\dstar & -78.1\dstar  \\
			 & Pct & --   & 29.3\% & 20.3\% & 19.5\% & 100.0\%\dstar & 51.2\% & 97.6\%\dstar &100.0\%\dstar \\
			 \cmidrule{1-10}
			\multirow{2}{*}{OoS} & AvgLL   & -21.5    & -24.9\dstar    & -21.7     & -24.0\dstar    & -24.3\dstar      & -21.9  & -30.7\dstar & -26.6\dstar    \\
			& Pct &-- & 57.7\%\star & 53.7\% & 56.1\% & 87.8\%\dstar & 62.6\%\dstar & 87.0\%\dstar & 76.4\%\dstar \\
			\bottomrule
	\end{tabular}
	
	\vskip 0.2 cm
	\noindent\parbox{\linewidth}{\sf \scriptsize
		\textit{Note:} ``IS'' means in-sample (i.e., evaluated on the training data), and ``OoS'' means out-of-sample (i.e., evaluated on the test data). ``AvgLL'' refers to the average log-likelihood at the subject level by each approach. ``Pct'' indicates the corresponding percentage of subjects for whom QCARE outperforms the benchmark. The remarks ``**'' and ``*'' indicate $ p $-values less than $ 0.01 $ and $ 0.05 $, respectively.}
}
\end{table}

\begin{figure}
	\centering
	\begin{subfigure}{0.41\textwidth}
		\includegraphics[width=\textwidth]{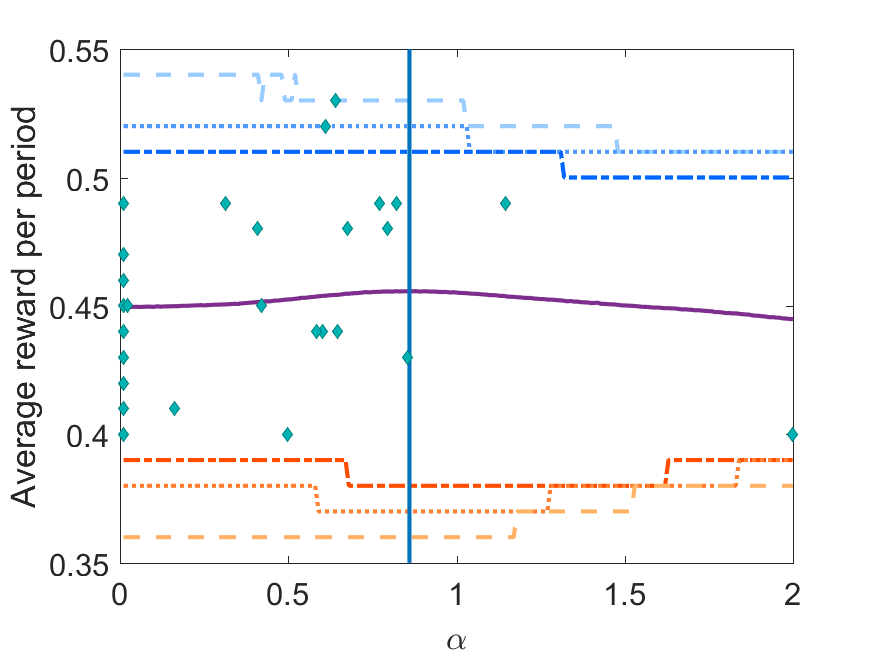}
		\subcaption{$\mu = (0.4, 0.5)$, $T=100$}
	\end{subfigure}
	\hspace{0.4 cm}
	\begin{subfigure}{0.41\textwidth}
		\includegraphics[width=\textwidth]{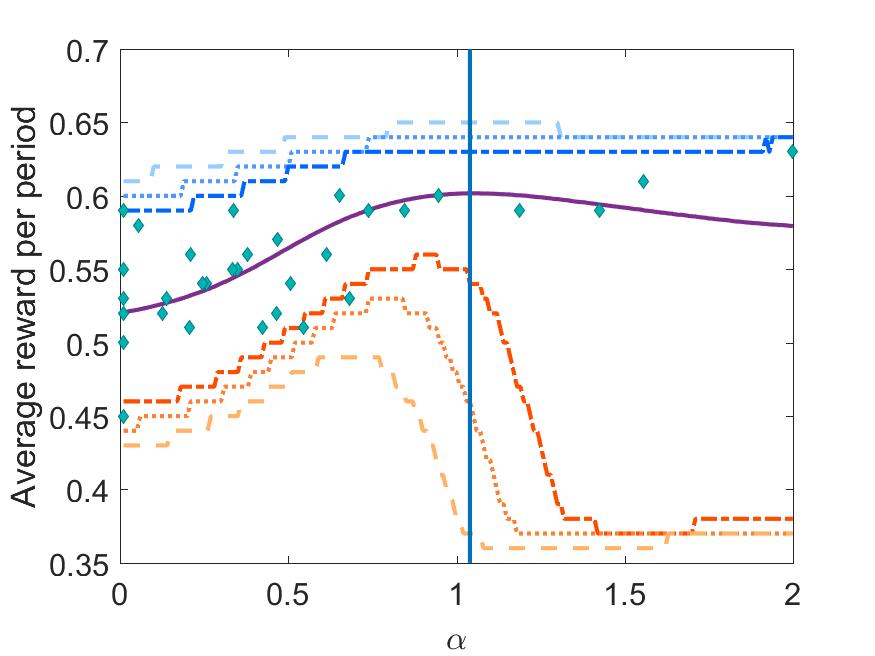}
		\subcaption{$\mu = (0.6,0.4)$, $T=100$}
	\end{subfigure}
	\begin{subfigure}{0.09\textwidth}
		\includegraphics[width=\textwidth]{figs/legend-beha.png}
	\end{subfigure}
	
	\begin{subfigure}{0.41\textwidth}
		\includegraphics[width=\textwidth]{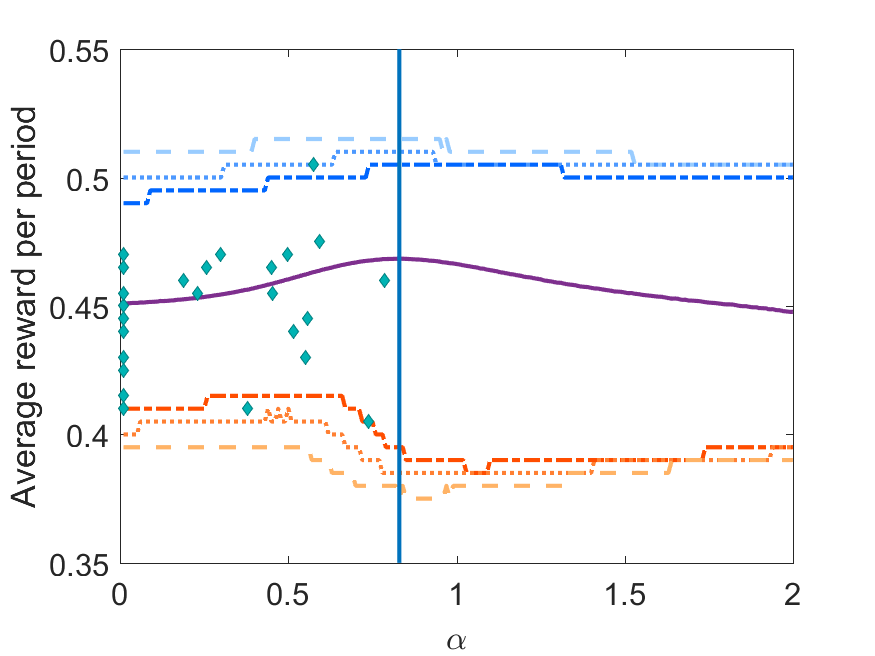}
		\subcaption{$\mu = (0.4, 0.5)$, $T=200$}
	\end{subfigure}
	\hspace{0.4 cm}
	\begin{subfigure}{0.41\textwidth}
		\includegraphics[width=\textwidth]{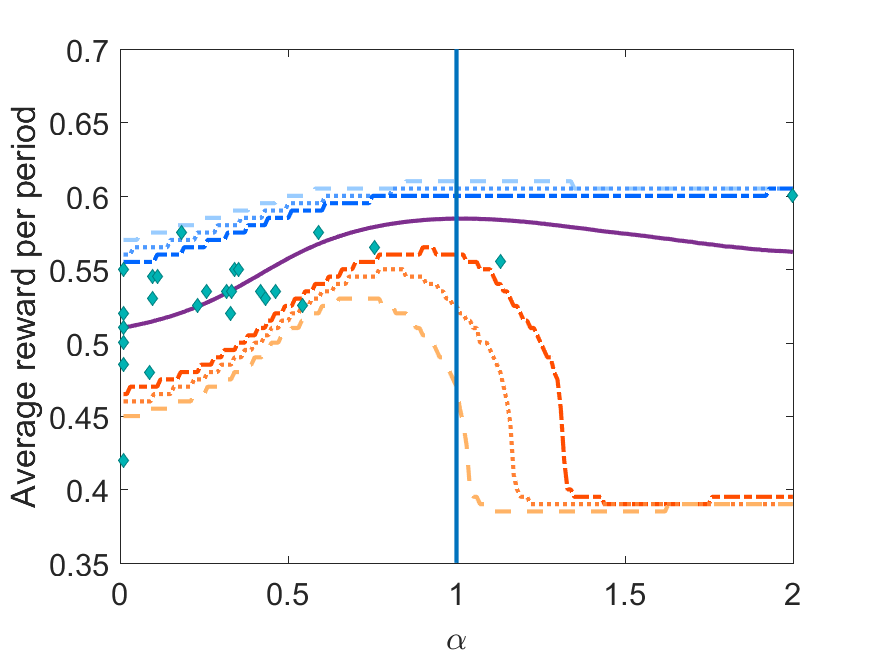}
		\subcaption{$\mu = (0.6,0.4)$, $T=200$}
	\end{subfigure}
	\begin{subfigure}{0.09\textwidth}
		\includegraphics[width=\textwidth]{figs/legend-beha.png}
	\end{subfigure}
	\caption{Comparison of Estimated $\alpha$ by Individual Subjects (Online Replication)}
	\label{fig:beha-individual-2}
\end{figure}

\noindent \textbf{Realized rewards and behavioral patterns.} 
We also investigate the relationship between the realized reward and the reduction rate of exploration, effectively replicating Figure~\ref{fig:beha-individual} using Figure~\ref{fig:beha-individual-2} based on online-replication data. As in the baseline experiment, most of the points lie within the percentile band implied by the QCARE model. In addition, the majority of points are positioned to the left of the optimal level. This observation suggests that the tendency towards over-exploration is robust in the online replication.

\subsection{Preregistered Replication}
\label{section:appendix:preregister experiment}

We conducted the preregistered replication to once again replicate our empirical analysis and confirm that our findings are robust.
The experimental procedure followed the same methodology outlined in Section~\ref{subsec:experiment-design}. The interfaces of the introduction and round pages, as well as the payment rules, were identical to those of the online replication.

\begin{figure}[htbp]
    \centering
    \begin{subfigure}[b]{0.33\textwidth}
        \centering
\includegraphics[width=\textwidth]{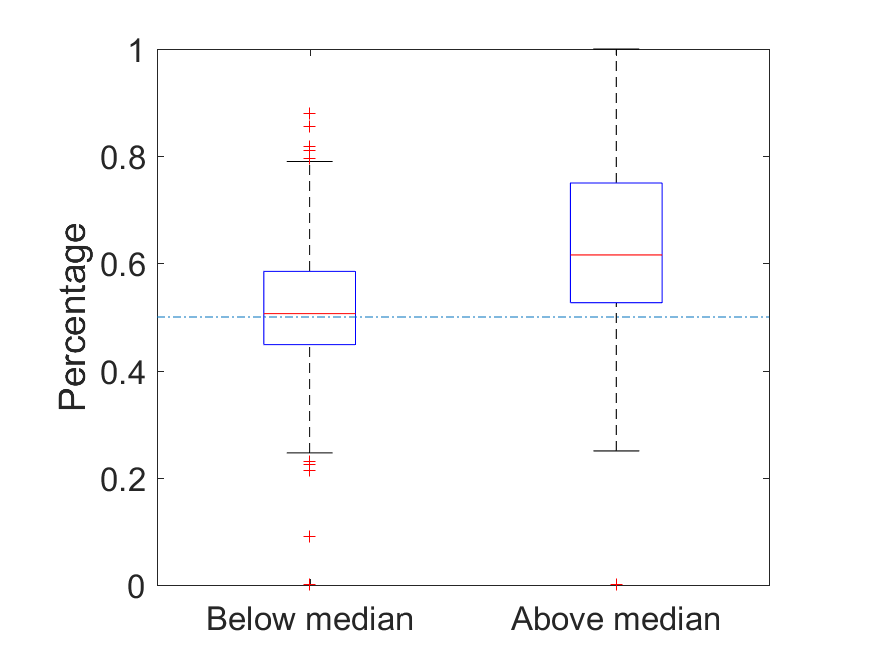}
        \caption{By empirical mean difference}
    \end{subfigure}
    \hspace{-2em}
    \begin{subfigure}[b]{0.33\textwidth}
        \centering
\includegraphics[width=\textwidth]{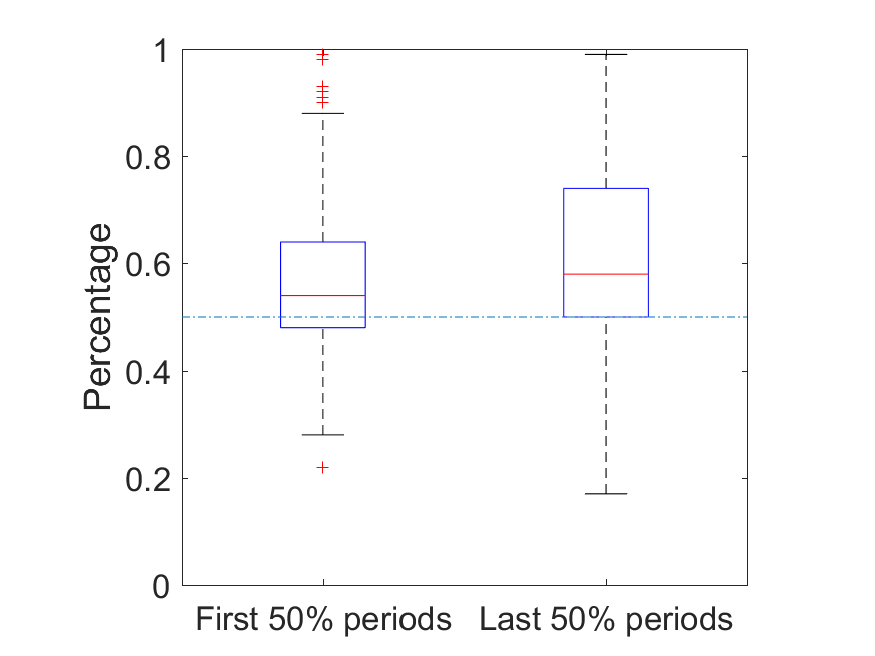}
        \caption{By period}
    \end{subfigure}
    \hspace{-2em}
    \begin{subfigure}[b]{0.33\textwidth}
        \centering
\includegraphics[width=\textwidth]{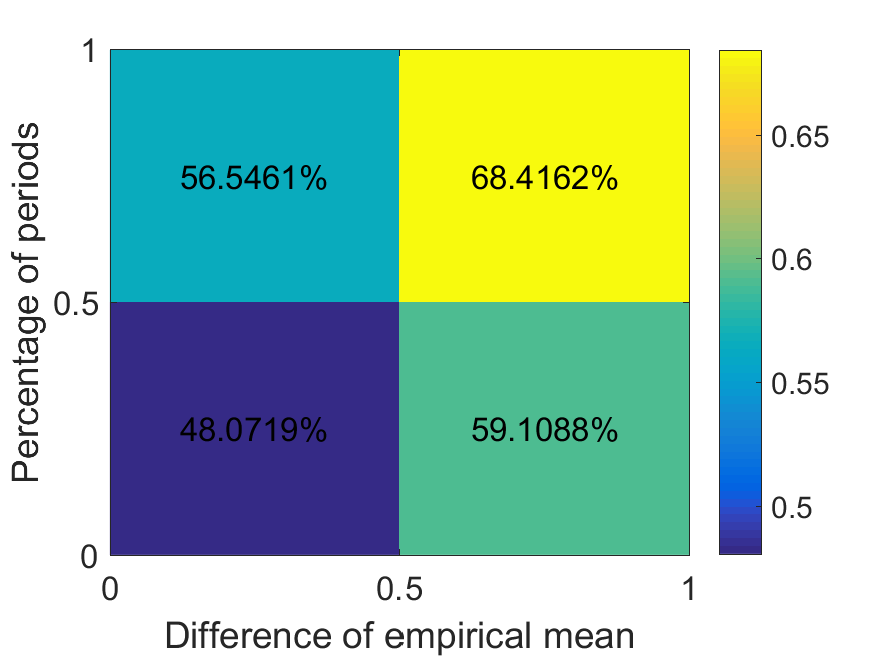}
        \caption{By both}
    \end{subfigure}
	\caption{Fraction of Periods Choosing the Leading Arm (Preregistered Replication)}
	\label{fig:box-num-larger-4}
\end{figure}

\noindent \textbf{Descriptive analysis.} 
The results from the preregistered replication show that only 4 of 303 subjects consistently chose the leading arm across all rounds, while 229 out of 303 subjects selected the leading arm in at least 50\% of the rounds. On average, subjects chose the leading arm $59.95\%$ of the time.\footnote{If taking the average at the individual-period level, the figure would be $ 60.33\% $.} That ratio is strictly less than one, yet significantly larger than $ 0.5 $ (Wilcoxon signed-rank test right tail, $ p $-value $ = $4.1155e-27).
Next, we use Figure~\ref{fig:box-num-larger-4} based on the preregistered replication as a replication for Figure~\ref{fig:box-num-larger}. In Panel (a), we observe a significantly higher average ratio when the difference exceeds the median ($ p $-value being 1.4501e-16 using a two-sided Wilcoxon signed-rank test). Panel (b) shows a significantly higher average ratio in the second half of the time horizon ($ p $-value being 2.3613e-10 using a two-sided Wilcoxon signed-rank test). Panel (c) indicates that the leading arm is more likely to be chosen when the empirical mean gap is wider or is based on more observations.
Overall, the plots in Figure~\ref{fig:box-num-larger-4} support the robustness of the findings in Section~\ref{subsec:behavioral-data-descriptive-analysis}.

\noindent \textbf{Comparison of prediction power.} 
As in Section~\ref{sec:structure-estimation-behavioral-data}, we compare the prediction power of the QCARE model against other benchmark models, but here using the preregistered-replication data. The detailed comparison results are provided in Table~\ref{tb:comp-predict-power-4}. The preregistered-replication results are broadly consistent with the earlier samples, although QCARE does not uniformly outperform every benchmark. QCARE achieves a significantly higher out-of-sample average log-likelihood than five out of seven benchmark methods. At the subject level, QCARE demonstrates statistically significant superiority over six of the seven competing methods. Thus, QCARE remains competitive with the strongest heuristic benchmark and performs favorably relative to the remaining models.
\begin{table}[htbp]
	\caption{ Comparison Between QCARE and Other Models (Preregistered Replication)}
	\label{tb:comp-predict-power-4}
	\centering{
		\begin{tabular}{cccccccccc}
			\toprule
			\footnotesize{Data} & \footnotesize{Metric}  & QCARE   & L-$n_j$ & HH-$n_j$ & ES-$\gamma_j$ & QCARE-0 & $\varepsilon_j$-greedy & TS-B  & TS-G \\
			\midrule
			 \multirow{2}{*}{IS}  & AvgLL  & -70.5      & -67.9     & -65.8    & -66.7     & -74.6\dstar    & -70.6  & -119.2\dstar & -79.2\dstar  \\
			 & Pct & --   & 24.4\% & 26.4\% & 17.8\% & 100.0\%\dstar & 52.5\% & 98.3\%\dstar &99.0\%\dstar \\
			 \cmidrule{1-10}
			\multirow{2}{*}{OoS} & AvgLL   & -23.4    & -25.1\dstar    & -22.0     & -24.6\star    & -24.6\dstar      & -23.4  & -29.3\dstar & -27.3\dstar    \\
			& Pct &-- & 57.8\%\dstar & 42.6\% & 55.1\%\star & 77.6\%\dstar & 55.1\%\star & 77.6\%\dstar & 73.9\%\dstar \\
			\bottomrule
	\end{tabular}
	
	\vskip 0.2 cm
	\noindent\parbox{\linewidth}{\sf \scriptsize
		\textit{Note:} ``IS'' means in-sample (i.e., evaluated on the training data), and ``OoS'' means out-of-sample (i.e., evaluated on the test data). ``AvgLL'' refers to the average log-likelihood at the subject level by each approach. ``Pct'' indicates the corresponding percentage of subjects for whom QCARE outperforms the benchmark. The remarks ``**'' and ``*'' indicate $ p $-values less than $ 0.01 $ and $ 0.05 $, respectively.}
}
\end{table}

\begin{figure}
	\centering
	\begin{subfigure}{0.41\textwidth}
		\includegraphics[width=\textwidth]{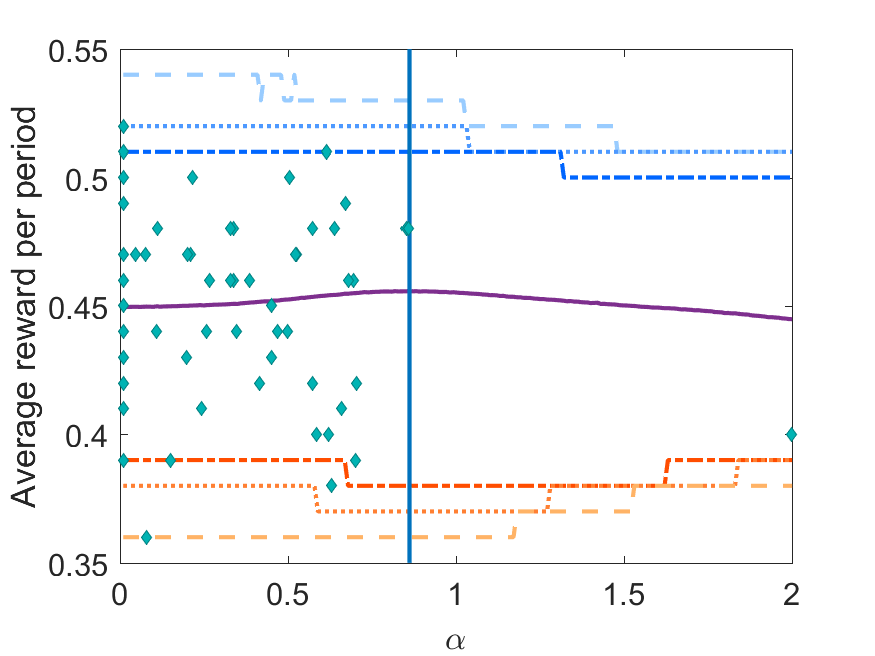}
		\subcaption{$\mu = (0.4, 0.5)$, $T=100$}
	\end{subfigure}
	\hspace{0.4 cm}
	\begin{subfigure}{0.41\textwidth}
		\includegraphics[width=\textwidth]{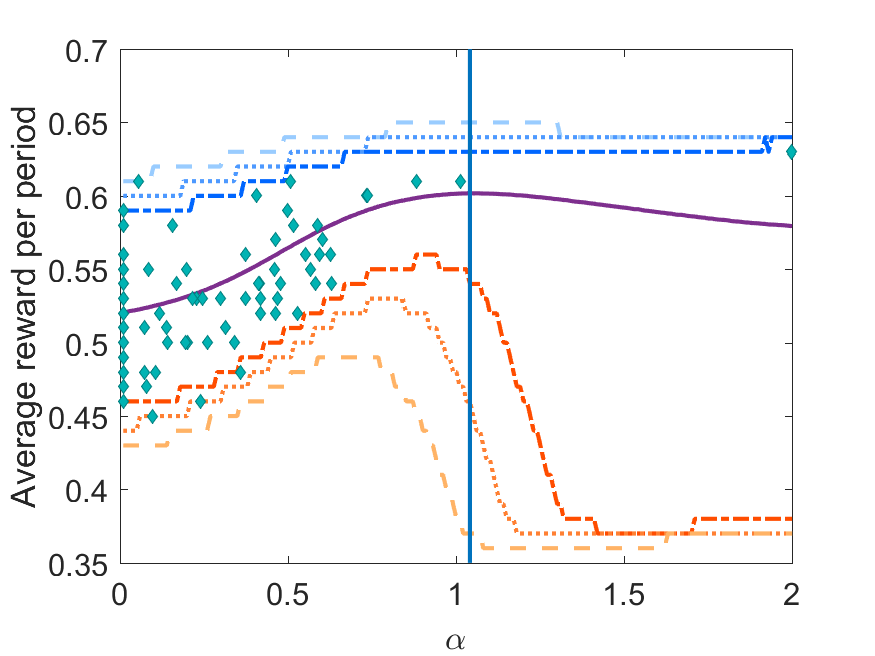}
		\subcaption{$\mu = (0.6,0.4)$, $T=100$}
	\end{subfigure}
	\begin{subfigure}{0.09\textwidth}
		\includegraphics[width=\textwidth]{figs/legend-beha.png}
	\end{subfigure}
	
	\begin{subfigure}{0.41\textwidth}
		\includegraphics[width=\textwidth]{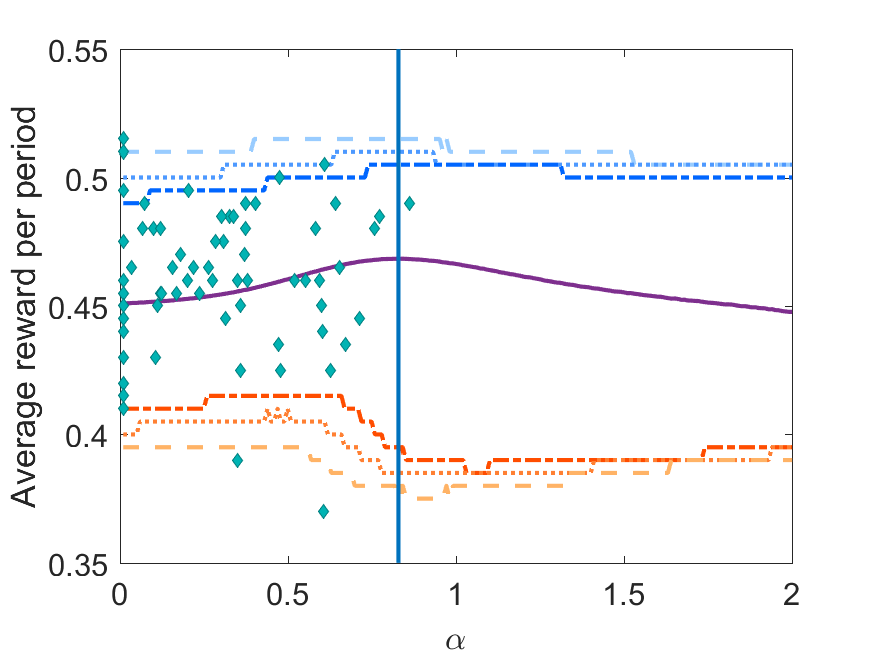}
		\subcaption{$\mu = (0.4, 0.5)$, $T=200$}
	\end{subfigure}
	\hspace{0.4 cm}
	\begin{subfigure}{0.41\textwidth}
		\includegraphics[width=\textwidth]{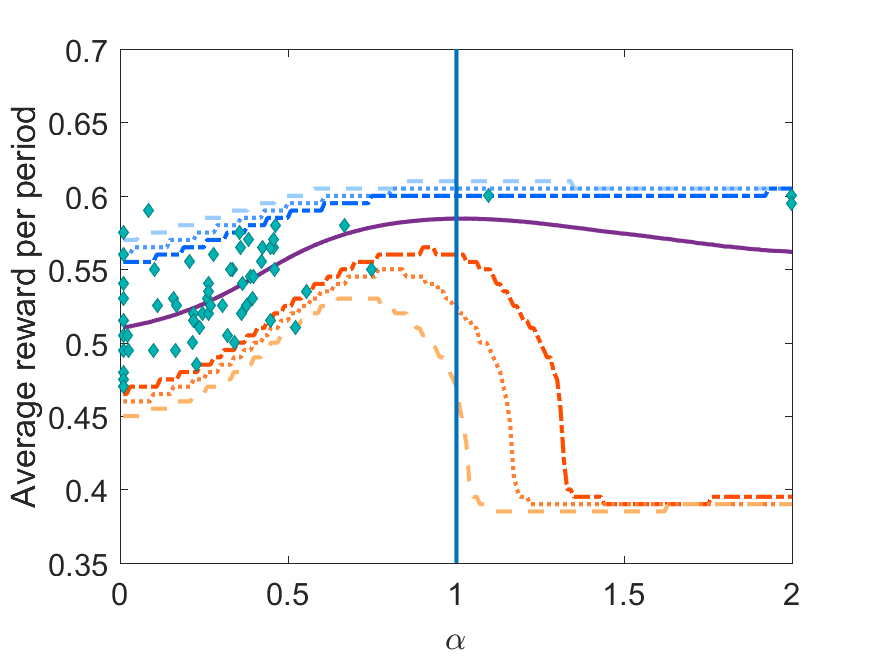}
		\subcaption{$\mu = (0.6,0.4)$, $T=200$}
	\end{subfigure}
	\begin{subfigure}{0.09\textwidth}
		\includegraphics[width=\textwidth]{figs/legend-beha.png}
	\end{subfigure}
	\caption{Comparison of Estimated $\alpha$ by Individual Subjects (Preregistered Replication)}
	\label{fig:beha-individual-4}
\end{figure}

\noindent \textbf{Realized rewards and behavioral patterns.} 
We also examine the relationship between the realized reward and the reduction rate of exploration, replicating Figure~\ref{fig:beha-individual} using data collected from the preregistered replication.  We report the results in Figure~\ref{fig:beha-individual-4}. Once again, it is clear that most green points are positioned to the left of the optimal level, suggesting that, as in previous results, the majority of subjects exhibit over-exploration behavior.

\subsection{Comparison and Pooled Validation Across Subject Pools}\label{section:appendix:subject pools}

To complement the baseline descriptive results in Section~\ref{subsec:behavioral-data-descriptive-analysis} and the sample-specific preliminary results reported in Appendices~\ref{section:appendix:online experiment} and~\ref{section:appendix:preregister experiment}, Table~\ref{tab:descriptive_stats_subject_pools} summarizes the key descriptive statistics across the baseline, online-replication, and preregistered-replication samples.

\begin{table}[htbp]
	\centering
		\caption{Descriptive Statistics Across Subject Pools}
		\label{tab:descriptive_stats_subject_pools}
		\begin{threeparttable}
			\begin{tabular}{lcccc}
				\toprule
				& & \multicolumn{3}{c}{Choice of the leading arm} \\
				\cmidrule(lr){3-5}
				Subject pool & Subjects
				& \begin{tabular}[b]{@{}c@{}}Always\\ chosen\end{tabular}
				& \begin{tabular}[b]{@{}c@{}}Chosen in at least $50\%$\\ of rounds\end{tabular}
				& \begin{tabular}[b]{@{}c@{}}Average choice rate\\ $[p$-value$]$\end{tabular} \\
				\midrule
				Laboratory (Baseline) & 123 & 2 (1.6\%) & 100 (81.3\%)           & 64.70\%\; $[8.68\times10^{-15}]$ \\
				Online replication & 123 & 1 (0.8\%) & \phantom{0}87 (70.7\%) & 59.76\%\; $[8.44\times10^{-10}]$ \\
				Preregistered replication & 303 & 4 (1.3\%) & 229 (75.6\%)           & 59.95\%\; $[4.12\times10^{-27}]$ \\
				\bottomrule
			\end{tabular}
			\begin{tablenotes}[flushleft]\footnotesize
				\item[] \textit{Note}: Figures in parentheses are shares of subjects in the corresponding pool. Bracketed figures are $p$-values from one-sided Wilcoxon signed-rank tests of whether subject-level leading-arm choice rates are centered above \(50\%\).
			\end{tablenotes}
		\end{threeparttable}
\end{table}

The summary shows a consistent qualitative pattern across the three samples. In all samples, only a very small number of subjects always choose the leading arm, indicating that subjects are not purely exploitative. At the same time, the average leading-arm choice rate is significantly above 50\% in all three samples, indicating that subjects are not choosing randomly or purely exploring either. Although the online-replication and preregistered-replication samples have slightly lower average leading-arm choice rates than the baseline sample, the main qualitative pattern remains the same. Subjects balance exploration and exploitation while still making substantial non-leading-arm choices.

\begin{figure}
	\centering
	\begin{subfigure}{0.41\textwidth}
		\includegraphics[width=\textwidth]{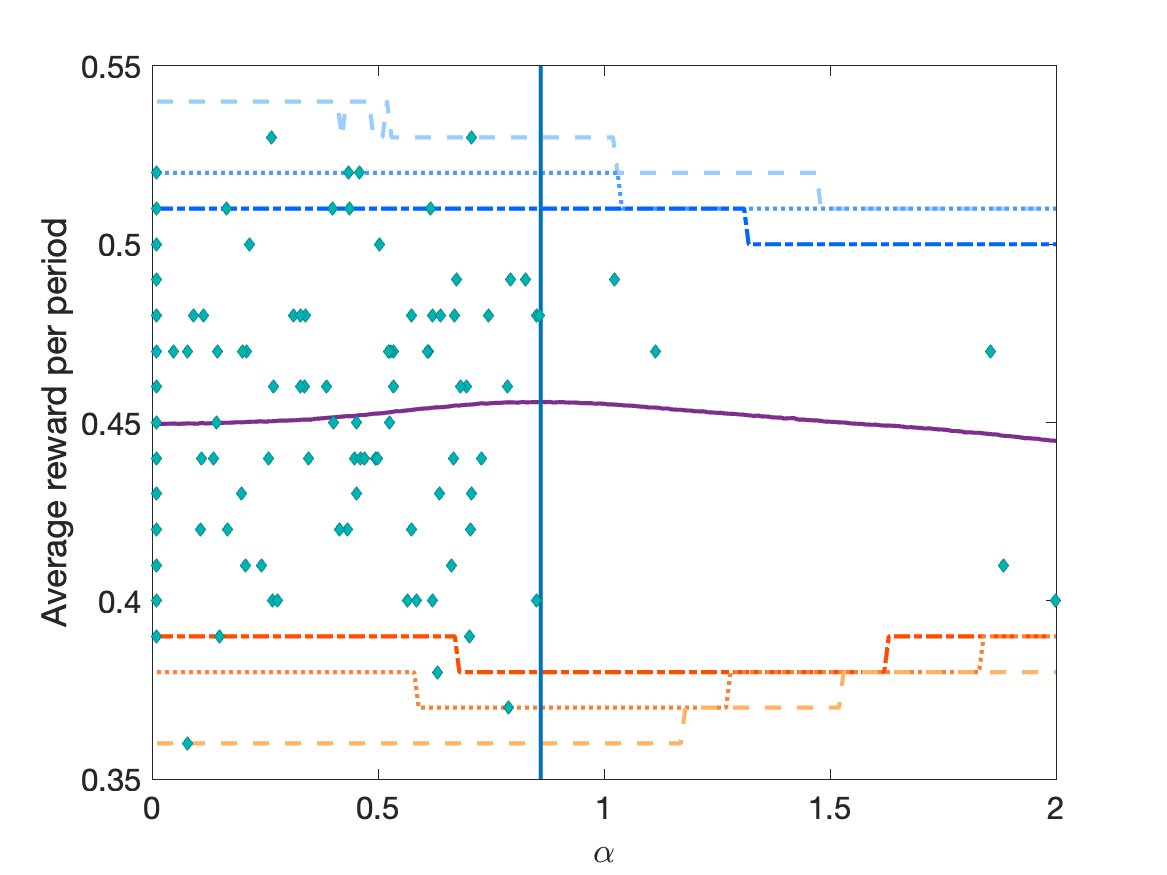}
		\subcaption{$\mu = (0.4, 0.5)$, $T=100$}
	\end{subfigure}
	\hspace{0.4 cm}
	\begin{subfigure}{0.41\textwidth}
		\includegraphics[width=\textwidth]{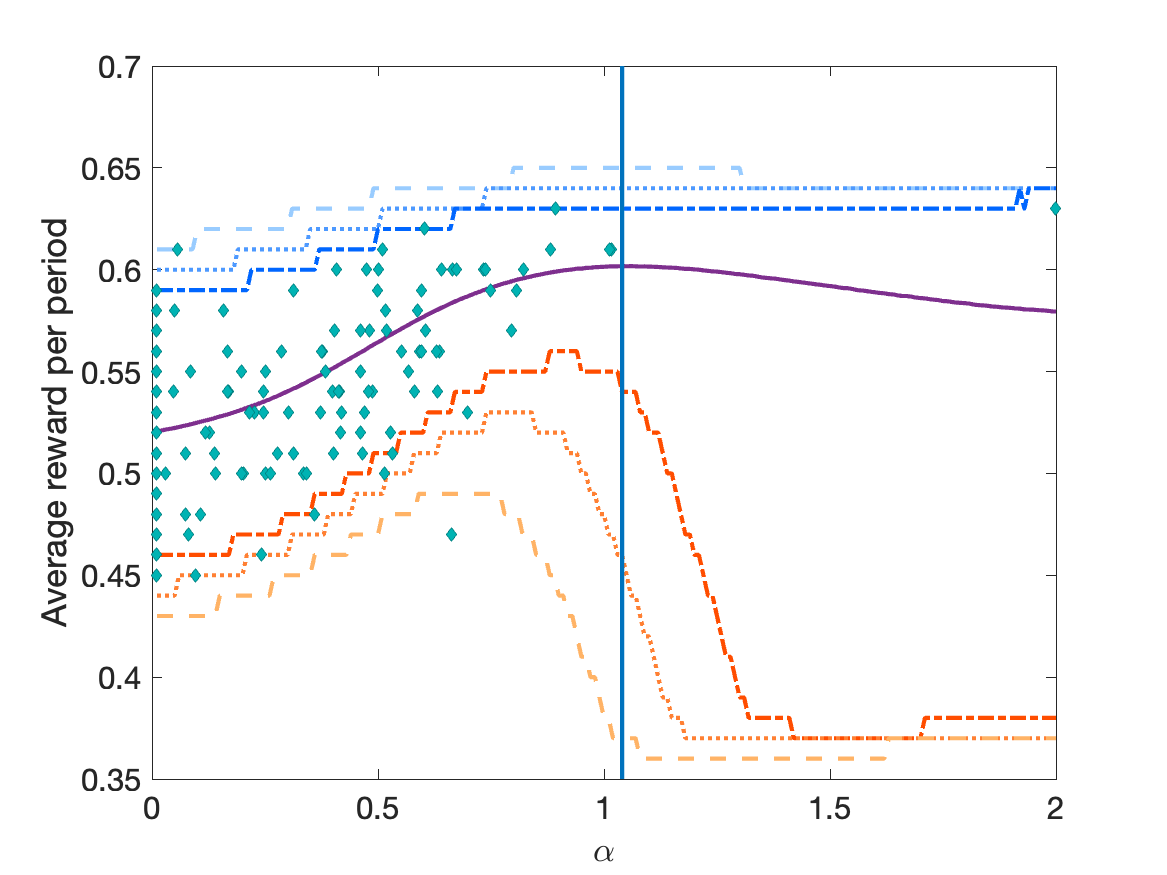}
		\subcaption{$\mu = (0.6,0.4)$, $T=100$}
	\end{subfigure}
	\begin{subfigure}{0.09\textwidth}
		\includegraphics[width=\textwidth]{figs/legend-beha.png}
	\end{subfigure}
	
	\begin{subfigure}{0.41\textwidth}
		\includegraphics[width=\textwidth]{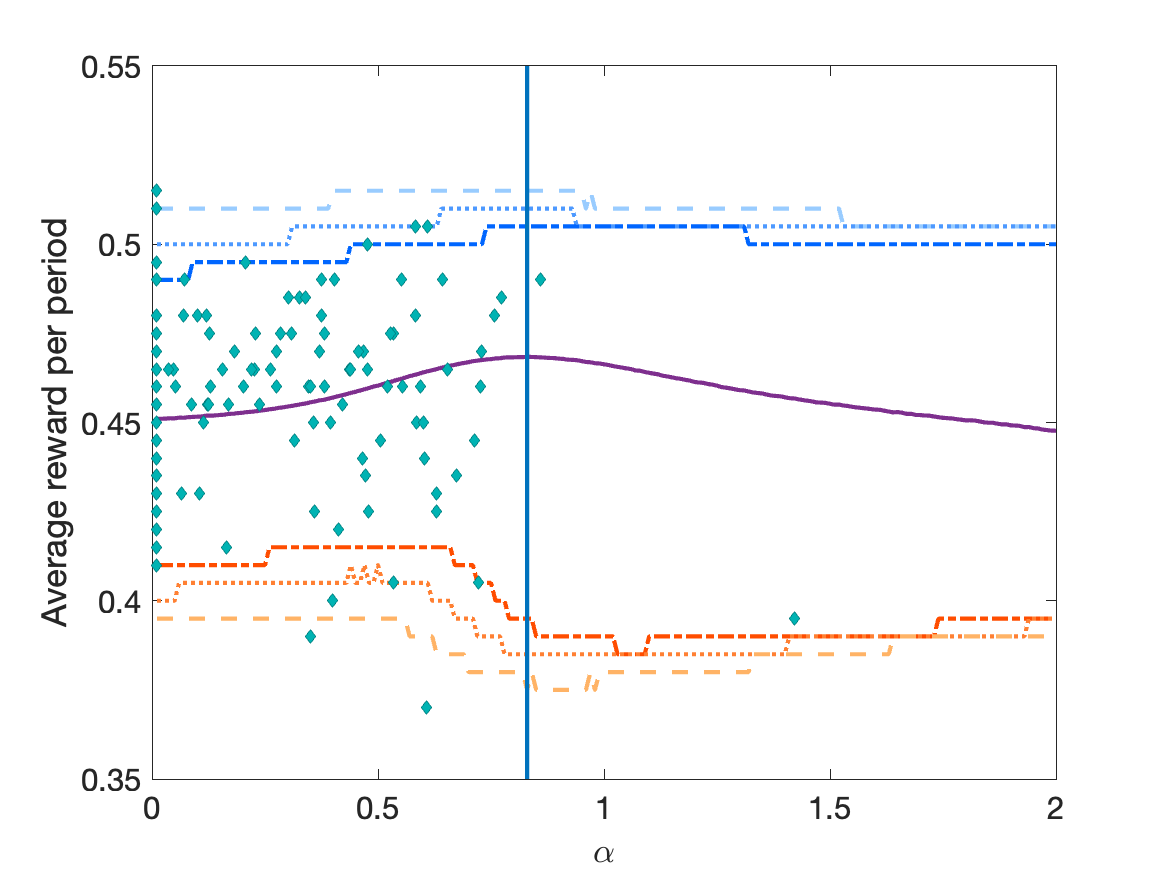}
		\subcaption{$\mu = (0.4, 0.5)$, $T=200$}
	\end{subfigure}
	\hspace{0.4 cm}
	\begin{subfigure}{0.41\textwidth}
		\includegraphics[width=\textwidth]{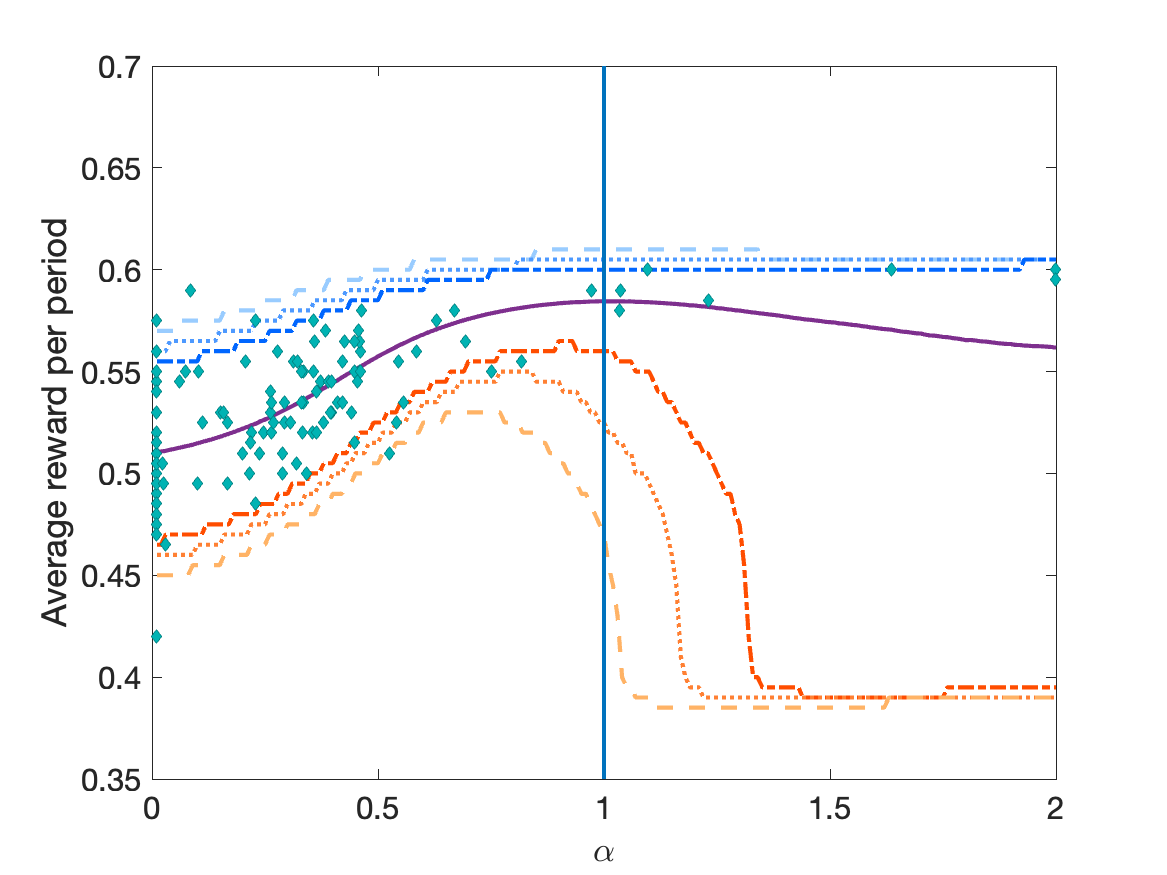}
		\subcaption{$\mu = (0.6,0.4)$, $T=200$}
	\end{subfigure}
	\begin{subfigure}{0.09\textwidth}
		\includegraphics[width=\textwidth]{figs/legend-beha.png}
	\end{subfigure}
	\caption{Comparison of Estimated \(\alpha\) by Individual Subjects (Pooled)}
\noindent\parbox{\linewidth}{\sf \scriptsize
\textit{Note:} The green dots correspond to subjects from the baseline and the two replication experiments, with 549 subjects in total.}
\label{fig:beha-individual-pooled}
\end{figure}

To complement the baseline validation figure in Section~\ref{sec:structure-estimation-behavioral-data} and the sample-specific validation figures reported above, Figure~\ref{fig:beha-individual-pooled} presents the corresponding validation results using the pooled two-arm samples from the baseline experiment, along with the online and preregistered replications. The pooled sample includes 549 subjects in total.

\textbf{Time-position effects.} This subsection provides the details of the diagnostic test referenced in the remark on time-position effects in Section~\ref{subsec: model_interpretation}. We examine whether subjects reduce exploration more aggressively as they approach the end of the horizon. For each subject, we split the decision sequence into the first and second halves of the horizon and estimate \(\alpha\) separately for the two segments. We then compare the first-half and second-half estimates using paired \(t\)-tests. The results are summarized in Table~\ref{tab:response_alpha_half_comparison}.

\begin{table}[htbp]
\centering
\caption{Diagnostic Comparison of First-Half and Second-Half Estimated \(\alpha\)}
\label{tab:response_alpha_half_comparison}
\begin{tabular}{lccccc}
\toprule
Sample & Subjects & First half & Second half & Change & \(p\)-value \\
\midrule
Laboratory (Baseline) & 123 & 0.437 & 0.471 & 0.034 & 0.379 \\
Online replication & 123 & 0.296 & 0.369 & 0.073 & 0.008 \\
Preregistered replication & 303 & 0.308 & 0.315 & 0.007 & 0.719 \\
\bottomrule
\end{tabular}
\end{table}

As can be seen, the magnitude of the estimated change in $ \alpha $ varies considerably across samples, and it is statistically significant only in the online replication experiment. The change in the preregistered replication, which is the one with the largest sample size, is particularly weak.
Hence, the evidence that $ \alpha $ differs between the two halves of the horizon does not appear robust. Therefore, the current empirical evidence seems to favor a simpler model with a fixed $\alpha$. We retain the parsimonious QCARE specification as the main model.

\vskip 0.2 cm
\subsection{Auxiliary Experiments}
\label{appendix:additional-experimental-checks}

\noindent \textbf{Three-arm extension.} In this auxiliary experiment, the rewards followed Bernoulli distributions with means 0.6, 0.4, and 0.2, respectively.
The show-up fee was 5 CNY, and each point was worth 0.4 CNY. In the following analysis, we replicate the empirical analysis of the baseline experiment to show that the findings from the two-arm setting remain robust in the three-arm extension.

Let us first conduct the descriptive analysis as in the main body. Our experimental data reveal that only 2 of 61 subjects chose the leading arm throughout the entire time horizon. On average, subjects chose the leading arm 62.52\% of the time.\footnote{If taking the average at the individual-period level, the figure would be $ 63.77\% $.} That ratio is strictly less than one, yet significantly larger than $ 0.5 $ (Wilcoxon signed-rank test right tail, $ p $-value $ = $4.9013e-6). Hence, at the aggregate level, the participants' choices reflect neither pure exploration nor pure exploitation, as expected. We also calculate the percentage of rounds in which each individual subject selected the second-best and third-best performing arms. The result shows that on average, subjects chose the second-performing (resp. third-performing) arm in 22.02\% (resp. 14.71\%) of rounds. This result indicates that subjects chose the second performing arm more often than the third performing arm. This difference is statistically significant by a two-sided Wilcoxon signed-rank test with a $ p $-value of 1.1683e-6. These findings suggest that human decision-making in this context does not follow the $\varepsilon$-greedy model, which predicts that these empirically suboptimal arms have an equal probability of being selected. 

\begin{figure}[htbp]
	\centering
	\begin{subfigure}[b]{0.33\textwidth}
		\centering
		\includegraphics[width=\textwidth]{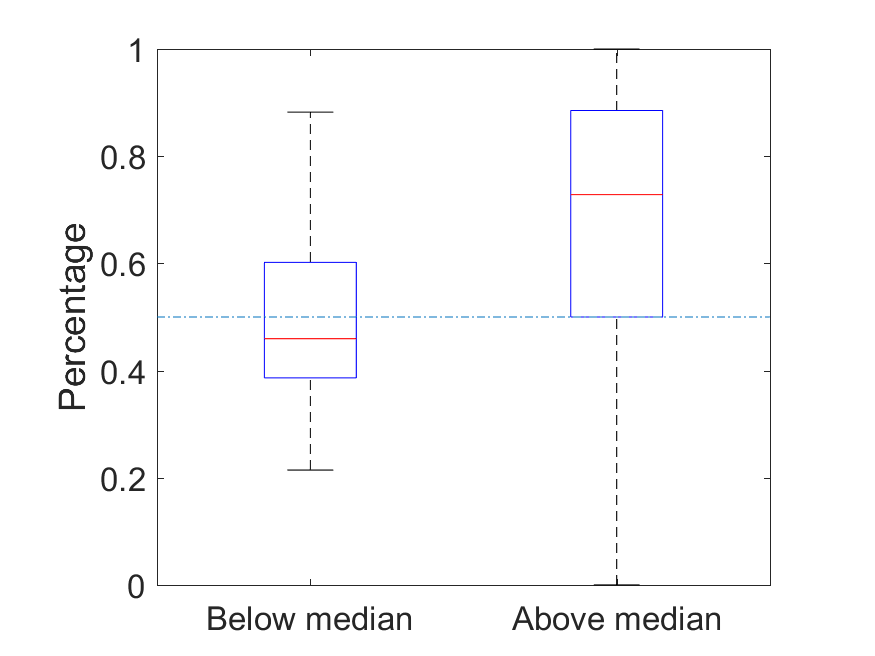}
		\caption{By empirical mean difference}
	\end{subfigure}
	\hspace{-2em}
	\begin{subfigure}[b]{0.33\textwidth}
		\centering
		\includegraphics[width=\textwidth]{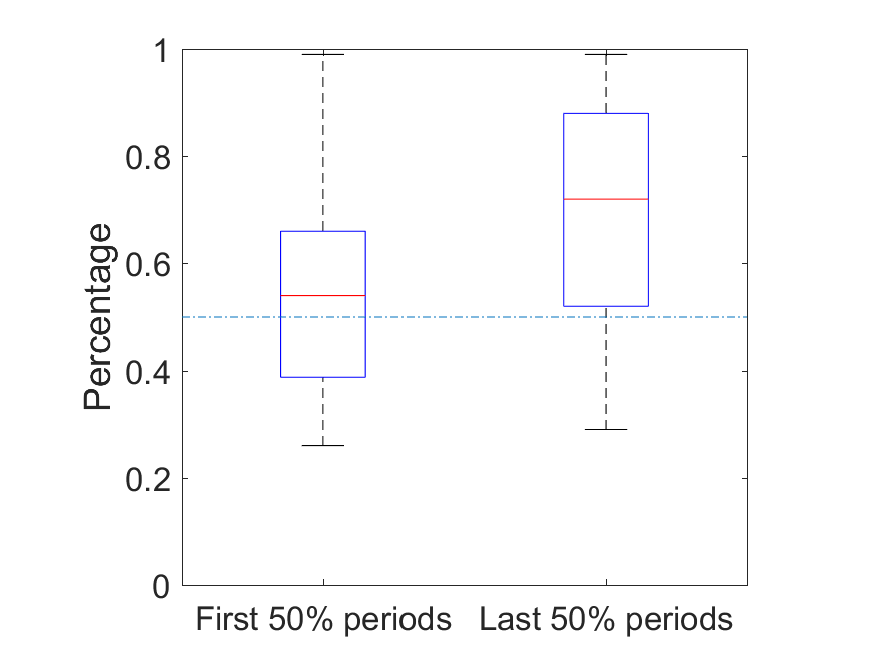}
		\caption{By period}
	\end{subfigure}
	\hspace{-2em}
	\begin{subfigure}[b]{0.33\textwidth}
		\centering
		\includegraphics[width=\textwidth]{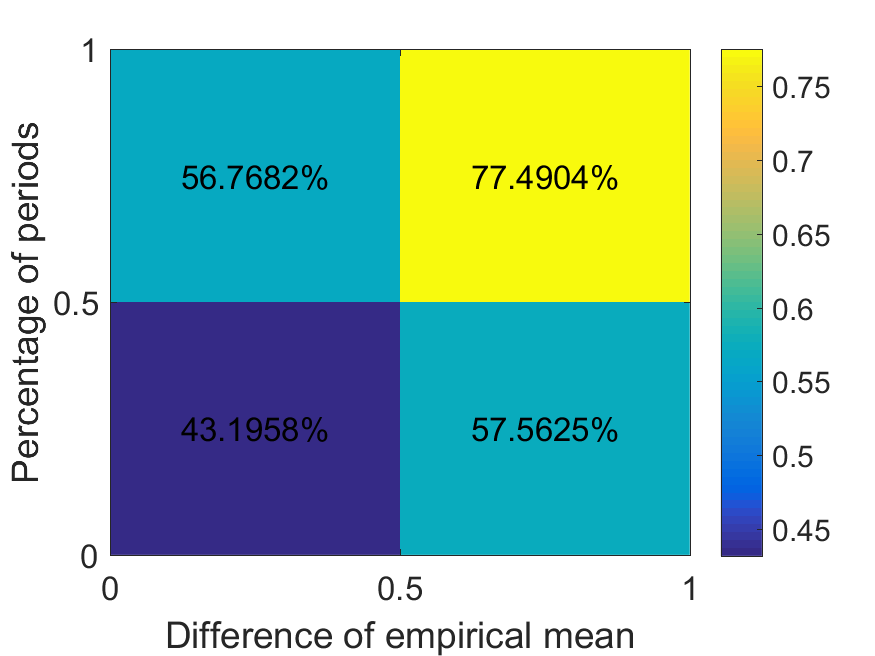}
		\caption{By both}
	\end{subfigure}
	\caption{Fraction of Periods Choosing the Leading Arm (Three-Arm Extension)}
	\label{fig:box-num-larger-3arm}
\end{figure}

Next, we use Figure~\ref{fig:box-num-larger-3arm} based on the three-arm-extension data as a replication for Figure~\ref{fig:box-num-larger} in the main body.
In Panel (a), we use the empirical mean difference between the leading arm and the second-best-performing arm. Overall, all the plots in Figure~\ref{fig:box-num-larger-3arm} indicate the findings in Section~\ref{subsec:behavioral-data-descriptive-analysis} are robust.

We next verify whether the estimated $\alpha$ can recover the true $\alpha$ value when subject $j$ indeed adheres to the QCARE model. Specifically, we simulate subject $j$'s behavior based on the QCARE model for a given $T$ and $\alpha$ in 10,000 samples. The mean vector of the three arms is $(0.6, 0.4, 0.2)$, which is aligned with that used in the laboratory experiment.
Table~\ref{tb:comp-alpha-estimation-3-arm} presents the summary of the estimated values of $\hat{\alpha}$. Again, the ground truth $\alpha$ is always within the 5th and 95th percentiles of $\hat{\alpha}$, and the average $\hat{\alpha}$ is very close to the ground truth $\alpha$, especially when the ground truth $ \alpha $ becomes larger.

\begin{table}[htbp]
	\caption{Comparison of Ground-Truth and Estimated $\alpha$ (Three-Arm Extension)}
	\label{tb:comp-alpha-estimation-3-arm}
	\centering{
		\begin{tabular}{cc|ccc}
			\toprule
			$ T $     &  $\alpha$ & average $\hat{\alpha}$ & \makecell{5th percentile \\of $\hat{\alpha}$} & \makecell{95th percentile \\of $\hat{\alpha}$} \\
			\midrule
			100   & 0.2 & 0.1978 & 0.0139 & 0.3637 \\
			200   & 0.2 & 0.1978 & 0.0916 & 0.2860 \\
			1000  & 0.2 & 0.2003 & 0.1771 & 0.2238 \\
			10000 & 0.2 & 0.2005 & 0.2004 & 0.2004 \\
			100   & 0.5 & 0.5041 & 0.4026 & 0.6047 \\
			200   & 0.5 & 0.5020 & 0.4570 & 0.5503 \\
			1000  & 0.5 & 0.5009 & 0.4881 & 0.5114 \\
			10000 & 0.5 & 0.5007 & 0.4958 & 0.5036 \\
			100   & 1   & 1.0193 & 0.9001 & 1.1721 \\
			200   & 1   & 1.0119 & 0.9234 & 1.1177 \\
			1000  & 1   & 1.0051 & 0.9467 & 1.0633 \\
			10000 & 1   & 1.0019 & 0.9622 & 1.0400 \\
			\bottomrule
	\end{tabular}}
\end{table}

Similar to Section~\ref{sec:structure-estimation-behavioral-data}, we compare the QCARE model with the same benchmark models in the main body, but under the three-arm extension. The detailed comparison results are presented in Table~\ref{tb:comp-predict-power-3-arm}. Similar to the two-arm setting, QCARE continues to display favorable performance, especially in terms of out-of-sample prediction power.

\begin{table}[htbp]
	\caption{ Comparison Between QCARE and Other Models (Three-Arm Extension)}
	\label{tb:comp-predict-power-3-arm}
	\centering{
		\begin{tabular}{cccccccccc}
			\toprule
			\footnotesize{Data} & \footnotesize{Metric}  & QCARE   & L-$n_j$ & HH-$n_j$ & ES-$\gamma_j$ & QCARE-0 & $\varepsilon_j$-greedy & TS-B  & TS-G \\
			\midrule
			\multirow{2}{*}{IS}  & AvgLL  & -91.7      & -91.8     & -95.0    & -85.9     & -111.1\dstar    & -92.8  & -163.2\dstar & -104.4\dstar  \\
			& Pct & --   & 50.8\% & 52.5\% & 18.0\% & 100.0\%\dstar & 55.7\% & 95.1\%\dstar &98.4\%\dstar \\
			\cmidrule{1-10}
			\multirow{2}{*}{OoS} & AvgLL   & -27.0    & -61.6\dstar    & -29.1\star     & -37.0\dstar    & -35.4\dstar      & -27.4  & -44.7\dstar & -32.1\dstar    \\
			& Pct &-- & 78.7\%\dstar & 62.3\%\star & 75.4\%\dstar & 82.0\%\dstar & 63.9\%\star & 86.9\%\dstar & 67.2\%\dstar \\
			\bottomrule
		\end{tabular}
		
		\vskip 0.2 cm
		\noindent\parbox{\linewidth}{\sf \scriptsize
			\textit{Note:} ``IS'' means in-sample (i.e., evaluated on the training data), and ``OoS'' means out-of-sample (i.e., evaluated on the test data). ``AvgLL'' refers to the average log-likelihood at the subject level by each approach. ``Pct'' indicates the corresponding percentage of subjects for whom QCARE outperforms the benchmark. The remarks ``**'' and ``*'' indicate $ p $-values less than $ 0.01 $ and $ 0.05 $, respectively.}
	}
\end{table}

\begin{figure}
	\centering
	\begin{subfigure}{0.41\textwidth}
		\centering
		\includegraphics[width=\textwidth]{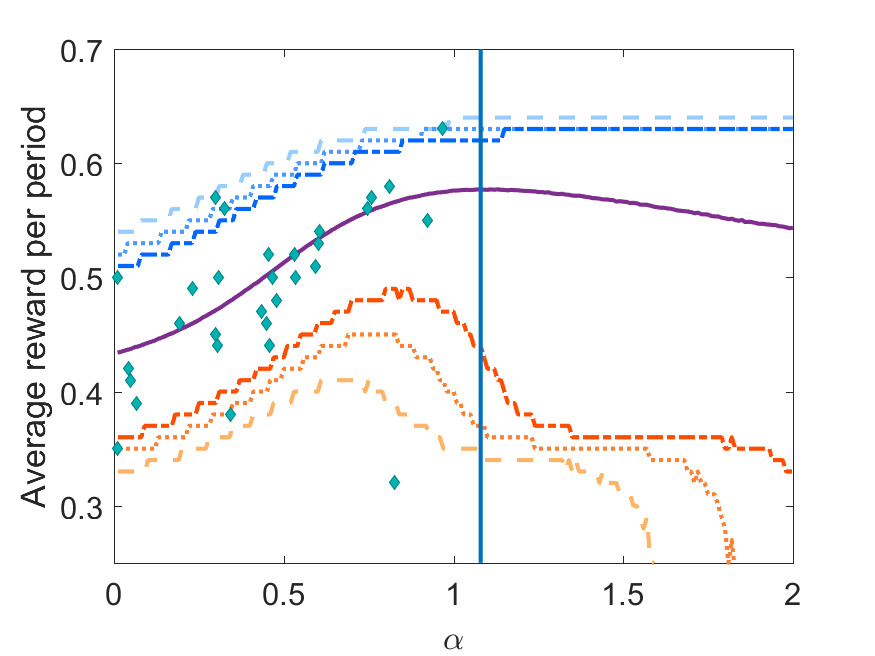}
		\caption{$\mu = (0.6, 0.4, 0.2)$, $T=100$}
	\end{subfigure}
	\hspace{0.4cm}
	\begin{subfigure}{0.41\textwidth}
		\centering
		\includegraphics[width=\textwidth]{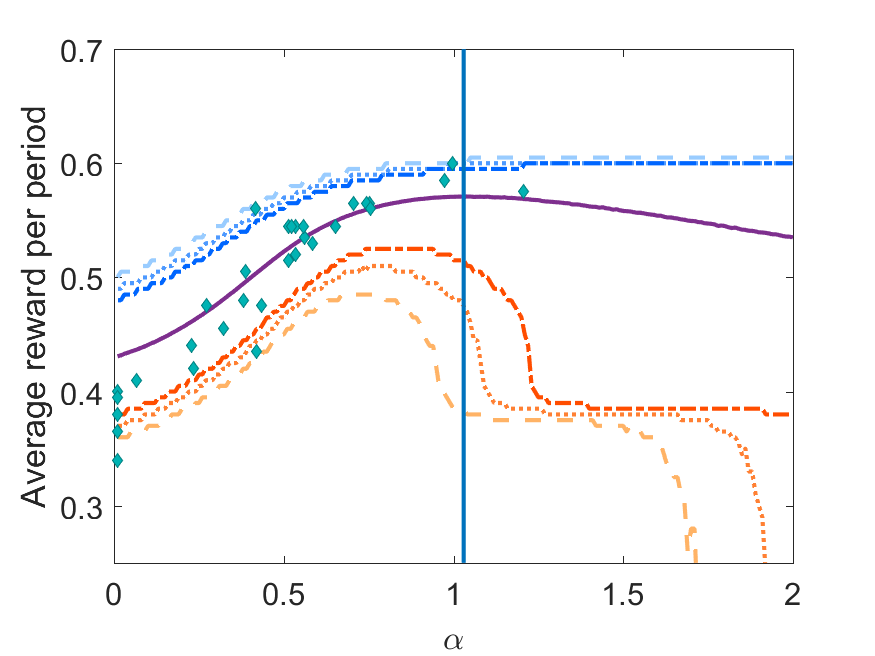}
		\caption{$\mu = (0.6, 0.4, 0.2)$, $T=200$}
	\end{subfigure}
	\begin{subfigure}{0.09\textwidth}
		\includegraphics[width=\textwidth]{figs/legend-beha.png}
	\end{subfigure}
	
	\caption{Comparison of Estimated $\alpha$ by Individual Subjects (Three-Arm Extension)}
	\label{fig:beha-individual-3-arm}
\end{figure}

Finally, we replicate the analysis in Section~\ref{sec:structure-estimation-behavioral-data} and examine the relationship between the realized reward and the reduction rate of exploration. We present a comprehensive illustration of reward distribution with different $\alpha$ values, as well as the subjects' behavior in Figure~\ref{fig:beha-individual-3-arm}. Similar to the two-arm setting, the green points fall within the percentile bands generated by the QCARE model, further validating the ability of the QCARE model to capture the actual human behavior.  In addition, the majority of the green points are positioned to the left of the blue vertical line, suggesting a tendency towards over-exploration.

\vskip 0.2 cm

\noindent\textbf{No-summary-table extension.} A possible concern is that displaying empirical means and pull counts might guide participants toward policies based on cumulative statistics. To examine this possibility, we conducted the no-summary-table extension with 29 participants using horizon \(T=100\) and reward-mean vector \((0.6,0.4)\). We removed the table containing the computed empirical means and pull counts while retaining the complete period-by-period history of each participant's actions and realized rewards. The interface also continued to display the current period and total horizon. 

\begin{figure}[htbp]
	\centering
	\begin{subfigure}{0.45\textwidth}
		\includegraphics[width=\textwidth]{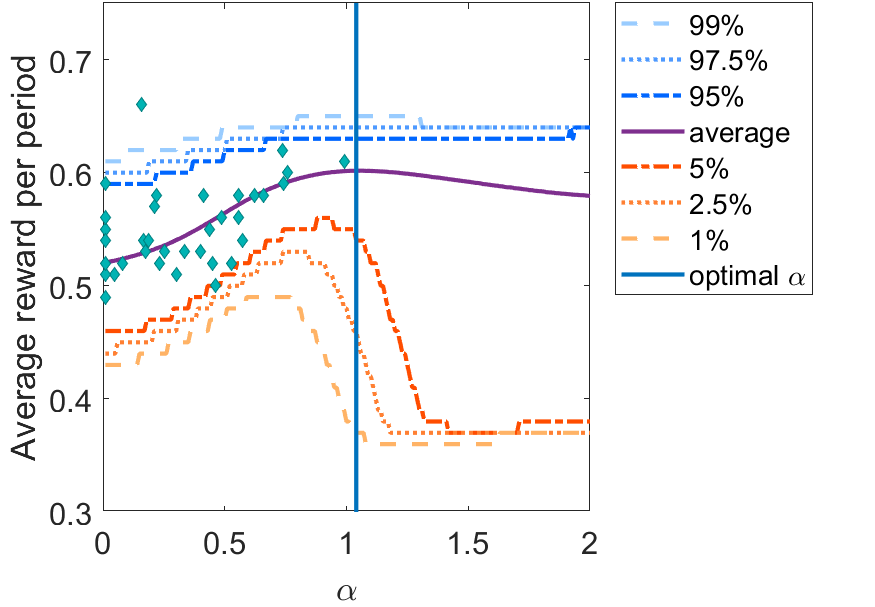}
		\subcaption{Practice rounds}
	\end{subfigure}
    \hspace{0.4 cm}
    \begin{subfigure}{0.45\textwidth}
		\includegraphics[width=\textwidth]{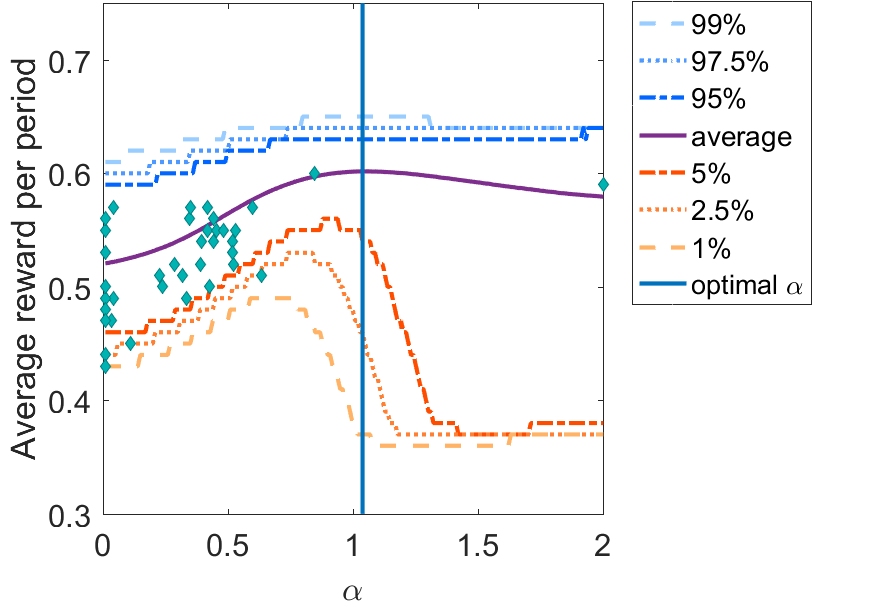}
		\subcaption{Higher stakes}
	\end{subfigure}
    
	\begin{subfigure}{0.45\textwidth}
		\includegraphics[width=\textwidth]{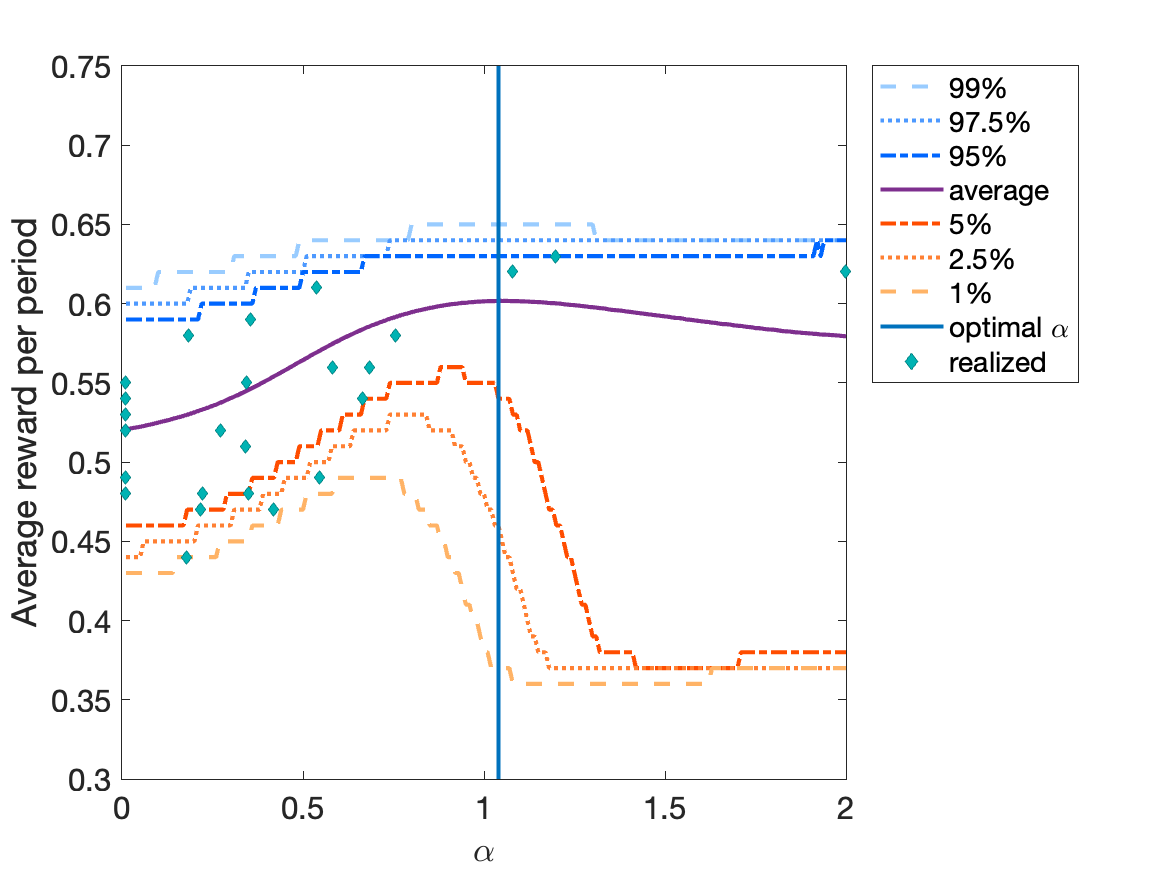}
		\subcaption{No-summary-table extension}
	\end{subfigure}
    \hspace{0.4 cm}
    \begin{subfigure}{0.45\textwidth}
		\includegraphics[width=\textwidth]{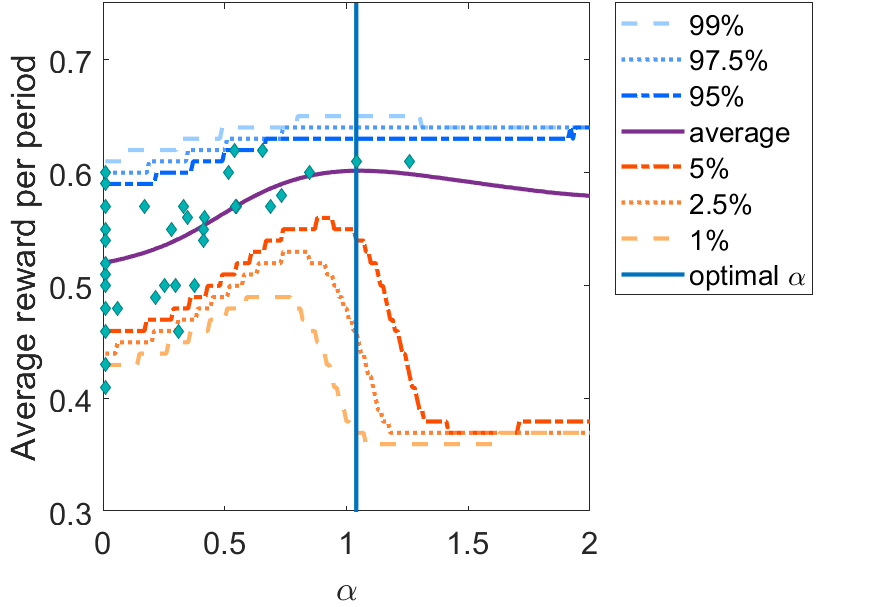}
		\subcaption{Remaining-period reminder}
	\end{subfigure}
	
	\caption{Individual-Level Behavior Patterns in the No-Summary-Table Extension and Targeted Robustness Checks}
	\label{fig:additional_indi_pattern}
\end{figure}

The average estimated \(\alpha\) was \(0.449\) in the no-summary-table extension and \(0.345\) in the control group of the targeted robustness checks described below; the difference was not statistically significant (\(p =0.401\)). We fit the same set of candidate behavioral models to this additional dataset. As shown in Table~\ref{tb:comp-predict-power-notable}, QCARE continues to perform well relative to the benchmark models. In particular, QCARE achieves the best out-of-sample average log-likelihood and outperforms the benchmark models for a large share of participants. This extension suggests that the predictive advantage of QCARE is not driven by the display of sufficient statistics on the interface. The behavioral pattern also remains consistent with over-exploration, as shown in Figure~\ref{fig:additional_indi_pattern}(c).

\begin{table}[htbp]
	\centering
		\caption{Comparison Between QCARE and Other Models (No-Summary-Table Extension)}
		\label{tb:comp-predict-power-notable}
		\centering
		\begin{tabular}{cccccccccc}
			\toprule
			\footnotesize{Data} & \footnotesize{Metric}  & QCARE   & L-$n_j$ & HH-$n_j$ & ES-$\gamma_j$ & QCARE-0 & $\varepsilon_j$-greedy & TS-B  & TS-G \\
			\midrule
			\multirow{2}{*}{IS}  & AvgLL  & -42.8      & -41.6     & -41.5    & -41.9     & -48.8\dstar    & -44.1\star  & -60.8\dstar & -49.1\dstar  \\
			& Pct & --   & 31.0\% & 37.9\% & 20.7\% & 100.0\%\dstar & 62.1\% & 100.0\%\dstar & 96.6\%\dstar \\
			\cmidrule{1-10}
			\multirow{2}{*}{OoS} & AvgLL   & -11.9    & -21.2    & -13.8\dstar     & -14.4\dstar    & -15.4\dstar      & -13.1\dstar  & -17.9\dstar & -15.8\dstar    \\
			& Pct &-- & 65.5\%\star & 65.5\%\star & 75.9\%\dstar & 79.3\%\dstar & 58.6\% & 89.7\%\dstar & 89.7\%\dstar \\
			\bottomrule
		\end{tabular}
		
		\vskip 0.1 cm
		\noindent\parbox{\linewidth}{\sf \scriptsize
			\textit{Note:} ``IS'' means in-sample (i.e., evaluated on the training data), and ``OoS'' means out-of-sample (i.e., evaluated on the test data). ``AvgLL'' refers to the average log-likelihood at the subject level by each approach. ``Pct'' indicates the corresponding percentage of subjects for whom QCARE outperforms the benchmark. The remarks ``**'' and ``*'' indicate $ p $-values less than $ 0.01 $ and $ 0.05 $, respectively.}
\end{table}

\begin{figure}[htbp]
	\centering
	\begin{subfigure}{0.7\textwidth}
		\includegraphics[width=\textwidth]{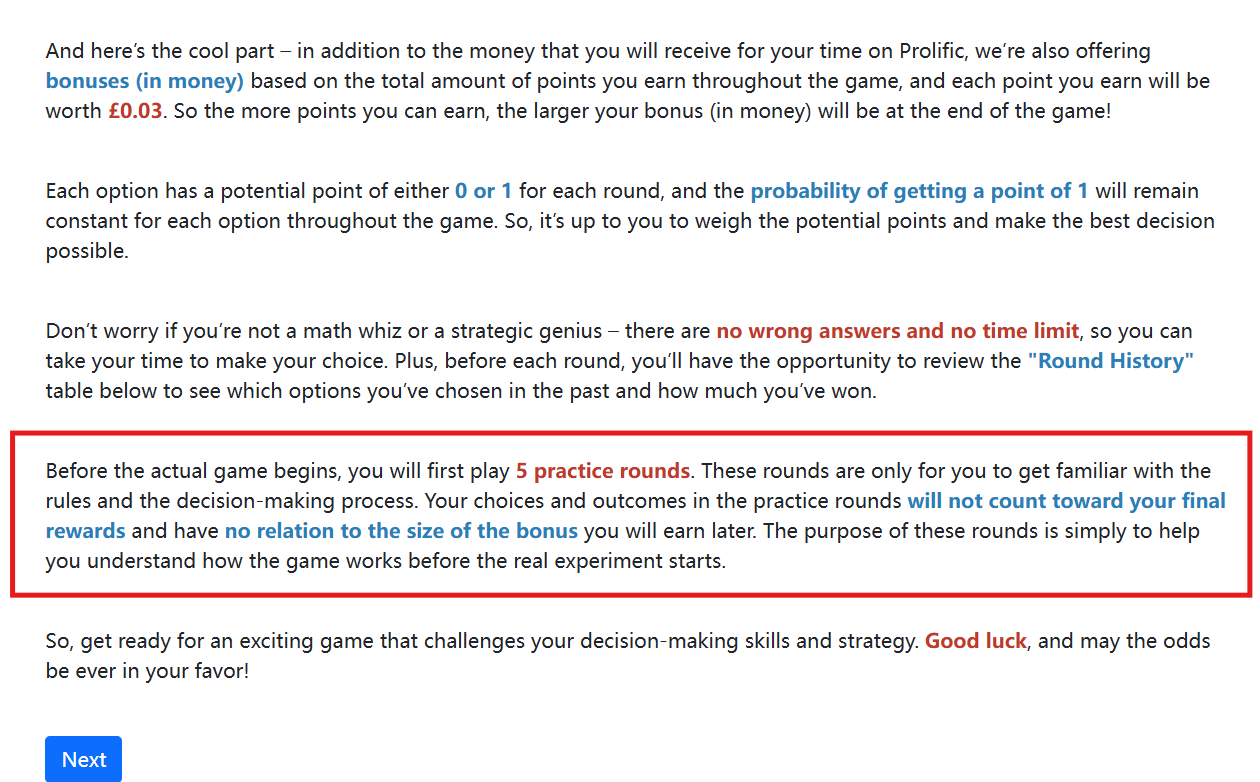}
		\subcaption{Practice-round page}
	\end{subfigure}
    
	\begin{subfigure}{0.45\textwidth}
		\includegraphics[width=\textwidth]{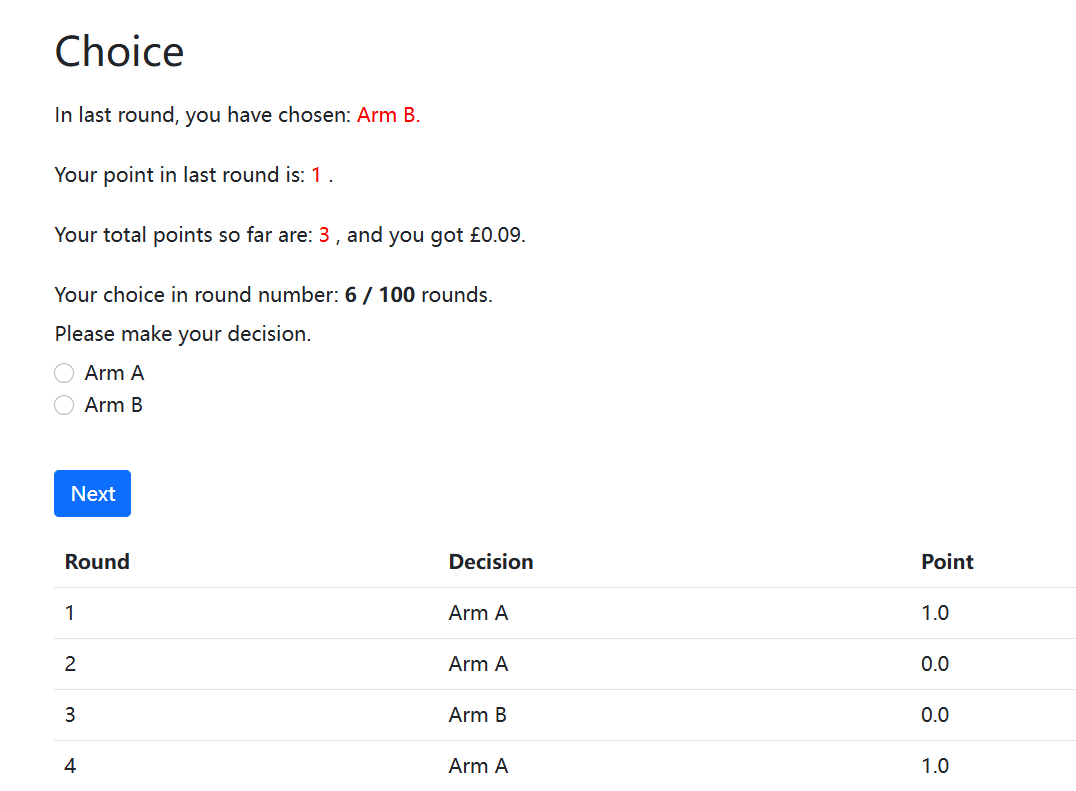}
		\subcaption{No-summary-table extension}
	\end{subfigure}
    \hspace{0.4 cm}
    \begin{subfigure}{0.45\textwidth}
		\includegraphics[width=\textwidth]{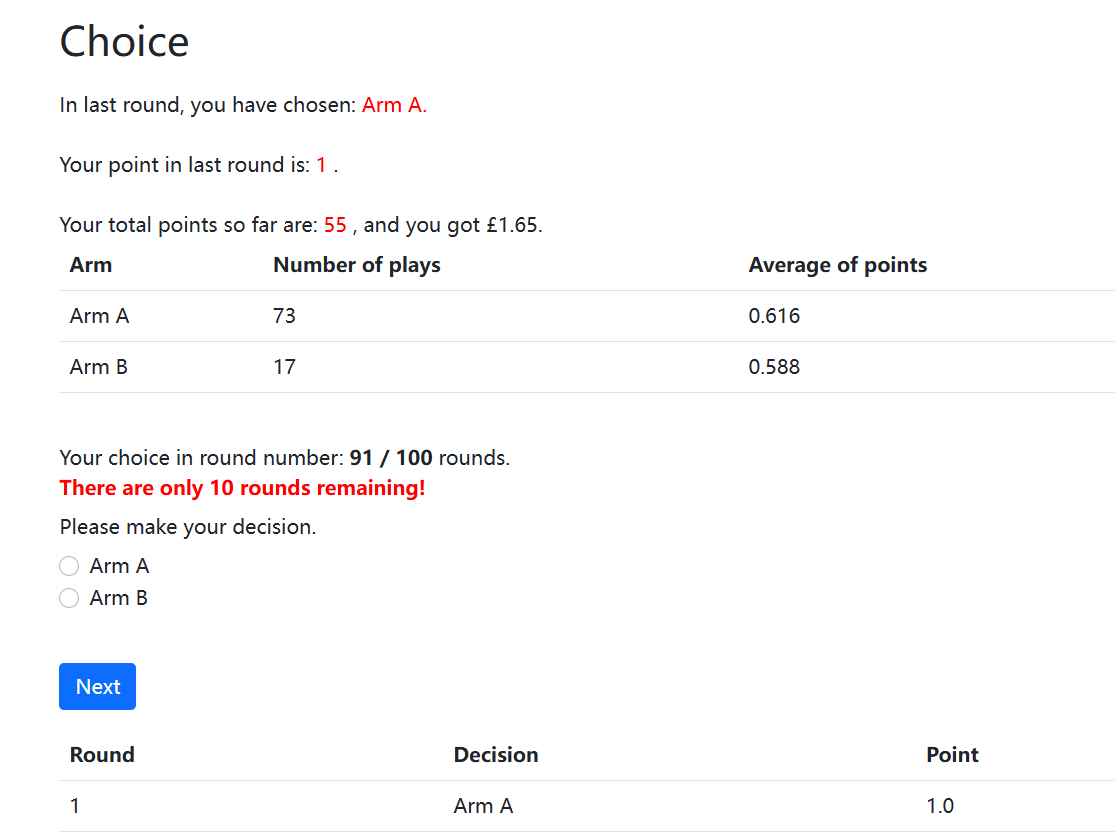}
		\subcaption{Remaining-period reminder}
	\end{subfigure}
	
	\caption{Screenshots From the No-Summary-Table Extension and Targeted Robustness Checks}
	\label{fig:screenshot_additional}
\end{figure}

\vskip 0.2 cm
\noindent\textbf{Targeted robustness checks.} We also conducted targeted robustness checks of the main behavioral pattern. These checks comprise a control group and three treatment conditions (practice rounds, higher stakes, and a remaining-period reminder) designed to test whether participants' behavior is sensitive to task familiarization, monetary stakes, or the salience of the remaining horizon. Due to budget constraints, all four groups used a single experimental configuration with horizon \(T=100\) and reward mean vector \((0.6,0.4)\).
Figure~\ref{fig:screenshot_additional} shows representative screenshots from these checks and the separate no-summary-table extension.

\paragraph{Practice rounds.} 
We first examine whether additional task familiarization affects participants' behavior. In this condition, participants completed five practice rounds before the formal experiment; they were explicitly informed that the practice-round outcomes would not affect their payoff in the main task. The results suggest that the practice rounds had little effect on formal-task behavior. Participants in the practice-round group earned an average of 54.925 total points, compared with 54.405 in the corresponding control group, and the difference is not statistically significant based on a two-sample \(t\)-test (\(p=0.5618\)). The average estimated \(\alpha\) was \(0.297\) in the practice-rounds condition and \(0.345\) in the control group; the difference was not statistically significant (\(p=0.563\)). As shown in Figure~\ref{fig:additional_indi_pattern}(a), the tendency toward over-exploration remains visible after the practice rounds.

\paragraph{Higher stakes.} 
We next examine whether increasing monetary incentives changes participants' behavior. In this condition, each point was converted to \pounds 0.06, which is twice the conversion rate in the original setting (\pounds 0.03). Participants in the higher-stakes condition earned fewer total points on average than those in the control group (52.225 versus 54.405), and the difference is statistically significant (\(p=0.0215\)). The average estimated \(\alpha\) was \(0.312\) in the higher-stakes condition and \(0.345\) in the control group; the difference was not statistically significant (\(p=0.721\)). Thus, increasing monetary stakes does not eliminate the over-exploration pattern; if anything, performance becomes lower in this setting. Figure~\ref{fig:additional_indi_pattern}(b) shows the corresponding individual-level pattern.

\paragraph{Remaining-period reminder.}
We finally examine whether participants' exploration behavior is partly driven by limited attention to the finite horizon. Participants were randomly assigned to either a reminder condition or a corresponding control condition. During periods 91--100, participants in the reminder condition were explicitly shown the number of periods remaining, whereas participants in the control condition completed the same task without this additional reminder. This comparison allows us to examine whether making the remaining horizon more salient shifts participants toward greater exploitation near the end of the task. The average estimated \(\alpha\) was \(0.296\) in the remaining-period-reminder condition and \(0.345\) in the control group; the difference was not statistically significant (\(p =0.576\)). Figure~\ref{fig:additional_indi_pattern}(d) shows the corresponding individual-level behavioral pattern.

We also estimate the following difference-in-differences linear probability model with participant and period fixed effects:
\begin{equation*}
    Y_{it}
    =
    \gamma_i
    +
    \lambda_t
    +
    \beta
    \bigl(
        \mathrm{Treat}_i
        \times
        \mathrm{Post}_t
    \bigr)
    +
    \varepsilon_{it}.
\end{equation*}
Here, \(Y_{it}\) equals one if participant \(i\) chooses the leading arm in period \(t\), and zero otherwise. The indicator \(\mathrm{Treat}_i\) equals one for participants assigned to the reminder condition and zero for those assigned to the control condition. The indicator \(\mathrm{Post}_t\) equals one for periods 91--100, during which the reminder was displayed, and zero for periods 1--90. The participant fixed effects \(\gamma_i\) absorb time-invariant individual heterogeneity, while the period fixed effects \(\lambda_t\) absorb changes in behavior that are common to both conditions in a given period. Accordingly, the main effects of \(\mathrm{Treat}_i\) and \(\mathrm{Post}_t\) are absorbed by these fixed effects. Standard errors are clustered at the participant level. The coefficient \(\beta\) measures whether the change in the probability of choosing the leading arm from periods 1--90 to periods 91--100 differs between the reminder and control conditions. The estimated interaction coefficient is
$
    \hat{\beta}
    =
    0.0049
$ 
$    \bigl(
        \text{clustered SE}=0.0536,\;
        t=0.091,\;
        p=0.928
    \bigr).
$
Thus, relative to the control condition, the reminder condition is associated with an estimated increase of approximately \(0.49\) percentage points in the probability of choosing the leading arm during the final ten periods. The estimate is close to zero and is not statistically significant. This robustness check therefore provides no evidence that making the remaining horizon more salient increases leading-arm choices near the end of the task.

Overall, the targeted robustness checks provide consistent evidence that the over-exploration pattern is robust across these design variations. Relative to the control group, participants continue to explore non-leading arms under practice rounds, higher stakes, and the remaining-period reminder.

\vskip 0.3 cm
\section{Proofs}
\label{section:appendix:D}
\setcounter{equation}{0}
\renewcommand{\theequation}{D.\arabic{equation}}
\setcounter{theorem}{0}
\renewcommand{\thetheorem}{D.\arabic{theorem}}

\vspace{0.3 cm}

\begin{proofref}{Proposition~\ref{prop:TS_QCARE_equa}}
We note that Gaussian Thompson Sampling draws samples independently from the posterior distribution $\mathcal{N}(\hat{\mu}_i, \frac{1}{k_i+1})$. QCARE with $\alpha = 0.5$ samples quantities $\theta_i(t)=\hat{\mu}_i(t)+\epsilon_{it}/\sqrt{k_i(t)+1}$, where $\epsilon_{it}\sim \mathcal{N}(0,1^2)$. By the property of normal distribution, these two policies are equivalent. 
\end{proofref}

\begin{proofref}{Proposition~\ref{prop:QCARE-basic-properties}} 
Let $ \kappa = (\kappa_1, \ldots, \kappa_N) $ and $ u = (u_1, \ldots, u_N) $ be given. Without loss, suppose $ \hat{i}^\ast = 1 $, i.e., $ u_1 \geq u_j $ for all $ j \neq 1 $.

(1). Recall that the arm-pulling score is such that $ \theta_{i} = u_{i} + \tfrac{\epsilon_{i}}{(\kappa_{i}+1)^\alpha} $ for every $ i \in [N] $, where $\epsilon_{i}$ follows the standard normal distribution independently.
Note that 
\begin{align*}
\hat{Q}(\kappa, u)  = \text{Pr}(\theta_1 \geq \max\{\theta_2, \ldots, \theta_N\})  = \text{Pr}_\epsilon \Big(\underbrace{u_{1} + \tfrac{\epsilon_{1}}{(\kappa_{1}+1)^\alpha} \geq \max\{u_2 + \tfrac{\epsilon_{2}}{(\kappa_{2}+1)^\alpha}, \ldots, u_N +\tfrac{\epsilon_{N}}{(\kappa_{N}+1)^\alpha} \}}_{E_1}\Big) 
\end{align*}
and
\begin{align*}
\hat{Q}(\kappa, u_1 +\delta , u_2, \ldots, u_N) =  \text{Pr}_\epsilon \Big(\underbrace{u_{1} + \tfrac{\epsilon_{1}}{(\kappa_{1}+1)^\alpha} + \delta \geq \max\{u_2 + \tfrac{\epsilon_{2}}{(\kappa_{2}+1)^\alpha}, \ldots, u_N + \tfrac{\epsilon_{N}}{(\kappa_{N}+1)^\alpha} \}}_{E_2}\Big) .
\end{align*}
Since $ u_{1} + \tfrac{\epsilon_{1}}{(\kappa_{1}+1)^\alpha} + \delta > u_{1} + \tfrac{\epsilon_{1}}{(\kappa_{1}+1)^\alpha}  $ with probability one, the proof is finished by noting the relationship between the following two events: $E_1 \subseteq E_2$.

(2).
Let \(F_{\theta_i}\) and \(f_{\theta_i}\) denote the cumulative distribution function and density of \(\theta_i\), respectively.
Note that 
	\begin{align*}
		 \hat{Q}(\kappa, u)  = & \text{Pr}(\theta_1 \geq \max\{\theta_2, \ldots, \theta_N\}) \\ 
		= & \text{Pr}_\epsilon \Big( \theta_1 \geq \theta_2, \ldots, \theta_1 \geq \theta_N
		\Big) \\
		= &  \int  \text{Pr}_\epsilon \Big(\theta_2\leq x, ..., \theta_N \leq x \Big) f_{\theta_1}(x) dx \\
		= &  \int  F_{\theta_2}(x)...F_{\theta_N}(x)  f_{\theta_1}(x) dx.
	\end{align*}
The last equality holds as $\theta_2, \ldots, \theta_N$ are independent of each other. Similarly, note that 
	\begin{align*}
		&\hat{Q}(\kappa_1, \ldots, \kappa_j+\delta, \ldots, \kappa_N, u_{1}, \ldots,  u_{N}) \\
		=  &\text{Pr}_\epsilon \Big(u_{1} + \tfrac{\epsilon_{1}}{(\kappa_{1}+1)^\alpha} \geq \max\{u_2 + \tfrac{\epsilon_{2}}{(\kappa_{2}+1)^\alpha}, \ldots, u_j +\tfrac{\epsilon_{j}}{(\kappa_{j} + \delta+1)^\alpha}, \ldots, u_N +\tfrac{\epsilon_{N}}{(\kappa_{N}+1)^\alpha} \}\Big) \\
		= & \int  F_{\theta_2}(x) \ldots F_{\theta'_j}(x) \ldots F_{\theta_N}(x)  f_{\theta_1}(x) dx,
	\end{align*}
where $ \theta'_j = u_{j} + \tfrac{\epsilon_{j}}{(\kappa_{j}+\delta+1)^\alpha} $.  Therefore, we have 
	\begin{align} \label{eq: diff_hat_q}
		&\hat{Q}(\kappa_1, \ldots, \kappa_j+\delta, \ldots, \kappa_N, u_{1}, \ldots,  u_{N})- \hat{Q}(\kappa, u) \nonumber \\
		=& \int  F_{\theta_2}(x)\ldots  F_{\theta'_j}(x)\ldots F_{\theta_N}(x)  f_{\theta_1}(x) dx-\int  F_{\theta_2}(x)\ldots F_{\theta_j}(x)\ldots F_{\theta_N}(x)  f_{\theta_1}(x) dx \nonumber \\
		= & \int  F_{\theta_2}(x) \ldots F_{\theta_{j-1}}(x)F_{\theta_{j+1}}(x)\ldots F_{\theta_N}(x) (F_{\theta'_j}(x)-F_{\theta_j}(x))  f_{\theta_1}(x) dx. 
	\end{align}
Next, we construct two functions to bound the terms $F_{\theta_2}(x) \ldots F_{\theta_{j-1}}(x)F_{\theta_{j+1}}(x) \ldots F_{\theta_N}(x)$ and $f_{\theta_1}(x)$.
We define two new functions $\chi(x)$ and $\phi(x)$ as follows: 
	\begin{equation*}
		\chi(x) =
		\left\{
		\begin{array}{ll}
			{F_{\theta_2}(x)...F_{\theta_{j-1}}(x)F_{\theta_{j+1}}(x)...F_{\theta_N}(x),} & { \textrm{if $x\leq u_j$}}, \\
			{F_{\theta_2}(2u_j-x)...F_{\theta_{j-1}}(2u_j-x)F_{\theta_{j+1}}(2u_j-x)...F_{\theta_N}(2u_j-x),} & { \textrm{if $x > u_j$}},
		\end{array}
		\right.
	\end{equation*}
    and
	\begin{equation*}
		\phi(x) =
		\left\{
		\begin{array}{ll}
			{f_{\theta_1}(x),} & { \textrm{if $x\leq u_j$}}, \\
			{f_{\theta_1}(2u_j-x),} & { \textrm{if $x > u_j$}}.
		\end{array}
		\right.
	\end{equation*}
It is evident that both $\chi(x)$ and $\phi(x)$ are axially symmetric about the line $x=u_j$. As CDFs are non-decreasing, when $x>u_j$, we have $F_{\theta_i}(2u_j-x)<F_{\theta_i}(u_j)<F_{\theta_i}(x)$ for all $ i $. Similarly, as $u_1>u_j$, it holds that $f_{\theta_1}(x)\geq f_{\theta_1}(2u_j-x)$ when $x>u_j$. Then, we have $F_{\theta_2}(x) \ldots F_{\theta_{j-1}}(x)F_{\theta_{j+1}}(x) \ldots F_{\theta_N}(x) \geq \chi(x)$ and $f_{\theta_1}(x) \geq \phi(x)$.
Thus, by Equation~\eqref{eq: diff_hat_q},
	\begin{align*}
		&\hat{Q}(\kappa_1, \ldots, \kappa_j+\delta, \ldots, \kappa_N, u_{1}, \ldots,  u_{N})- \hat{Q}(\kappa, u) \\
		= & \int_{-\infty}^{u_j}  F_{\theta_2}(x) \ldots F_{\theta_{j-1}}(x)F_{\theta_{j+1}}(x)\ldots F_{\theta_N}(x) (F_{\theta'_j}(x)-F_{\theta_j}(x))  f_{\theta_1}(x) dx \\ 
        & + \int_{u_j}^{\infty}  F_{\theta_2}(x) \ldots F_{\theta_{j-1}}(x)F_{\theta_{j+1}}(x)\ldots F_{\theta_N}(x) (F_{\theta'_j}(x)-F_{\theta_j}(x))  f_{\theta_1}(x) dx  \\
		\geq & \int_{-\infty}^{u_j} \chi(x)  (F_{\theta'_j}(x)-F_{\theta_j}(x)) \phi(x) dx + \int_{u_j}^{\infty} \chi(x)  (F_{\theta'_j}(x)-F_{\theta_j}(x)) \phi(x) dx \\
		\geq & 0.
	\end{align*}
The first inequality holds since when $x \leq u_j$, we get $F_{\theta_2}(x) \ldots F_{\theta_{j-1}}(x)F_{\theta_{j+1}}(x) \ldots F_{\theta_N}(x) = \chi(x)$ and $f_{\theta_1}(x) = \phi(x)$ by definitions. When $x > u_j$, as $F_{\theta'_j}(x)-F_{\theta_j}(x)>0$ and the facts that $F_{\theta_2}(x) \ldots F_{\theta_{j-1}}(x)F_{\theta_{j+1}}(x) \ldots F_{\theta_N}(x) \geq \chi(x)$ and $f_{\theta_1}(x) \geq \phi(x)$, the inequality holds.
The last inequality holds as $\chi(x)$ and $\phi(x)$ are both axially symmetric about the line $x=u_j$. We know that $\chi(x)\geq 0$, $\phi(x)> 0$. Figure~\ref{fig:proof_cdf} plots the CDFs of $\theta_j$ and $\theta'_j$, and the difference term $F_{\theta'_j}(x)-F_{\theta_j}(x)$ is highlighted in color. We can see that the term $F_{\theta'_j}(x)-F_{\theta_j}(x)$ is centrally symmetric about point $(u_j,\frac{1}{2})$, and $F_{\theta'_j}(x)-F_{\theta_j}(x)$ is negative when $x<u_j$ and positive when $x>u_j$.  The blue part and yellow part are symmetric about the point $(u_j,\frac{1}{2})$. Thus, we have $\hat{Q}(\kappa_1, \ldots, \kappa_j+\delta, \ldots, \kappa_N, u_{1}, \ldots,  u_{N})- \hat{Q}(\kappa, u) \geq 0$.
    \begin{figure}[htbp]
		\centering
		\includegraphics[width=0.5 \textwidth]{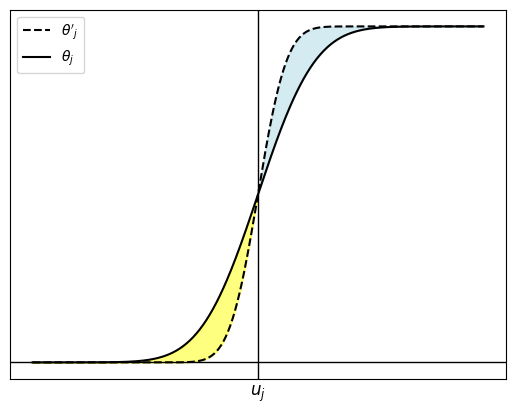}
		\caption{CDFs of $\theta_j$ and $\theta'_j$}
		\label{fig:proof_cdf}
    \end{figure}
\end{proofref}

\subsection{Proof of Theorem~\ref{thm:smaller0.5}}\label{appendix:upper-bound-alpha-lt-half}

In this subsection, we perform regret analysis that is customized to the QCARE policy. Fix \(\alpha\in(0,0.5)\) and \(N>1\) as required by the theorem. Theorem~\ref{thm:smaller0.5} is a direct consequence of the combination of the propositions below.

\begin{proposition}\label{prop:ub_smaller0.5_upper}
	There exists \(C_{\alpha,N} \in (0, \infty)\) such that, for every instance \(\Theta\) and \(T \geq 1\), 
	\begin{equation*}
	\mathcal{R}^{\Theta}(T)
	\,\leq\,
	C_{\alpha,N} \big(\log(eNT)\big)^{1/(2\alpha)}
	\sum_{i:\Delta_i>0}
	\Delta_i^{1-1/\alpha}.
	\end{equation*}
	In addition, $ \mathcal{R}(T)
	\le
	C_{\alpha,N}T^{1-\alpha}\sqrt{\log(eT)} = \widetilde{O}(T^{1-\alpha}) $.
\end{proposition}

\begin{proposition}\label{prop:ub_smaller0.5_lower_instance}
	There exists \(c_{\alpha,N} \in (0, \infty)\) such that for every instance \(\Theta\) and sufficiently large $ T $, 
	\begin{equation*}
	\mathcal{R}^{\Theta}(T) \,\geq\, c_{\alpha,N}
	\bigl(\log(eNT)\bigr)^{1/(2\alpha)}
	\sum_{i:\Delta_i>0}
	\Delta_i^{1-1/\alpha}.
	\end{equation*}
\end{proposition}

\begin{proposition}\label{prop:ub_smaller0.5_lower_worst}
	$
	\mathcal{R}(T) = \Omega(T^{1-\alpha}).
	$
\end{proposition}

In the analysis below, we provide proofs of Propositions~\ref{prop:ub_smaller0.5_upper} and~\ref{prop:ub_smaller0.5_lower_instance}, where we directly exploit the structures of QCARE. For Proposition~\ref{prop:ub_smaller0.5_lower_worst}, we are able to extend it to a general class of Markovian MAB policies. We thus defer its analysis and proof to Appendix~\ref{appendix:formal-general-framework}.

\vskip 0.3 cm
\textbf{Proof of Proposition~\ref{prop:ub_smaller0.5_upper}.} Recall that \(i^*\), \(\Delta_i\), \(k_i(t)\), \(\hat\mu_i(t)\), and the QCARE score rule are defined in the main text. Relabel the arms so that \(i^\ast=1\), then \(\Delta_i=\mu_1-\mu_i>0\) for \(i=2,\ldots,N\).
For the proof, let \(X_{i,s}\) be the \(s\)-th potential reward of arm \(i\), so that for each arm \(i\), the sequence \(\{X_{i,s}:s\ge1\}\) is an i.i.d. sequence supported on \([0,1]\) with mean \(\mu_i\), and these sequences are independent across arms. Write \(\hat\mu_{i,s}:=\sum_{\ell=1}^{s}X_{i,\ell}/(s+1)\), with \(\hat\mu_{i,0}=0\). Thus, when arm \(i\) has been pulled \(k_i(t)\) times, the QCARE score can be written as \(\theta_i(t)=\hat\mu_{i,k_i(t)}+(k_i(t)+1)^{-\alpha}\epsilon_{i,t}\).

The proof has three steps. First, we obtain a uniform concentration envelope for the shrunken empirical means. Second, we show that the optimal arm cannot be delayed for too long. Third, we show that, after both the optimal arm and a fixed suboptimal arm have been sufficiently sampled, the probability of selecting that suboptimal arm is negligible.

\begin{lemma}[Uniform concentration envelope]
\label{lem:uniform_concentration_envelope_alpha_lt_half}
For each arm \(i\), define
\[
\xi_i
=
\sup_{0\le s\le T}
(s+1)^\alpha\left|\hat\mu_{i,s}-\mu_i\right|.
\]
There exist constants \(b_\alpha>1/2\) and \(C_\alpha<\infty\), depending only on \(\alpha\), such that
$\mathbb E\left[\exp\left(b_\alpha \xi_i^2\right)\right]\le C_\alpha $
uniformly over \(T\), \(i\), and the reward distribution of arm \(i\).
\end{lemma}

\begin{proofref}{Lemma~\ref{lem:uniform_concentration_envelope_alpha_lt_half}}
For \(s\ge1\), let \(\bar X_{i,s}:=\frac{1}{s}\sum_{\ell=1}^s X_{i,\ell}\).
Since \(X_{i,\ell}\in[0,1]\), we have
\[
\left|\hat\mu_{i,s}-\bar X_{i,s}\right|
=
\left|\frac{s}{s+1}\bar X_{i,s}-\bar X_{i,s}\right|
\le \frac{1}{s+1}.
\] 
Fix \(x\ge4\). The case \(s=0\) is harmless because
\( (s+1)^\alpha|\hat\mu_{i,s}-\mu_i|\le 1, \text{when }s=0. \)
If, for some \(s\ge1\),
\( (s+1)^\alpha|\hat\mu_{i,s}-\mu_i|>x, \)
then
\[
{\color{black}
|\bar X_{i,s}-\mu_i|
>
\frac{x}{(s+1)^\alpha}-\frac{1}{s+1}
\ge
\frac{3x}{4(s+1)^\alpha},
}
\]
where the last inequality uses \(x\ge4\) and \(\alpha<1\). By Hoeffding's inequality,
\[
{\color{black}
\mathbb P\left(
(s+1)^\alpha|\hat\mu_{i,s}-\mu_i|>x
\right)
\le
2\exp\left\{
-2s\left(\frac{3x}{4(s+1)^\alpha}\right)^2
\right\}.
}
\]
Since \(s+1\le 2s\) for \(s\ge1\), and since \(\beta:=1-2\alpha>0\),
there exists a constant \(\mathcal{A}_\alpha>1/2\), depending only on \(\alpha\), such that
\( \mathbb P\left( (s+1)^\alpha|\hat\mu_{i,s}-\mu_i|>x \right) \le 2\exp\{-\mathcal{A}_\alpha x^2s^\beta\}. \)
Taking a union bound over \(s\ge1\) and comparing the resulting series with an integral, we obtain constants \(\widetilde C_\alpha<\infty\) and \(c_\alpha>1/2\) such that
\[
{\color{black}
\mathbb P(\xi_i>x)
\le
\widetilde C_\alpha \exp\{-c_\alpha x^2\},
\qquad x\ge0.
}
\]
The extension from \(x\ge4\) to all \(x\ge0\) is obtained by increasing \(\widetilde C_\alpha\). Choosing any \(b_\alpha\in(1/2,c_\alpha)\) and integrating the tail of \(\exp(b_\alpha \xi_i^2)\) gives
\( \mathbb E\left[\exp(b_\alpha \xi_i^2)\right]\le C_\alpha. \)
This proves the claim.
\end{proofref}

For the next step, condition on the reward table up to horizon \(T\), \(\mathcal X_T:=\{X_{i,\ell}:1\le i\le N,\ 1\le \ell\le T\}\).
Let \(Q_T:=\prod_{i=1}^N\Phi\bigl(-(\xi_i+1)\bigr)\), where \(\Phi\) is the standard normal distribution function.

\begin{lemma}[Entry into concentration]
\label{lem:entry_into_concentration_alpha_lt_half}
There exists a constant \(C^{\mathrm{ent}}_{\alpha,N}<\infty\) such that, for every \(u\in\{1,\ldots,T\}\),
\[
{\color{black}
\mathbb E\left[
\sum_{t=1}^T
\mathbf 1\{a(t)\neq 1,\ k_1(t)<u\}
\right]
\le
C^{\mathrm{ent}}_{\alpha,N}u.
}
\]
\end{lemma}

\begin{proofref}{Lemma~\ref{lem:entry_into_concentration_alpha_lt_half}}
Let \(\mathcal F_{t-1}\) denote the sigma-field generated by the observed actions and rewards through period \(t-1\), together with the policy randomization used through period \(t-1\):
\[
{\color{black}
\mathcal F_{t-1}:=\sigma\!\left(\left\{a(s),r_{a(s)}(s),\epsilon_{i,s}:i\in[N],\ 1\le s\le t-1\right\}\right).
}
\]
We first prove a forcing bound. For every \(t\le T\),
\( \mathbb P\left(a(t)=1\mid \mathcal F_{t-1},\mathcal X_T\right)
\ge Q_T. \)
Let \(m_i=k_i(t)\). Consider the event
\[
{\color{black}
E_t
=
\left\{
\epsilon_{1,t}\ge \xi_1+1,\ 
\epsilon_{j,t}\le -(\xi_j+1)\ \text{for all }j\ge2
\right\}.
}
\]
Conditional on \(\mathcal F_{t-1}\vee\sigma(\mathcal X_T)\), the current Gaussian shocks are still independent standard normal random variables. Therefore,
\( \mathbb P(E_t\mid \mathcal F_{t-1},\mathcal X_T)=Q_T. \)
On \(E_t\), by the definition of \(\xi_1\),
\[
{\color{black}
\theta_1(t)
\ge
\hat\mu_{1,m_1}+(m_1+1)^{-\alpha}(\xi_1+1)
\ge
\mu_1+(m_1+1)^{-\alpha}.
}
\]
Similarly, for every \(j\ge2\),
\[
{\color{black}
\theta_j(t)
\le
\hat\mu_{j,m_j}-(m_j+1)^{-\alpha}(\xi_j+1)
\le
\mu_j-(m_j+1)^{-\alpha}
<
\mu_1.
}
\]
Thus \(E_t\subseteq\{a(t)=1\}\), which proves the forcing bound.

Next we show that \(Q_T^{-1}\) has bounded expectation. A standard Gaussian lower-tail bound implies that, for any \(b'_\alpha>1/2\), after increasing the multiplicative constant if necessary,
\( \Phi(-(z+1))^{-1} \le C_\alpha \exp\{b'_\alpha z^2\}, z\ge0. \)
Choose \(b'_\alpha\in(1/2,b_\alpha)\), where \(b_\alpha\) is the constant in Lemma~\ref{lem:uniform_concentration_envelope_alpha_lt_half}. Since the reward tables are independent across arms, the random variables \(\xi_1,\ldots,\xi_N\) are independent. Hence,
\[
{\color{black}
\mathbb E[Q_T^{-1}]
\le
C_{\alpha,N}
\prod_{i=1}^N
\mathbb E\left[\exp\{b'_\alpha \xi_i^2\}\right]
\le
C_{\alpha,N}
<\infty.
}
\]

Now fix \(m\in\{0,\ldots,u-1\}\), and define
\( V_m = \sum_{t=1}^T \mathbf 1\{k_1(t)=m,\ a(t)\neq1\}. \)
Conditional on \(\mathcal X_T\), whenever the process is in a period with \(k_1(t)=m\), the probability of pulling arm \(1\) is at least \(Q_T\). Once arm \(1\) is pulled, the process leaves the state \(k_1(t)=m\) forever. Therefore, the number of failures before leaving this state is stochastically dominated by a geometric waiting time with success probability \(Q_T\). Equivalently, for every integer \(r\ge0\),
\( \mathbb P(V_m\ge r\mid \mathcal X_T) \le (1-Q_T)^r. \)
Thus
\( \mathbb E[V_m\mid \mathcal X_T]\le Q_T^{-1}, \)
and consequently
\( \mathbb E[V_m]\le C_{\alpha,N}. \)
Summing over \(m=0,\ldots,u-1\), we obtain
\[
{\color{black}
\mathbb E\left[
\sum_{t=1}^T
\mathbf 1\{a(t)\neq 1,\ k_1(t)<u\}
\right]
\le
C^{\mathrm{ent}}_{\alpha,N}u.
}
\]
This completes the proof.
\end{proofref}

We now control the residual probability of selecting a fixed suboptimal arm after enough samples have been collected.

\begin{lemma}[Residual exploration after concentration]
\label{lem:residual_exploration_alpha_lt_half}
There exists a constant \(\mathcal{A}_\alpha<\infty\), depending only on \(\alpha\), such that the following holds for every arm \(i\ge2\) with \(\Delta_i>0\). Define
\[
{\color{black}
n_i
=
\left\lceil
\mathcal{A}_\alpha
\left(
\tfrac{\log(eNT)}{\Delta_i^2}
\right)^{1/(2\alpha)}
\right\rceil.
}
\]
If \(n_i<T\), let \(G_i:=\{\xi_1\le \frac{\Delta_i}{4}(n_i+1)^\alpha,\ \xi_i\le \frac{\Delta_i}{4}(n_i+1)^\alpha\}\).
After choosing \(\mathcal{A}_\alpha\) sufficiently large,
\( \mathbb P(G_i^c)\le (eNT)^{-10}, \)
and on \(G_i\), whenever \(k_1(t)\ge n_i\) and \(k_i(t)\ge n_i\),
\[
{\color{black}
\mathbb P\left(a(t)=i\mid \mathcal F_{t-1},\mathcal X_T\right)
\le
(eNT)^{-10}.
}
\]
\end{lemma}

\begin{proofref}{Lemma~\ref{lem:residual_exploration_alpha_lt_half}}
On \(G_i\), for every \(s\ge n_i\),
$
|\hat\mu_{1,s}-\mu_1|
\le
\frac{\xi_1}{(s+1)^\alpha}
\le
\frac{\Delta_i (n_i+1)^\alpha}{4(s+1)^\alpha}
\le
\frac{\Delta_i}{4},
$
and similarly,
\( |\hat\mu_{i,s}-\mu_i| \le \frac{\Delta_i}{4}. \)
Therefore, whenever \(k_1(t)\ge n_i\) and \(k_i(t)\ge n_i\),
\( \hat\mu_{1,k_1(t)}-\hat\mu_{i,k_i(t)} \ge
\mu_1-\mu_i-\frac{\Delta_i}{2}
= \frac{\Delta_i}{2}. \)

We next bound \(\mathbb P(G_i^c)\). By Lemma~\ref{lem:uniform_concentration_envelope_alpha_lt_half}, there exist constants \(\widetilde C_\alpha<\infty\) and \(c_\alpha>0\) such that \( \mathbb P(\xi_j>x)\le \widetilde C_\alpha e^{-c_\alpha x^2}, j=1,i. \)
Moreover, by the definition of \(n_i\),
\[
{\color{black}
\frac{\Delta_i^2}{16}(n_i+1)^{2\alpha}
\ge
c \mathcal{A}_\alpha^{2\alpha}\log(eNT)
}
\]
for a numerical constant \(c>0\). Choosing \(\mathcal{A}_\alpha\) sufficiently large gives
\( \mathbb P(G_i^c) \le \mathbb P(\xi_1>\frac{\Delta_i}{4}(n_i+1)^\alpha)+\mathbb P(\xi_i>\frac{\Delta_i}{4}(n_i+1)^\alpha) \le (eNT)^{-10}. \)

It remains to control the conditional probability of selecting arm \(i\) on \(G_i\). Conditional on \(\mathcal F_{t-1}\vee\sigma(\mathcal X_T)\),
\( \theta_i(t)-\theta_1(t) \sim N\left( \hat\mu_{i,k_i(t)}-\hat\mu_{1,k_1(t)},\
(k_i(t)+1)^{-2\alpha}+(k_1(t)+1)^{-2\alpha} \right). \)
On \(G_i\), if \(k_1(t)\ge n_i\) and \(k_i(t)\ge n_i\), then
\( \hat\mu_{1,k_1(t)}-\hat\mu_{i,k_i(t)} \ge \frac{\Delta_i}{2}. \)
Hence,
\[
\mathbb P\left(a(t)=i\mid \mathcal F_{t-1},\mathcal X_T\right)
\le
\mathbb P\left(\theta_i(t)\ge \theta_1(t)\mid \mathcal F_{t-1},\mathcal X_T\right)
\le
\Phi\left(
-
\tfrac{\Delta_i/2}{
\sqrt{(k_i(t)+1)^{-2\alpha}+(k_1(t)+1)^{-2\alpha}}}
\right)
\le
\Phi\left(-c'\Delta_i(n_i+1)^\alpha\right),
\]
for some numerical constant \(c'>0\). By the definition of \(n_i\),
\( \Delta_i(n_i+1)^\alpha \ge c'' \mathcal{A}_\alpha^\alpha\sqrt{\log(eNT)} \)
for another numerical constant \(c''>0\). Choosing \(\mathcal{A}_\alpha\) large enough and using the Gaussian upper-tail bound, we obtain
\( \mathbb P\left(a(t)=i\mid \mathcal F_{t-1},\mathcal X_T\right) \le (eNT)^{-10}. \)
This completes the proof.
\end{proofref}

\medskip

We now prove Proposition~\ref{prop:ub_smaller0.5_upper}. Fix an arm \(i\ge2\) with \(\Delta_i>0\). If \(n_i\ge T\), since \(n_i=\lceil\mathcal A_\alpha \left(\frac{\log(eNT)}{\Delta_i^2}\right)^{1/(2\alpha)}\rceil\), we have
\[
\mathbb E[k_i(T+1)]
\le T
\le n_i
\le 1+\mathcal A_\alpha \left(\frac{\log(eNT)}{\Delta_i^2}\right)^{1/(2\alpha)}
\le C_{\alpha,N}
\left(
1+\left(\tfrac{\log(eNT)}{\Delta_i^2}\right)^{1/(2\alpha)}
\right)
\]
after increasing \(C_{\alpha,N}\), if necessary, so that
\(C_{\alpha,N}\ge\max\{1,\mathcal A_\alpha\}\).
If $ n_i < T $, decompose
\( k_i(T+1) \le N_i^{(1)}+N_i^{(2)}+N_i^{(3)}, \)
where
\[
N_i^{(1)}
=
\sum_{t=1}^T
\mathbf 1\{a(t)=i,\ k_i(t)<n_i\},
\]
\[
N_i^{(2)}
=
\sum_{t=1}^T
\mathbf 1\{a(t)=i,\ k_i(t)\ge n_i,\ k_1(t)<n_i\}, \quad \text{ and } \quad 
N_i^{(3)}
=
\sum_{t=1}^T
\mathbf 1\{a(t)=i,\ k_i(t)\ge n_i,\ k_1(t)\ge n_i\}.
\]
The first term is bounded pathwise by \( N_i^{(1)}\le n_i. \)
For the second term,
\( N_i^{(2)} \le \sum_{t=1}^T \mathbf 1\{a(t)\neq1,\ k_1(t)<n_i\}. \)
Therefore, by Lemma~\ref{lem:entry_into_concentration_alpha_lt_half},
\( \mathbb E[N_i^{(2)}] \le C^{\mathrm{ent}}_{\alpha,N}n_i. \)
For the third term, use the event \(G_i\) from Lemma~\ref{lem:residual_exploration_alpha_lt_half}. On \(G_i\), the conditional probability of selecting arm \(i\) in each relevant period is at most \((eNT)^{-10}\). On \(G_i^c\), we use the trivial bound \(N_i^{(3)}\le T\). Hence,
\( \mathbb E[N_i^{(3)}] \le T(eNT)^{-10}+T\mathbb P(G_i^c) \le 2. \)
Combining the three bounds gives
\( \mathbb E[k_i(T+1)] \le C_{\alpha,N}(1+n_i). \)
Since \(n_i=\lceil \mathcal{A}_\alpha \left(\frac{\log(eNT)}{\Delta_i^2}\right)^{1/(2\alpha)} \rceil\), that translates to
\[
{\color{black}
\mathbb E[k_i(T+1)]
\le
C_{\alpha,N}
\left(
1+
\left(
\tfrac{\log(eNT)}{\Delta_i^2}
\right)^{1/(2\alpha)}
\right).
}
\]

We now convert this pull-count bound into a regret bound. Since rewards are supported on \([0,1]\), we have \(0<\Delta_i\le1\). Also,
$
\left(
\frac{\log(eNT)}{\Delta_i^2}
\right)^{1/(2\alpha)}
\ge1.
$
Thus the additive \(1\) can be absorbed after multiplying by \(\Delta_i\). That is, 
$
\Delta_i
\left(
1 + \left(
\tfrac{\log(eNT)}{\Delta_i^2}
\right)^{1/(2\alpha)}
\right)
\le
2(\log(eNT))^{1/(2\alpha)}
\Delta_i^{1-1/\alpha}.
$ Therefore,
\[
{\color{black}
\mathcal R^\Theta(T)
=
\sum_{i:\Delta_i>0}\Delta_i\,\mathbb E[k_i(T+1)]
\le
C_{\alpha,N}
\sum_{i:\Delta_i>0}  (\log(eNT))^{1/(2\alpha)}
\Delta_i^{1-1/\alpha}
}
\]
with an appropriate rescaling of $ C_{\alpha, N} $.
This proves the instance-dependent bound.

It remains to derive the instance-independent bound. Invoking the trivial bound $ \Delta_i \mathbb E[k_i(T+1)] \le \Delta_i T$, the instance-dependent regret upper bound can be strengthened to
\[
\mathcal R^\Theta(T)
		=
		\sum_{i:\Delta_i>0}\Delta_i\,\mathbb E[k_i(T+1)]
		\le
		C_{\alpha,N}
		\sum_{i:\Delta_i>0} \min\{\Delta_i T, (\log(eNT))^{1/(2\alpha)}
		\Delta_i^{1-1/\alpha}\}.
	\]
Let \(L_T:=(\log(eNT))^{1/(2\alpha)}\) and define
\( f_T(\Delta) = \min\left\{ T\Delta,\ L_T\Delta^{1-1/\alpha} \right\}. \) The first branch \(T\Delta\) is increasing in \(\Delta\), whereas the second branch \(L_T\Delta^{1-1/\alpha}\) is decreasing because \(1-1/\alpha<0\). Hence the maximum of \(f_T\) is achieved, up to constants, at the balancing point satisfying
\( T\Delta=L_T\Delta^{1-1/\alpha}. \)
This gives
\( \Delta_T=\left(\frac{L_T}{T}\right)^\alpha. \)
At this point,
\[
{\color{black}
f_T(\Delta_T)
=
T\Delta_T
=
T^{1-\alpha}L_T^\alpha
=
T^{1-\alpha}\sqrt{\log(eNT)}.
}
\]
Therefore, for all \(\Delta>0\),
\( f_T(\Delta) \le C_\alpha T^{1-\alpha}\sqrt{\log(eNT)}. \)
Since \(N\) is fixed, summing over at most \(N-1\) suboptimal arms and absorbing fixed \(N\)-dependent factors into the constant gives
\( \mathcal R(T) \le C_{\alpha,N}T^{1-\alpha}\sqrt{\log(eT)}. \)
This completes the proof of Proposition~\ref{prop:ub_smaller0.5_upper}.
\hfill $\Box$

\vskip 0.3 cm
\begin{proofref}{Proposition~\ref{prop:ub_smaller0.5_lower_instance}}
Fix an MAB instance \(\Theta\). Use the reward-table representation introduced above, in which \(\{X_{i,s}\}_{s\ge1}\) is the i.i.d. reward sequence of arm \(i\), independent across arms, and
\(\hat{\mu}_{i,s}:=\sum_{\ell=1}^sX_{i,\ell}/(s+1)\). We extend the horizon-independent QCARE rule indefinitely and evaluate its first \(T\) periods. Let \(\mathcal X_\infty:=\{X_{j,s}:j\in[N],\ s\ge1\}\) denote the complete reward table, which is independent of the Gaussian shocks used by QCARE.

We first record a statewise lower bound on the QCARE choice probability. Fix a state \(S=(\kappa,u)\) and an arm \(i\), and write \(Q_i(S)=\mathbb P\{a(t)=i\mid S(t)=S\}\). Recall that \(\epsilon_{1,t},\ldots,\epsilon_{N,t}\) are independent standard normal shocks used to generate the QCARE scores at this state. On the event
\(\epsilon_{i,t}\ge(\max_j u_j-u_i)(\kappa_i+1)^\alpha\) and \(\epsilon_{j,t}\le0\) for every \(j\ne i\), the score of arm \(i\) is at least \(\max_j u_j\), whereas every other arm has score at most its empirical mean. Therefore,
\begin{equation}
\label{eq:qcare-statewise-lower-bound}
{\color{black}
Q_i(S)
\ge
2^{-(N-1)}\bar\Phi\!\left((\max_j u_j-u_i)(\kappa_i+1)^\alpha\right),
}
\end{equation}
where \(\bar\Phi(x):=1-\Phi(x)\). Ties occur with probability zero because the Gaussian shocks are continuous. 

Fix a suboptimal arm \(i\ne i^\ast\). For every integer \(s_0\ge1\), define
\[
{\color{black}
\mathcal C_i(s_0):=
\left\{
|\hat{\mu}_{j,s}-\mu_j|\le \frac{\Delta_i}{4}
\text{ for every }j\in[N]\text{ and }s\ge s_0
\right\}.
}
\]
Since \(\hat{\mu}_{j,s}\to\mu_j\) almost surely for every arm \(j\) and \(N\) is finite, \(\mathbb P\bigl(\bigcup_{s_0\ge1}\mathcal C_i(s_0)\bigr)=1\). Hence, there exists a finite integer \(s_i^\circ\ge1\) such that, after abbreviating \(\mathcal C_i:=\mathcal C_i(s_i^\circ)\),
\(\mathbb P(\mathcal C_i)\ge1/2\). On \(\mathcal C_i\), whenever every arm has been pulled at least \(s_i^\circ\) times,
\begin{equation}
\label{eq:empirical-gap-upper-bound}
{\color{black}
\max_j\hat{\mu}_j(t)-\hat{\mu}_i(t)
\le
\left(\mu^*+\frac{\Delta_i}{4}\right)
-
\left(\mu_i-\frac{\Delta_i}{4}\right)
=\frac{3\Delta_i}{2}.
}
\end{equation}

Let \(H_i^{\mathrm{ent}}\) be the first period at which every arm has been pulled at least \(s_i^\circ\) times. The process reaches this concentration region in finite expected time. Indeed, if an arm \(j\) has been pulled fewer than \(s_i^\circ\) times, then \(\max_\ell u_\ell-u_j\le1\), and \eqref{eq:qcare-statewise-lower-bound} gives
\[
{\color{black}
Q_j(S)
\ge
p_i^{\mathrm{ent}}
:=
2^{-(N-1)}\bar\Phi\!\left((s_i^\circ)^\alpha\right)
>0.
}
\]
Conditional on any realized reward table, successively waiting for the at most \(Ns_i^\circ\) required pulls is therefore stochastically dominated by a sum of geometric waiting times with success probability \(p_i^{\mathrm{ent}}\). In particular,
\[
{\color{black}
\mathbb E[H_i^{\mathrm{ent}}\mid\mathcal X_\infty]
\le
1+\frac{Ns_i^\circ}{p_i^{\mathrm{ent}}}
<\infty.
}
\]
Because \(\mathcal C_i\) is determined by the reward table, the same bound holds conditional on \(\mathcal C_i\).

Let \(\Lambda_T:=\log(eNT)\) and set
$
{\color{black}
\rho_\alpha:=\left(\frac{4}{9}\right)^{1/(2\alpha)}.
}
$
Define the target pull count
\begin{equation}
\label{eq:target-pull-count}
{\color{black}
\ell_i(T):=
\left\lfloor
\rho_\alpha
\left(\frac{\Lambda_T}{\Delta_i^2}\right)^{1/(2\alpha)}
\right\rfloor.
}
\end{equation}
For sufficiently large \(T\), we have \(\ell_i(T)\ge s_i^\circ\). On \(\mathcal C_i\), after every arm has entered the concentration region and while \(k_i(t)<\ell_i(T)\), equations \eqref{eq:qcare-statewise-lower-bound} and \eqref{eq:empirical-gap-upper-bound} imply
$
Q_i(S(t))
\ge
2^{-(N-1)}\bar\Phi\!\left(\frac{3\Delta_i}{2}\ell_i(T)^\alpha\right).
$
By the definitions of \(\rho_\alpha\) and \(\ell_i(T)\), $
\frac{3\Delta_i}{2}\ell_i(T)^\alpha
\le
\frac{3}{2}\rho_\alpha^\alpha\sqrt{\Lambda_T}
=
\sqrt{\Lambda_T}. $
The standard Gaussian lower bound on the upper tail gives \(\bar\Phi(x)\ge c_0(1+x)^{-1}e^{-x^2/2}\) for \(x\ge0\). Therefore,
\begin{equation}
\label{eq:qcare-target-selection-lower-bound}
{\color{black}
Q_i(S(t))
\ge
q_T
:=
c_0 2^{-(N-1)}(eN)^{-1/2} \frac{T^{-1/2}}{1+\sqrt{\Lambda_T}}.
}
\end{equation}

Let \(H_i^{\mathrm{tar}}(T)\) be the additional number of periods after \(H_i^{\mathrm{ent}}\) required for arm \(i\) to reach \(\ell_i(T)\) pulls. Conditional on \(\mathcal C_i\), the current Gaussian shocks remain independent of the reward table and the past. Thus, \eqref{eq:qcare-target-selection-lower-bound} and the geometric waiting-time argument give
\[
{\color{black}
\mathbb E[H_i^{\mathrm{tar}}(T)\mid\mathcal C_i]
\le
\frac{\ell_i(T)}{q_T}
\le
\tfrac{1}{c_0 2^{-(N-1)}(eN)^{-1/2}}\ell_i(T)(1+\sqrt{\Lambda_T})\sqrt T
=o(T).
}
\]
Combining this bound with the finite expected time required to enter the concentration region and applying Markov's inequality gives, for all sufficiently large \(T\),
$
\mathbb P\{k_i(T+1)\ge \ell_i(T)\mid\mathcal C_i\}
\ge
\frac12.
$ Since \(\mathbb P(\mathcal C_i)\ge1/2\), we obtain
\[
{\color{black}
\mathbb E[k_i(T+1)]
\ge
\frac{\ell_i(T)}{4}.
}
\]
For all sufficiently large \(T\), the quantity inside the floor in \eqref{eq:target-pull-count} is at least two, and hence
\[
{\color{black}
\ell_i(T)
\ge
\frac{\rho_\alpha}{2}
\left(\frac{\log(eNT)}{\Delta_i^2}\right)^{1/(2\alpha)}.
}
\]
Consequently,
$
\mathbb E[k_i(T+1)]
\ge
\frac{\rho_\alpha}{8}
\left(\frac{\log(eNT)}{\Delta_i^2}\right)^{1/(2\alpha)}.
$
Taking the maximum of the finitely many arm-specific thresholds gives a common finite threshold
$
T_0=T_0(\alpha,N,\Theta)<\infty.
$
Finally, summing over the suboptimal arms yields
\begin{align*}
\mathcal{R}^{\Theta}(T)
&=
\sum_{\substack{i\ne i^\ast,\ \Delta_i>0}}
\Delta_i\,\mathbb E[k_i(T+1)] \ge
\frac{\rho_\alpha}{8}
\bigl(\log(eNT)\bigr)^{1/(2\alpha)}
\sum_{\substack{i\ne i^\ast,\ \Delta_i>0}}
\Delta_i^{1-1/\alpha}.
\end{align*}
The conclusion follows by taking \(c_{\alpha,N}:=\rho_\alpha/8\).
\end{proofref}

\subsection{The General Framework for Analyzing the Regret Lower Bounds}
\label{appendix:formal-general-framework}

This appendix develops the general Markovian-policy framework summarized in Section~\ref{sec:key_steps_generalizations}. We first define the structural properties used in the analysis. We then establish three sets of results, which lead to (i) complete learning and the \(o(T)\) regret for every fixed instance; (ii) the worst-case lower bound \(\Omega(T^{1-\alpha})\) when $ \alpha \in (0,0.5) $; and (iii) the lower bound \(\Omega(T^{1-o(1)})\) on every interior Bernoulli instance, and hence in the worst case, when $ \alpha > 0.5 $, respectively.  In the latter part of this section, we verify these properties for QCARE with general shock variances.

Recall that \(k_i(t)\) and \(\hat{\mu}_i(t)\) denote the pull count and empirical reward of arm \(i\) up to period \(t-1\). We consider MAB policies that depend only on the state
\[
{\color{black}
\bmS(t):=(k_1(t),\ldots,k_N(t),\hat{\mu}_1(t),\ldots,\hat{\mu}_N(t))\in\calS:=\mathbb{Z}_+^N\times[0,1]^N.}
\]
A typical state is denoted by \(S=(\kappa,u)\), where \(\kappa=(\kappa_1,\ldots,\kappa_N)\in\mathbb{Z}_+^N\) and \(u=(u_1,\ldots,u_N)\in[0,1]^N\). An admissible Markovian policy can be represented by a probability function \(Q=(Q_1,\ldots,Q_N):\calS\to\Delta_N\), where
\[
{\color{black}
Q_i(S)=\Pr{a(t)=i\mid \bmS(t)=S}}
\]
is the probability of pulling arm \(i\) at state \(S\). With a slight abuse of notation, we also write \(Q_i(S)=Q_i(S^i;S^{-i})\), where \(S^i=(\kappa_i,u_i)\) is the state of arm \(i\) and \(S^{-i}\) collects the states of the remaining arms. The family of MAB policies we consider is equivalent to the \textit{Sequentially Randomized Markov Experiments} studied by \cite{kuang2024weak}. It is easy to verify that QCARE belongs to this family, and we refer the reader to \cite{kuang2024weak} for many other members of this family.

We now define the relevant properties of the probability function.

\begin{definition}\label{def:anonymous}
	The function $ Q  $ is \textit{anonymous} if it is invariant to relabeling. That is, for every permutation $ \pi: [N] \to [N] $, $ i \in [N] $, and state $ S = (\kappa_1, \ldots, \kappa_N, u_{1}, \ldots, u_{N}) \in \calS $, $$ Q_i(\kappa_{\pi(1)}, \ldots, \kappa_{\pi(N)}, u_{\pi(1)}, \ldots, u_{\pi(N)}) = Q_{\pi(i)} (\kappa_1, \ldots, \kappa_N, u_{1}, \ldots, u_{N}). $$
\end{definition}

This property means that the probability function needs to be symmetric with respect to relabeling of the arms. In other words, it cannot ``cheat'' by depending on the arm identity information. This is a minimal requirement for $ Q $, and we assume it throughout the section.

\begin{definition}\label{def:mixed}
	The function $ Q  $ is \textit{recurrent} if for every arm $ 1 $'s state $ S^1 = (\kappa_1, u_1) \in \mathbb{Z}_+ \times [0,1] $, there exists $ \varepsilon > 0 $ (possibly dependent on $S^1$) so that 
	\begin{align*}
	Q_1 (S^1; S^{-1}) \geq \varepsilon \quad \text{ for all } S^{-1}.
	\end{align*}
\end{definition}

This property means that after an arm is pulled, there is always a non-diminishing probability that it is re-pulled (regardless of the history of other arms). The term ``recurrent'' is formally justified by the result below.

\begin{remark}
	{\sf Note that Definition \ref{def:mixed} only nominally describes a condition for $ Q_1 $ (i.e., arm 1). Since the probability function $ Q $ is anonymous, similar conditions apply to other arms, too. To be more precise, for every arm $ i $'s state $ S^i = (\kappa_i, u_i) \in \mathbb{Z}_+ \times [0,1] $, there exists $ \varepsilon > 0 $ (possibly dependent on $S^i$) so that $ Q_i (S^i; S^{-i}) \geq \varepsilon $ for all $ S^{-i}.$ The definitions below are stated with this anonymity property understood.}
\end{remark}

\begin{lemma}\label{lemma:k1_inf}
	If the probability function $Q$ is recurrent, then with probability one, every arm is pulled infinitely many times, i.e., $ k_i (t) \to \infty $ as $ t \to \infty $ for every $ i \in [N] $.
\end{lemma}

\begin{proofref}{Lemma \ref{lemma:k1_inf}}
			Fix an arm $ i \in [N] $. Let $ t_0 := 0 $ and for $ j \in \mathbb{Z}_+ $, let 
			$
			t_j := \min \{t : k_i(t) \geq j\}
			$
			be the first period after the $j^{th}$ pull of the arm $ i $. Define the process $\{W_j\}_{j\geq 0}$, where $W_0 :=0$ and $W_j:=t_j-t_{j-1}$ for every $j\geq 1$. This process represents the inter-pulling times between the $ (j-1) $th and $ j $th pulls of arm $ i $.
			It suffices to show that  $W_j < \infty$ almost surely for every $j\geq 1$.

			Fix an arbitrary $ j\geq 1 $.  By induction, assume that $ W_\ell < +\infty $ almost surely for all $ \ell \in \{0, \ldots, j-1\} $. In other words, $ t_{j-1} < +\infty  $. We wish to show that $ t_{j} < \infty $ almost surely, too.  Pick arbitrary $ u \in [0,1] $ and it suffices to show that  $ t_{j} < \infty $ almost surely conditional on   $ \{\hat{\mu}_i (t_{j-1}) = u\} $.
			
			Between the $ (j-1) $th and $ j $th pull of arm $ i $, the state value for this arm is kept at $ (j-1,u) $.
			Because the $ Q $ function is recurrent, there exists $ \varepsilon >0 $ such that
			$$  \Pr{a(t)= i \mid S^i(t) = (j-1,u)}  \,\geq\, \inf_{S^{-i}} Q_i\big((j-1,u);S^{-i}\big) \,\geq\, \varepsilon > 0.
			 $$
			 For all $ m \in \{0, 1, \ldots\} $, let $ E_m := \big\{a(t_{j-1} + \ell)  \neq i \text{ for all } \ell \in \{0,\ldots,m\} \big\} $ be the event that arm $ i $ keeps being skipped starting from period $ t_{j-1} $ to period $ t_{j-1} + m $. Therefore, $ \Pr{E_m \mid \hat{\mu}_i (t_{j-1}) = u} \leq (1-\varepsilon)^{m+1} $. As a consequence, $ \sum_{m=0}^{\infty} \Pr{E_m \mid \hat{\mu}_i (t_{j-1}) = u}  < \infty $, and due to the Borel-Contelli lemma, the events $ \{E_m\} $ only happen finitely often. In other words, arm $ i $ is eventually pulled within finitely many periods.
\end{proofref}
	
  \vspace{0.3 cm}
\begin{definition}\label{def:consistent}
	The function $ Q  $ is \textit{leading-arm convergent} if for every sequence $(\kappa, u)$ satisfying $ \kappa_i = \omega(1)$ for every arm $ i \in [N] $ and  $ u_1 - u_i = \Omega(1) $ for every arm $ i \geq 2 $, it holds that $$ 1 - Q_1((\kappa, u))  = o(1). $$
\end{definition}

Roughly speaking, a policy with a leading-arm convergent probability function $Q$ concentrates its pulling probability on the leading arm when there is sufficiently strong evidence that this arm is optimal, i.e., every arm has been pulled infinitely many times and the advantage of the leading arm is still non-diminishing. The two properties in Definitions \ref{def:mixed} and \ref{def:consistent} are fairly mild and intuitive. 
That said, they are sufficient to guarantee a reasonable performance of the policy, which is stated in the result below.

\begin{theorem}\label{thm:main2}
	If the probability function $Q$ is recurrent and leading-arm convergent, then for every problem instance $ \Theta $, $ \Pr{a(t) = i^\ast} \to 1 $ as $ t \to \infty $. As a result, the regret satisfies $\mathcal{R}^\Theta(T)= o(T)$ and the policy is long-run-average optimal.
\end{theorem}

\begin{proofref}{Theorem \ref{thm:main2}}
Following Lemma~\ref{lemma:k1_inf} and the recurrence property of \(Q\), we have \(k_i(t)\to\infty\) almost surely for every \(i\in[N]\). Hence, by the law of large numbers, \(\hat{\mu}_i(t)\to\mu_i\) almost surely for all arms \(i\). Without loss of generality, suppose that the unique optimal arm is \(i^*=1\). Then, almost surely, arm 1 is eventually the leading arm, and its empirical mean gap from every suboptimal arm is eventually bounded away from zero.

Since \(Q\) is leading-arm convergent, it follows that \(1-Q_1(\mathbf{S}(t))\to 0\) almost surely. Moreover, \(0\leq 1-Q_1(\mathbf{S}(t))\leq 1\). Therefore, by the bounded convergence theorem,
\(\Pr{a(t)\neq 1}=\mathbb{E}\!\left[1-Q_1(\mathbf{S}(t))\right]\to 0.\)
As a result, 
\(\sum_{t=1}^{T}\Pr{a(t)\neq 1}=o(T).\) 
Therefore, the expected number of pulls of suboptimal arms satisfies
\[ \color{black} { \sum_{i=2}^{N}\mathbb{E}[k_i(T+1)] = \sum_{t=1}^{T}\Pr{a(t)\neq 1} = o(T). }\]

Recall that \(\Delta_i=\mu^*-\mu_i\). Then
\[ \color{black} {\mathcal{R}^{\Theta}(T) = \mathbb{E}\left[\sum_{t=1}^{T}\left(\mu^*-\mu_{a(t)}\right)\right]
\leq (\max_i \Delta_i)\sum_{i=2}^{N}\mathbb{E}[k_i(T+1)] = o(T). }\]
Thus, the average regret vanishes and
\( \frac{1}{T}\mathbb{E}\left[\sum_{t=1}^{T}r_{a(t)}(t)\right]\to \mu^*. \)
Since no feasible MAB policy can have long-run average expected reward exceeding \(\mu^*\), the policy induced by \(Q\) is long-run-average optimal. 
\end{proofref}	

\noindent \textbf{Conditions for ``over'' and ``under'' exploration.} Let us now proceed to two key novel measures we use to quantify ``over'' and ``under'' exploration. Let us start with the one corresponding to over-exploration.

\vspace{0.3 cm}
\begin{definition}\label{def:exploratory}
	Let $ \alpha \in (0,0.5) $ be fixed. The probability function  $ Q  $ is \textit{$\alpha$-exploratory} if there exist constants $ \varepsilon > 0 $, $\delta \in (0,1) $ and a sequence of numbers $ \{\Delta_1, \Delta_2, \ldots\} $, where $ \Delta_T = \Omega(T^{-\alpha}) $, such that for all sufficiently large $ T $, $ u_1 \leq \Delta_T $, and $ \kappa $ such that $ \sum_{i>1} \kappa_i < \delta T $, it holds that 
	\begin{align*}
	1 -   Q_1\left(\kappa_1,\kappa_2,\ldots,\kappa_N,u_1, 0,0,\ldots,0\right) \geq \varepsilon.
	\end{align*}
\end{definition}

Roughly speaking, a policy with an $ \alpha $-exploratory function $ Q $ allows situations in which, even though all arms have been pulled $ \Omega(T) $ times and their empirical reward differences are on the order of $ \Omega(T^{-\alpha}) $, the arm-pulling probabilities do not concentrate on the leading arm. In other words, an $ \alpha $-exploratory function $ Q $ with a small $ \alpha $ is a sign of heavy exploration.  Intuitively, if the exploration is too heavy, then the policy may fail to concentrate on any particular arm despite overwhelming evidence that a certain arm is the best. In the result below, we formally characterize the consequence of such ``over-exploration.'' 

\begin{theorem}\label{thm:over explore}
	Let $ \alpha \in (0,0.5) $ and $ N> 1 $ be fixed. If the probability function $ Q $ is $\alpha$-exploratory, then the regret satisfies $\mathcal{R}(T)= \Omega(T^{1-\alpha}). $
\end{theorem}

	\begin{proofref}{Theorem \ref{thm:over explore}}
Fix $ N > 1 $ and $ \alpha \in (0,0.5) $. Suppose the $ Q $ function is $ \alpha $-exploratory. Pick constants $ \varepsilon > 0 $, $ \delta \in (0,1) $, and a sequence $  \{\Delta_T\} $ satisfying $\Delta_T = \Omega(T^{-\alpha}) $ as required in Definition \ref{def:exploratory}. Therefore, for all sufficiently large $ T \geq T_0 $, $ u_1 \leq \Delta_T $, and $ \kappa $ such that $ \sum_{i>1} \kappa_i < \delta T $, it holds that 
		$
		1 -   Q_1\left(\kappa_1,\kappa_2,\ldots,\kappa_N,u_1, 0,0,\ldots,0\right) \geq \varepsilon.
		$
		
		We construct a problem instance for every $ T \geq T_0 $, which we denote by $ \Theta_T $. Under this instance, the reward distribution for each arm is \textit{deterministic}. The rewards of the arms are $\mu_1= \min\{1, \Delta_T\} = \Omega(T^{-\alpha})$ and  $\mu_2= \mu_3 = \cdots = \mu_N=0$. Let $ \delta_2 := (1+\delta)/2 $ so that $ \delta_2 \in (\delta, 1) $. For every $  t \geq \delta_2 T $, we introduce the event $$ A_t := \Big\{\sum_{i > 1}k_i(t)\leq \delta T\Big\}. $$
		Let us discuss the consequences of this event for future reference. First, since $ \sum_{i > 1}k_i(t) <  \delta T  $, it holds that $ k_1(t) = (t-1) - [ \sum_{i > 1}k_i(t)] \geq \delta_2 T - 1 - \delta T = \frac{1-\delta}{2} T -1  > 0 $ for $ T > 2/(1-\delta) $. In other words, arm 1 has been pulled at least once before period $ t $. Because the rewards are deterministic, we must have $ \hat{\mu}_1(t) = \mu_1 \times k_1(t)/(k_1(t)+1) \leq \Delta_T $ and $ \hat{\mu}_2(t) = \cdots = \hat{\mu}_N(t) = 0 $. (Here we use the convention that if an arm is not pulled, the reward state is valued at zero.) Since the probability function $Q$ is $\alpha$-exploratory, it holds that for all $ t \geq \delta_2 T $,
		\begin{align*}
		\Pr{a(t) \neq 1 |A_t} = \ &1-Q_1(  k_1(t), \ldots, k_N(t),\hat{\mu}_1(t), \ldots,  \hat{\mu}_N(t) | k \in A_t) \\
		= \ &1-Q_1( k_1(t), \ldots, k_N(t),\hat{\mu}_1(t), 0, \ldots,  0 | k \in A_t)  \geq \varepsilon.
		\end{align*}
		
		Let $ c_0 := 1/2 $ be an absolute constant. We break our proof into two cases depending on the values of $ \{\Pr{A_t}\} $. 
		
		\noindent \underline{Case 1:  $ \Pr{A_t} \geq c_0 $ for all $ t  > \delta_2 T $.} In this case,
		\begin{align*}
		\mathcal{R}(T) \,\geq\,  \mathcal{R}^{\Theta_T}(T)
		=&\ \mu_1 \sum_{t=1}^{T} \Pr{a(t) \neq 1}   
		\geq  \mu_1 \sum_{t= \lceil  \delta_2 T \rceil + 1  }^{T} \Pr
		{a(t) \neq 1 | A_t} \,\Pr{A_t} \\
		= &\Omega(T^{-\alpha}) \cdot (1-\delta_2) T \cdot \varepsilon c_0 = \Omega(T^{1-\alpha}).
		\end{align*}
		
		\noindent \underline{Case 2: There exists $ \tilde{t} \in (\delta_2 T, T] $ such that $ \Pr{A_{\tilde{t}}} < c_0 $.}  In this case, note that  $\sum_{i > 1}k_i(t)$ is non-decreasing in $t$. As a result, the events $ \{A_t\}_t $ are non-increasing. (That is, $ A_t $ is false implies $ A_\ell $ is false for all $ \ell \geq t $.) Therefore, $ \Pr{A_{T+1}} \leq \Pr{A_{\tilde{t}}} < c_0 $. In other words, $ \Pr{\overline{A_{T+1}}} \geq 1 - c_0 $. Let us evaluate the cumulative regret below. Due to Markov's inequality,
		\begin{align*}
		\mathcal{R}(T) \,\geq\, \mathcal{R}^{\Theta_T}(T) =&\mu_1 \sum_{i>1}  \mathbb{E}\left[k_{i}(T+1)\right] \\ \geq &\mu_1 \left(\delta T\right) \Pr{\sum_{i > 1}k_i(T+1) >  \delta T} \\
		= &\mu_1 \left(\delta T\right) \Pr{\overline{A_{T+1}}} \\
		\geq &\mu_1 \left(\delta T\right) (1- c_0) =\Omega(T^{-\alpha}) \left(\delta T\right) (1- c_0) =  \Omega(T^{1-\alpha}).
		\end{align*}
		That finishes the proof. 
		\end{proofref}
	
To elaborate on the consequence of under-exploration, we need a technical condition below to rule out some extreme cases.
	
\vspace{0.3 cm}
\begin{definition}\label{def:unradical}
	The function  $ Q  $ is \textit{unradical} if for every $ \kappa \in \mathbb{Z}_+^N $ and $ u \in [0,1]^N $ such that $ u_1 \leq \min \{u_2, \ldots, u_N\} $, $ Q_1(\kappa,u) \leq \frac{1}{2}. $
\end{definition}

The property means that the weakest-performing arm (i.e., the arm with the lowest empirical mean) is pulled with a probability no larger than $ 1/2 $. This is fairly intuitive since one would expect that the better-performing arms would be pulled with higher probabilities. So we call this property ``unradical.'' The next property provides a quantitative measure of how a policy could potentially ``under explore.''

\begin{definition}\label{def:irreversible}
    Let $ \alpha\in (0.5, +\infty) $ and $ N> 1 $ be fixed. The function  $ Q  $ is \textit{$\alpha$-irreversible} if for \textit{any} of the following quantities: (i) arbitrarily small constant $ \delta >0 $, (ii) arm $ j >1 $, (iii) $ u \in [0,1]^N $ satisfying $ u_j - u_1 = \Omega(1) $, and  (iv) $ \kappa \in \mathbb{Z}_+^N  $ satisfying  $ \kappa_1, \kappa_j \geq (\log T)^{\frac{1}{2\alpha-\delta}} $ for sufficiently large $ T $, it holds that
	\begin{align*}
	Q_1(\kappa, u) \leq  \frac{(\log T)^{\frac{1}{2\alpha -\delta}}}{T}.
	\end{align*}
\end{definition}

Roughly speaking, a policy with an $\alpha$-irreversible function $ Q $ allows situations where two arms have both been pulled only $ (\log T)^{\frac{1}{2\alpha - \delta}} $ times, while the worse-performing arm's pulling probability has already diminished to the order of $ O \left(\frac{(\log T)^{\frac{1}{2\alpha - \delta}}}{T}\right). $ After a preliminary calculation, this loosely translates to the probability of pulling the non-leading arms reducing super-exponentially fast with the pull count. Hence, that corresponds to a sign of aggressive concentration (or lack of exploration). Intuitively, if the concentration is too aggressive, then the policy may be mistakenly stuck in the inferior arm in a ``bad'' event where the inferior arm performs better empirically during the initial periods. This intuition is also where the name ``irreversibility'' comes from. In the result below, we formally characterize the performance consequences of ``under-exploration.''

\begin{theorem}\label{thm:under explore}
Fix $N>1$ and $\alpha>0.5$. Consider any Bernoulli instance $\Theta$ with $0<\mu_i<\mu_1<1$ for every $i>1$. For every policy with an unradical
and $\alpha$-irreversible probability function $Q$, it pulls
suboptimal arms $T-o(T)$ times with probability $\Omega(T^{-o(1)})$.
Consequently, $\mathcal{R}^{\Theta}(T)=\Omega(T^{1-o(1)})$, and the same
lower bound holds for $\mathcal{R}(T)$.
\end{theorem}

\begin{proofref}{Theorem~\ref{thm:under explore}}
Pick $\delta>0$ sufficiently small that $2\alpha-2\delta>1$. Define
\[ c:=\left\lceil2(\log T)^{\frac{1}{2\alpha-2\delta}}\right\rceil, \quad \text{ and } \quad
t_c:=\inf\{t\geq1:k_1(t)\geq c\}.
\]
In short, $ t_c $ is the first period after the $c$th pull of arm 1, and by convention, we let $t_c:=+\infty$ if that never occurs. Let $E_1$ be the event that either $t_c=+\infty$, or $t_c<+\infty$ and $\sum_{i=2}^{N}k_i(t_c)>0.9c$. Let $E_2$ be the event that the first $c$ rewards of arm 1 are all zero. Let $E_3$ be the event that, for every $i>1$ and every integer $n\geq0.9c/(N-1)$, the empirical reward based on the first $n$ rewards and denominator $n+1$ is at least $\mu_i/2$. Set
$E:=E_1\cap E_2\cap E_3$. We break the rest of the proof into three steps.

\noindent\underline{Step 1: The misleading event has probability
$\Omega(T^{-o(1)})$.} Let $\lambda_1:=-\log(1-\mu_1)>0$. Since
$c/\log T\to0$,
\[
\Pr{E_2}=(1-\mu_1)^c=\exp(-\lambda_1c)=T^{-o(1)}.
\]
On event $E_2$, arm 1 is empirically weakest until its $c$th pull. Unradicality and
a sequential coupling show that the number of suboptimal-arm pulls before
that pull stochastically dominates the number $X_c$ of failures before $c$
successes in independent Bernoulli trials with success probability $1/2$.
Thus, even after conditioning on $E_3$,
\[
\Pr{\overline{E_1}\mid E_2,E_3}
\leq\Pr{|X_c-c|\geq0.1c}
\leq\frac{200}{c},
\]
because $\mathbb{E}[X_c]=c$ and $\operatorname{Var}(X_c)=2c$. This coupling
is pointwise in the suboptimal-arm reward sequences.

Put $\underline\mu:=\min_{i>1}\mu_i$ and
$n_0:=\lceil0.9c/(N-1)\rceil$. For all sufficiently large $T$, failure of
$E_3$ at arm $i$ and sample size $n$ implies that the ordinary sample average
of the first $n$ rewards is below $3\mu_i/4$. Hence Hoeffding's inequality and
a union bound yield
\[
\Pr{\overline{E_3}}
\leq
\frac{N-1}{1-\exp(-\underline\mu^2/8)}
\exp\left(-\frac{\underline\mu^2n_0}{8}\right)
\leq C_3e^{-c_3c},
\]
where $C_3:=(N-1)/(1-\exp(-\underline\mu^2/8))$ and
$c_3:=0.9\underline\mu^2/(8(N-1))$. Since $E_2$ and $E_3$ concern
independent arm-reward sequences,
\[
\Pr{E}\geq(1-\mu_1)^c
\left(1-\frac{200}{c}-C_3e^{-c_3c}\right)=\Omega(T^{-o(1)}).
\]

\noindent\underline{Step 2: Conditional on $E$, the policy remains committed to
suboptimal arms for nearly the entire horizon.}

Let $m_0:=\min\{0.1,\frac18\min_{i>1}\mu_i\}>0$ and
$M:=\lfloor m_0c\rfloor$. If $t_c\geq T+1$, then
$\sum_{i=2}^{N}k_i(T+1)\geq T-c$. Otherwise, choose $j>1$ with
$k_j(t_c)=\max_{i>1}k_i(t_c)$. On $E_1$,
$k_j(t_c)>0.9c/(N-1)$. Until arm 1 receives another $M$ pulls, it holds that
\[
\hat\mu_1(t)\leq\frac{M}{c+1}\leq\frac18\min_{i>1}\mu_i
\quad \text{ and } \quad 
\hat\mu_j(t)\geq\frac12\min_{i>1}\mu_i.
\]
Writing $a:=1/(2\alpha-2\delta)$ and $b:=1/(2\alpha-\delta)$, $a>b$
implies that both $c/(\log T)^b$ and
$(0.9c/(N-1))/(\log T)^b$ tend to infinity. The conditions of
$\alpha$-irreversibility therefore hold, so during this period
\[
Q_1(\bmS(t))\leq\frac{(\log T)^{\frac{1}{2\alpha-\delta}}}{T}=:g.
\]
The number of suboptimal-arm pulls before another $M$ pulls of arm 1
stochastically dominates a negative-binomial variable $X$ with $M$ successes
and success probability $g$. Let
$
\delta':=\frac{1}{2\alpha-2\delta}-\frac{1}{2\alpha-\delta}>0.
$
For all sufficiently large $T$,
\[
m_T:=\mathbb{E}[X]=\frac{M(1-g)}{g}
\geq\frac{m_0}{2}T(\log T)^{\delta'}\geq2T.
\]
Since $\operatorname{Var}(X)=M(1-g)/g^2$, Chebyshev's inequality gives
\[
\Pr{X<T}\leq\frac{\operatorname{Var}(X)}{(m_T-T)^2}
\leq\frac{4}{M(1-g)}\leq\frac8M.
\]
Consequently, conditional on $E$, with probability $1-O(1/c)$,
$
\sum_{i=2}^{N}k_i(T+1)\geq T-c-M=T-o(T).
$

\noindent\underline{Step 3: Regret lower bound.} Combining the two cases for $t_c$ gives
\[
\Pr{\sum_{i=2}^{N}k_i(T+1)\geq T-c-M}
\geq\Pr{E}\left(1-O\left(\frac1c\right)\right)
=\Omega(T^{-o(1)}).
\]
Finally,
\[
\mathcal{R}^{\Theta}(T)
\geq\left(\min_{i>1}\Delta_i\right)(T-c-M)
\Pr{\sum_{i=2}^{N}k_i(T+1)\geq T-c-M}
=\Omega(T^{1-o(1)}).
\]
Since $\mathcal{R}(T)\geq\mathcal{R}^{\Theta}(T)$, the claimed
instance-independent lower bound follows. 
\end{proofref}
	
The theorem reveals that without sufficient exploration, the aforementioned ``bad'' event happens with probability $ \Omega(T^{-o(1)}) $ and the regret conditional on that bad event is on the order of $ \Omega(T) $. That makes the overall regret nearly linear in $ T $. Next, we reveal the connection between the QCARE policy family and our general framework by showing that QCARE satisfies all properties defined above.

\begin{proposition}
\label{prop: QCARE satisfy Q}
For every fixed \(N>1\), \(\sigma>0\), and \(\alpha>0\), QCARE with general variance $ \sigma^2 $ is anonymous, recurrent, leading-arm convergent, and unradical. When \(\alpha\in(0,0.5)\), it is \(\alpha\)-exploratory. When \(\alpha\in(0.5,+\infty)\), it is \(\alpha\)-irreversible.
\end{proposition}

\textbf{Proof of Proposition~\ref{prop: QCARE satisfy Q}.} Let \(\epsilon_1,\ldots,\epsilon_N\) be i.i.d. standard normal random variables and define the score for QCARE with general variance
 \( {\color{black} \theta_i^\sigma(S)=u_i+\sigma(\kappa_i+1)^{-\alpha}\epsilon_i.} \)  Let \(Q_i^\sigma(S)\) denote the probability that arm \(i\) has the largest score, and let \(Q^\sigma=(Q_1^\sigma,\ldots,Q_N^\sigma)\) denote the resulting Markovian policy. Define
\( {\color{black} \beta_i(\kappa):=(\kappa_i+1)^{-\alpha}.} \)
At a fixed state \(S=(\kappa,u)\), the score difference between arms \(j\) and \(i\) is Gaussian:
\begin{equation}
{\color{black}
\theta_j^\sigma(S)-\theta_i^\sigma(S)
\sim
N\left(u_j-u_i,\sigma^2\left[\beta_i(\kappa)^2+\beta_j(\kappa)^2\right]\right).}
\label{eq:pairwise_sigma_comparison}
\end{equation}
Thus, whenever \(u_i>u_j\),
\begin{equation}
{\color{black}
\Pr{\theta_j^\sigma(S)\ge \theta_i^\sigma(S)}
=
\Phi\left(
-
\frac{u_i-u_j}
{\sigma\sqrt{(\kappa_i+1)^{-2\alpha}+(\kappa_j+1)^{-2\alpha}}}
\right).}
\label{eq:pairwise_sigma_probability}
\end{equation}
Next, we verify the properties in Proposition~\ref{prop: QCARE satisfy Q} one by one.

The function \(Q^\sigma\) is anonymous because the score rule has the same functional form for every arm and the Gaussian shocks are i.i.d. \(N(0,\sigma^2)\). Relabeling the arms simply relabels the vector of scores, and therefore relabels the choice probabilities in the same way.

\begin{lemma}[QCARE is recurrent]
\label{lem:recurrence_sigma}
Fix \(\sigma>0\) and \(\alpha>0\). For each arm \(i\), there exist constants \(c_{\sigma,N},C_\sigma>0\) such that for every state \((\kappa,u)\in\mathbb Z_+^N\times[0,1]^N\),
\begin{equation}
{\color{black}
Q_i^\sigma(\kappa,u)
\ge
c_{\sigma,N}\exp\left\{-C_\sigma(\kappa_i+1)^{2\alpha}\right\}.}
\label{eq:recurrence_sigma_bound}
\end{equation}
\end{lemma}

\begin{proofref}{Lemma~\ref{lem:recurrence_sigma}}
Write the score as \(u_i+\sigma(\kappa_i+1)^{-\alpha}\epsilon_i\). Consider the event
\[ {\color{black}
E_i=
\left\{
\epsilon_i\ge \frac{(\kappa_i+1)^\alpha}{\sigma},
\qquad
\epsilon_j\le 0\quad \text{for all }j\ne i
\right\}.}
\]
On \(E_i\),
\(\theta_i^\sigma(S) =u_i+\sigma(\kappa_i+1)^{-\alpha}\epsilon_i \ge u_i+1\ge 1,\)
while for every \(j\ne i\),
\( {\color{black} \theta_j^\sigma(S) = u_j+\sigma(\kappa_j+1)^{-\alpha}\epsilon_j \le u_j\le 1.} \)
Hence \(E_i\subseteq\{a=i\}\) up to a null event, and therefore \(Q_i^\sigma(\kappa,u)\ge \Pr{E_i}\). By independence,
\[
{\color{black}
\Pr{E_i}
=
\bar\Phi\left(\frac{(\kappa_i+1)^\alpha}{\sigma}\right)2^{-(N-1)},}
\]
where \(\bar\Phi(x)=1-\Phi(x)\). Let \(x=(\kappa_i+1)^\alpha/\sigma\). The Gaussian lower-tail bound gives
\( {\color{black} \bar\Phi(x)\ge c(1+x)^{-1}e^{-x^2/2}.} \)
For any fixed \(\eta>0\), there exists \(c_\eta>0\) such that
\( {\color{black} (1+x)^{-1}\ge c_\eta e^{-\eta x^2}, x\ge0.} \)
Therefore,
\[
{\color{black}
\bar\Phi(x)
\ge
cc_\eta e^{-(1/2+\eta)x^2}
=
cc_\eta
\exp\left\{
-\frac{(1/2+\eta)(\kappa_i+1)^{2\alpha}}{\sigma^2}
\right\}.}
\]
After multiplying by \(2^{-(N-1)}\) and absorbing constants into \(c_{\sigma,N}\) and \(C_\sigma\), we obtain \eqref{eq:recurrence_sigma_bound}.
\end{proofref}

\begin{lemma}[QCARE is leading-arm convergent]
\label{lem:leading_sigma}
Fix \(\sigma>0\) and \(\alpha>0\). Consider any sequence of states \((\kappa,u)=(\kappa(T),u(T))\in\mathbb Z_+^N\times[0,1]^N\) such that \(\kappa_i(T)=\omega(1)\) for every arm \(i\in[N]\) and \(u_1(T)-u_i(T)=\Omega(1)\) for every arm \(i\ge2\). Then
$
1-Q_1^\sigma(\kappa(T),u(T))=o(1).
$
\end{lemma}

\begin{proofref}{Lemma~\ref{lem:leading_sigma}}
Recall that the score rule for QCARE with general variance is given by
\( {\color{black} \theta_i^\sigma=u_i+\sigma\beta_i(\kappa(T))\epsilon_i,} \)
where \(\epsilon_1,\ldots,\epsilon_N\) are independent standard normal random variables.
For every fixed $ i $, let $p_i :=\Pr{\theta_i^\sigma\ge \theta_1^\sigma}$. By Equations~\eqref{eq:pairwise_sigma_comparison} and \eqref{eq:pairwise_sigma_probability}, for each \(i\ge2\), we have
\begin{equation}
{\color{black}
p_i =\Pr{\theta_i^\sigma\ge \theta_1^\sigma}
=
\Phi\left(
-
\tfrac{u_1-u_i}
{\sigma\sqrt{(\kappa_1+1)^{-2\alpha}+(\kappa_i+1)^{-2\alpha}}}
\right).}
\label{eq: pi_property}
\end{equation}
Since \(u_1-u_i=\Omega(1)\) and \(\kappa_1,\kappa_i\to\infty\), the denominator is \(o(1)\) while the numerator is bounded below by a positive constant eventually. Thus the positive standardized gap diverges to \(+\infty\), the argument of \(\Phi\) goes to \(-\infty\), and thus \(p_i = o(1)\). The union bound thus implies that $ 1-Q_1^\sigma(\kappa,u) \le \sum_{i=2}^N
	p_i = o(1) $.
\end{proofref}

\begin{lemma}[QCARE is \(\alpha\)-exploratory for \(0<\alpha<0.5\)]
\label{lem:exploratory_sigma}
Let \(N>1\), fix \(\sigma>0\), and consider QCARE with general variance when \(0<\alpha<0.5\). Then \(Q^\sigma\) is \(\alpha\)-exploratory.
\end{lemma}

\begin{proofref}{Lemma~\ref{lem:exploratory_sigma}}
Fix \(\delta\in(0,1)\), for example \(\delta=1/2\), and set
\( {\color{black} \Delta_T=d_{\sigma,N}T^{-\alpha},} \)
where \(d_{\sigma,N}>0\) is a fixed constant, chosen below and allowed to depend on \(\alpha,\sigma,N\) and the fixed \(\delta\). Then \(\Delta_T=\Omega(T^{-\alpha})\). Consider states of the form
\( {\color{black}
(u_1,u_2,\ldots,u_N)=(u_1,0,\ldots,0), u_1\le \Delta_T,} \)
and suppose \(\sum_{i>1}\kappa_i<\delta T\). Then for some challenger \(j>1\),
$
{\color{black}
\kappa_j\le \frac{\delta T}{N-1}.}
$
Since the event \(\theta_j^\sigma\ge\theta_1^\sigma\) rules out selecting arm 1, up to the zero-probability tie event,
\( {\color{black} 1-Q_1^\sigma(\kappa,u_1,0,\ldots,0) \ge \Pr{\theta_j^\sigma\ge\theta_1^\sigma}.} \)
Based on Equation~\eqref{eq: pi_property} and \(u_j=0\),
$
{\color{black}
\Pr{\theta_j^\sigma\ge \theta_1^\sigma}
=
\Phi\left(
-\frac{u_1}
{\sigma\sqrt{(\kappa_1+1)^{-2\alpha}+(\kappa_j+1)^{-2\alpha}}}
\right).}
$
The denominator is bounded below by \(\sigma(\kappa_j+1)^{-\alpha}\). Thus, for all sufficiently large \(T\),
\[
{\color{black}
\frac{u_1}
{\sigma\sqrt{(\kappa_1+1)^{-2\alpha}+(\kappa_j+1)^{-2\alpha}}}
\le
\frac{d_{\sigma,N}T^{-\alpha}(\kappa_j+1)^\alpha}{\sigma}
\le
\frac{d_{\sigma,N}}{\sigma}
\left(\frac{2\delta}{N-1}\right)^\alpha .}
\]
Choose \(d_{\sigma,N}>0\) small enough that the last display is at most \(1\). Then
\[
1-Q_1^\sigma(\kappa,u_1,0,\ldots,0)\ge \Phi(-1)>0.
\]
The constants \(\delta\), \(d_{\sigma,N}\), and \(\Phi(-1)\) are fixed and independent of \(T\). Thus Definition~\ref{def:exploratory} is satisfied.
\end{proofref}

\begin{lemma}[QCARE is \(\alpha\)-irreversible for \(\alpha>0.5\)]
\label{lem:irreversible_sigma}
Let \(N>1\), fix \(\sigma>0\), and consider QCARE with general variance when \(\alpha>0.5\). Then \(Q^\sigma\) is \(\alpha\)-irreversible.
\end{lemma}

\begin{proofref}{Lemma~\ref{lem:irreversible_sigma}}
Fix an arbitrarily small \(\delta\in(0,2\alpha)\), an arm \(j>1\), and a sequence of states with \(u_j-u_1=\Omega(1)\) and
\( {\color{black}
\kappa_1,\kappa_j\ge (\log T)^{\frac{1}{2\alpha-\delta}}} \)
for sufficiently large \(T\), as in Definition~\ref{def:irreversible}. There exists \(c_0>0\) such that \(u_j-u_1\ge c_0\) for all sufficiently large \(T\). Since selecting arm 1 implies \(\theta_1^\sigma\ge\theta_j^\sigma\),
\( {\color{black} Q_1^\sigma(\kappa,u) \le \Pr{\theta_1^\sigma\ge\theta_j^\sigma}.} \)
Based on Equation~\eqref{eq: pi_property}, 
$
{\color{black}
\Pr{\theta_1^\sigma\ge \theta_j^\sigma}
= 1- p_j =
\Phi\left(
-
\frac{u_j-u_1}
{\sigma\sqrt{(\kappa_1+1)^{-2\alpha}+(\kappa_j+1)^{-2\alpha}}}
\right).}
$
Let \(L_T=(\log T)^{1/(2\alpha-\delta)}\). Since \(\kappa_1,\kappa_j\ge L_T\), the denominator is at most a constant multiple of \(\sigma L_T^{-\alpha}\). Hence the standardized gap is at least \(C_{\sigma,c_0}(\log T)^{\alpha/(2\alpha-\delta)}\) for some constant \(C_{\sigma,c_0}>0\). Since \(2\alpha/(2\alpha-\delta)>1\), the Gaussian upper-tail bound implies that, for all sufficiently large \(T\),
\[
{\color{black}
Q_1^\sigma(\kappa,u)
\le
\exp\left\{-C'_{\sigma,c_0}(\log T)^{\frac{2\alpha}{2\alpha-\delta}}\right\}
\le
\frac{(\log T)^{\frac{1}{2\alpha-\delta}}}{T}.}
\]
This is the \(\alpha\)-irreversible property.
\end{proofref}

\begin{lemma}[QCARE is unradical]
\label{lem: unradical}
Let \(N>1\), \(\sigma>0\), and \(\alpha>0\) be fixed. Under QCARE with general variance, the probability function \(Q^\sigma\) is unradical.
\end{lemma}

\begin{proofref}{Lemma~\ref{lem: unradical}} 
By anonymity, it suffices to verify the condition for an arbitrary weakest arm. We use arm \(N\) for notational clarity. Suppose
\( {\color{black} u_N\le \min\{u_1,\ldots,u_{N-1}\}.} \)
Pick any \(j<N\). Since selecting arm \(N\) implies \(\theta_N^\sigma\ge\theta_j^\sigma\),
\( {\color{black} Q_N^\sigma(\kappa,u) \le \Pr{\theta_N^\sigma\ge\theta_j^\sigma}.} \)
The difference \(\theta_N^\sigma-\theta_j^\sigma\) is Gaussian with mean \(u_N-u_j\le0\). Therefore, $
\Pr{\theta_N^\sigma\ge\theta_j^\sigma}\le \frac{1}{2},$
and hence \(Q_N^\sigma(\kappa,u)\le1/2\). By anonymity and relabeling, this verifies the unradicality condition in Definition~\ref{def:unradical}.
\end{proofref}

The preceding arguments establish anonymity, recurrence, leading-arm convergence, unradicality, the \(\alpha\)-exploratory property for \(0<\alpha<0.5\), and the \(\alpha\)-irreversible property for \(\alpha>0.5\). Combining the preceding arguments establishes all properties stated in Proposition~\ref{prop: QCARE satisfy Q}.

It remains to note that Proposition~\ref{prop: QCARE satisfy Q}, combined with Theorems~\ref{thm:main2}, \ref{thm:over explore}, and \ref{thm:under explore}, implies that QCARE with general variance achieves complete learning and satisfies the regret lower bounds $ \mathcal{R}(T) = \Omega(T^{1-\alpha}) $ and $ \mathcal{R}(T) = \Omega(T^{1-o(1)}) $ when $ 0 < \alpha < 0.5 $ and $ \alpha > 0.5 $, respectively. When specializing to $ \sigma^2 = 1 $, that implies Theorem~\ref{thm:UB_sublinear}, Proposition~\ref{prop:ub_smaller0.5_lower_worst}, and Theorem~\ref{thm:lb_greater0.5}, respectively.

\subsection{Reduced-Form Representation of QCARE}
\label{appendix:reduced-form-qcare}

Recall that \(Q_i^\sigma(S)\) defined in Appendix~\ref{appendix:formal-general-framework} denotes the probability of selecting arm \(i\), and \(\widehat Q(S)\) denotes the probability of selecting the current leading arm; and \(Q_{\mathrm{non\text{-}lead}}(S):=1-\widehat Q(S)\) denotes the probability of selecting a non-leading arm.
For the reduced-form calculations, write the score under this policy as
\( {\color{black} \theta_i^\sigma(S)=u_i+\sigma(\kappa_i+1)^{-\alpha}\epsilon_i.} \) 

We first consider the two-arm case. By the anonymity property in Proposition~\ref{prop: QCARE satisfy Q}, we may relabel the arms without changing the choice structure. Thus, without loss of generality, suppose \(u_1>u_2\), so arm 1 is the leading arm and arm 2 is the non-leading arm. 

\begin{theorem}[Two-Arm Reduced Form]
\label{thm:reduced-form-two-arm} 
For any state \(S=((\kappa_1,\kappa_2),(u_1,u_2))\) such that $ u_1 > u_2 $
\[
{\color{black}
Q_{\mathrm{non\text{-}lead}}(S)
=
1-\Phi\left(
\frac{u_1-u_2}
{\sigma\sqrt{(\kappa_1+1)^{-2\alpha}+(\kappa_2+1)^{-2\alpha}}}
\right).
}
\]
Moreover, along any sequence of states such that
\(\kappa \to\infty\) and \(u_1-u_2\) is bounded away from zero, the non-leading arm choice probability satisfies:
\[
{\color{black}
\begin{aligned}
\frac{1}{2}\log\log\left(\frac{1}{Q_{\mathrm{non\text{-}lead}}(S)}\right)
=\, & \alpha\log\left(\min\{\kappa_1+1,\kappa_2+1\}\right)
+ \log\left(\tfrac{u_1-u_2}{\sigma}\right)
-\frac{1}{2}\log 2 \\
&-\frac{1}{2}
\log\left[
\left(\tfrac{\min\{\kappa_1+1,\kappa_2+1\}}{\kappa_1+1}\right)^{2\alpha}
+
\left(\tfrac{\min\{\kappa_1+1,\kappa_2+1\}}{\kappa_2+1}\right)^{2\alpha}
\right]
+o(1) \\
=\,  & \alpha\log\left(\min\{\kappa_1+1,\kappa_2+1\}\right)+O(1).
\end{aligned}
}
\]
Therefore, in this log-log reduced form representation, \(\alpha\) corresponds to the slope with respect to the logarithm of arm-level pull count.
\end{theorem}

\begin{proofref}{Theorem~\ref{thm:reduced-form-two-arm}}
Based on Equation~\eqref{eq: pi_property}, 
\[
Q_{\mathrm{non\text{-}lead}}(S)
= p_2 = 
\Phi\left(
-
\frac{u_1-u_2}
{\sigma\sqrt{(\kappa_1+1)^{-2\alpha}+(\kappa_2+1)^{-2\alpha}}}
\right)
= 1-\Phi\left(
\frac{u_1-u_2}
{\sigma\sqrt{(\kappa_1+1)^{-2\alpha}+(\kappa_2+1)^{-2\alpha}}}
\right).
\] 
This proves the exact expression for \(Q_{\mathrm{non\text{-}lead}}(S)\). Let
$
z:=\frac{u_1-u_2}{\sigma\sqrt{\beta_1^2+\beta_2^2}}
$ for shorthand notation.
Along any sequence of states such that \(\min\{\kappa_1+1,\kappa_2+1\}\to\infty\) and \(u_1-u_2\) is bounded away from zero, we have \(z\to\infty\). Applying the standard tail approximation
\(1-\Phi(z)=\frac{1}{\sqrt{2\pi}z}e^{-z^2/2}(1+o(1))\) and taking the logarithm twice:
\[
{\color{black}
\log\left(\frac{1}{Q_{\mathrm{non\text{-}lead}}(S)}\right) = \frac{z^2}{2} + \log z + \frac{1}{2}\log(2\pi) + o(1).
}
\]
Since \(z^2\) dominates \(\log z\), we have
\(\log\log(1/Q_{\mathrm{non\text{-}lead}}(S))=\log(z^2/2)+o(1)=2\log z-\log 2+o(1)\). Consequently,
$$
{\color{black}
\frac12\log\log(1/Q_{\mathrm{non\text{-}lead}}(S))=\log z-\frac12\log 2+o(1).
}
$$
Substituting \(z\) gives:
$
\log z = \log\left(\frac{u_1-u_2}{\sigma}\right) - \frac{1}{2}\log(\beta_1^2+\beta_2^2).
$
By factoring out \(\min\{\kappa_1+1,\kappa_2+1\}^{-2\alpha}\) from \(\beta_1^2+\beta_2^2=(\kappa_1+1)^{-2\alpha}+(\kappa_2+1)^{-2\alpha}\):
\[
{\color{black}
\begin{aligned}
\log(\beta_1^2+\beta_2^2)
=\, & -2\alpha\log\left(\min\{\kappa_1+1,\kappa_2+1\}\right) +\log\left[
\left(\tfrac{\min\{\kappa_1+1,\kappa_2+1\}}{\kappa_1+1}\right)^{2\alpha}
+ \left(\tfrac{\min\{\kappa_1+1,\kappa_2+1\}}{\kappa_2+1}\right)^{2\alpha}
\right].
\end{aligned}
}
\]
Combining these results yields the desired expression.
\end{proofref}

The preceding theorem gives an exact formula and refined expansion for the two-arm case. We next extend the conclusion to a multi-arm setting under QCARE with general variance.

\begin{proposition}[Multi-Arm Reduced Form]
\label{prop:reduced-form-fixed-N}
Consider states \(S=(\kappa,u)\) with \(\Delta_j(S):=u_1-u_j >0\) for all $ j > 1 $, 
and \(j^*=j^*(S)\in\arg\max_{j\ge2}p_j(S)\). Along the sequence where \(\kappa \to\infty\) and \(\{\Delta_j(S)\}\) are bounded away from zero, it holds that
\begin{equation*}
    \frac{1}{2}\log\log\left(\frac{1}{Q_{\mathrm{non\text{-}lead}}(S)}\right)
    =
    \alpha\log\left(\min\{\kappa_1+1,\kappa_{j^*}+1\}\right)+O(1).
\end{equation*}
\end{proposition}

\begin{proofref}{Proposition~\ref{prop:reduced-form-fixed-N}}
	For each \(j\geq 2\), based on Equation~\eqref{eq: pi_property}, 
	$
	p_j(S)
	=
	1-\Phi\left(
	\frac{\Delta_j(S)}
	{\sigma\sqrt{
			(\kappa_1+1)^{-2\alpha}
			+
			(\kappa_j+1)^{-2\alpha}
	}}
	\right).
	$ 
	Because arm \(1\) is the unique leading arm,
	$
	\left\{
	\theta^{\sigma}_1(S)
	<
	\max_{j\geq 2}\theta^{\sigma}_j(S)
	\right\}
	=
	\bigcup_{j=2}^N
	\left\{
	\theta^{\sigma}_j(S)
	>
	\theta^{\sigma}_1(S)
	\right\}.
	$
	Due to the union bound, it follows that
	\[
	\max_{j\geq 2}p_j(S)
	\leq
	Q_{\mathrm{non\text{-}lead}}(S)
	\leq
	\sum_{j=2}^N p_j(S)
	\leq
	(N-1)\max_{j\geq 2}p_j(S).
	\]

	Finally, consider any sequence of states such that
	\(\min_i\{\kappa_i\}\to\infty\). Since
	\(\Delta_j(S)\geq\underline{\Delta}>0\) for every \(j\geq2\), the
	formula for \(p_j(S)\) above implies that \(p_j(S)\to0\) for every
	\(j\geq2\). Because \(N\) is fixed, this also implies
	\(p_{j^*}(S)\to0\). The preceding sandwich bound gives
	\[
	\log\left(\frac{1}{p_{j^*}(S)}\right)-\log(N-1)
	\leq
	\log\left(\frac{1}{Q_{\mathrm{non\text{-}lead}}(S)}\right)
	\leq
	\log\left(\frac{1}{p_{j^*}(S)}\right).
	\]
	Therefore,
	\[
	\begin{aligned}
	0
	&\leq
	\log\log\left(\frac{1}{p_{j^*}(S)}\right)
	-
	\log\log\left(\frac{1}{Q_{\mathrm{non\text{-}lead}}(S)}\right)\leq
	-\log\left(
	1-
	\frac{\log(N-1)}
	{\log(1/p_{j^*}(S))}
	\right)
	=o(1),
	\end{aligned}
	\]
	and hence
	$
	\log\log\left(\frac{1}{Q_{\mathrm{non\text{-}lead}}(S)}\right)
	=
	\log\log\left(\frac{1}{p_{j^*}(S)}\right)+o(1).
	$ Note that \(j^*(S)\) is valued in $ [N] $. Along any infinite subsequence on which \(j^*(S)=j\) for some fixed \(j\), the quantity
	\(p_{j^*}(S)=p_j(S)\) is the non-leading-arm choice probability in the induced two-arm comparison between arms \(1\) and
	\(j\). Moreover,
	$
	\Delta_j(S)\geq\underline{\Delta}
	$ and $
	\min\{\kappa_1+1,\kappa_j+1\}\to\infty.
	$
	Theorem~\ref{thm:reduced-form-two-arm} therefore gives, along this
	subsequence,
	\[
	\frac{1}{2}\log\log\left(\frac{1}{p_{j^*}(S)}\right)
	=
	\alpha\log\left(
	\min\{\kappa_1+1,\kappa_{j^*}+1\}
	\right)+O(1).
	\]
	Because there are only finitely many possible values of \(j^*(S)\), the preceding bound holds along the
	original sequence of states.
\end{proofref}

\end{appendices}
\end{document}